%% file: opus.tex
\documentclass{article}

\newtheorem{fed}{\textbf{Definition}}[section]
\newtheorem{thm}[fed]{\textbf{Theorem}}
\newtheorem{lemma}[fed]{\textbf{Lemma}}
\newtheorem{ex}[fed]{\textbf{Example}}
\newtheorem{rem}[fed]{\textbf{Remark}}
\newtheorem{prop}[fed]{\textbf{Proposition}}
\newtheorem{cor}[fed]{\textbf{Corollary}}
\usepackage{amssymb,bbm,graphicx,epsfig,psfrag,epic,eepic,latexsym}
\usepackage{amsmath}
\usepackage{mathrsfs} 
\usepackage[dvips]{color}       

\begin{document}
\title{The Arnold-Givental conjecture and moment Floer homology}
\author{Urs Frauenfelder\footnote{Partially supported by Swiss National Science
Foundation}}
\maketitle

\begin{abstract}
In this article we prove the Arnold-Givental conjecture for a class of 
Lagrangian submanifolds in Marsden-Weinstein quotients which are
fixpoint sets of some antisymplectic involution. For these 
Lagrangians the Floer homology cannot in general be defined by standard
means due to the bubbling phenomenon. To overcome this difficulty we 
consider moment Floer homology whose boundary operator is defined by
counting solutions of the symplectic vortex equations on the strip which
have better compactness properties than the original Floer equations. 
\end{abstract}
\tableofcontents

\newpage

\section[Introduction]{Introduction}

Assume that $(M,\omega)$ is a $2n$-dimensional
compact symplectic manifold, 
$L \subset M$ is a compact Lagrangian submanifold of $M$, and
$R \in \mathrm{Diff}(M)$ is an antisymplectic involution, i.e.
$$R^*\omega=-\omega, \quad R^2=\mathrm{id}$$
whose fixpoint set is the Lagrangian
$$L=\mathrm{Fix}(R).$$
The Arnold-Givental conjecture gives a lower bound on the number of
intersection points of $L$ with a Hamiltonian isotopic Lagrangian
submanifold which intersects $L$ transversally in terms of the Betti
numbers of $L$, more precisely
\\ \\
\textbf{Conjecture (Arnold-Givental)}
\emph{For $t \in [0,1]$ let $H_t \in C^\infty(M)$
be a smooth family of Hamiltonian function of $M$ and denote by
$\phi_H$ the one-time map of the flow of the Hamiltonian vector field
$X_{H_t}$
of $H_t$. Assume that $L$ and $\phi_H(L)$ intersect transversally. Then
the number of intersection points of $L$ and $\phi_H(L)$ can be
estimated from below by the sum of the $\mathbb{Z}_2$ Betti numbers of
$L$, i.e.
$$\#(L \cap \phi_H(L)) \geq \sum_{k=0}^n b_k(L;\mathbb{Z}_2).$$
}
\\ \\ 
Using the fact that the diagonal $\Delta$ is a Lagrangian submanifold
of $(M \times M, \omega \oplus -\omega)$ which equals the fixpoint
set of the antisymplectic involution of $M \times M$ given by
interchanging the two factors, the Arnold conjecture with
$\mathbb{Z}_2$ coefficients for arbitrary
compact symplectic manifolds is a Corollary of the Arnold-Givental
conjecture.
\\ \\
\textbf{Corollary (Arnold conjecture) } \emph{Assume that for $t \in
[0,1]$ there is a smooth family 
$H_t \in C^\infty(M)$ of Hamiltonian functions such
that $1$ is not an eigenvalue of $d \phi_H(x)$ for each fixpoint
$x$ of $\phi_H$, then the number of fixpoints of $\phi_H$ can be
estimated from below by the sum of the $\mathbb{Z}_2$ Betti numbers of
$M$, i.e.
$$\#\mathrm{Fix}(\phi_H) \geq \sum_{k=0}^{2n}b_k(M,\mathbb{Z}_2).$$ 
}
\\
Up to now, the Arnold-Givental conjecture could only be proven 
under some additional assumptions. Givental proved it in 
\cite{givental} for $(\mathbb{C}P^n,\mathbb{R}P^n)$, Oh proved it
in \cite{oh3} for real forms of compact Hermitian spaces with suitable
conditions on the Maslov indices, Lazzarini proved it in
\cite{lazzarini} in the negative monotone case under suitable
assumptions on the minimal Maslov number,
and Fukaya, Oh, Ohta, Ono proved it in
\cite{fukaya-oh-ohta-ono} for semipositive Lagrangian submanifolds. In
my thesis, see \cite{frauenfelder2}, I introduced moment Floer homology
and proved the Arnold-Givental conjecture for a class of Lagrangians in 
Marsden-Weinstein quotients which satisfy some monotonicity
condition. In this paper, this monotonicity condition will be removed.

\subsection[Outline of the paper]{Outline of the paper}

In section~\ref{agi} we give a heuristic argument why the
Arnold-Givental conjecture should hold. This argument is based on Floer
theory. We will ignore questions of bubbling and of
transversality. The main point is that due to the antisymplectic
involution some contribution to the boundary operator should cancel.

To describe the cancelation process one has to choose the almost complex
structure in a non-generic way. For such a choice of almost complex
structure one cannot in general expect to achieve transversality. 
To overcome this problem we consider in section~\ref{transi} abstract 
perturbations as in \cite{fukaya-ono2} and prove some transversality
result specific for problems of Arnold-Givental type. 
We will need that later to
compute moment Floer homology. In my thesis, see \cite{frauenfelder2},
this technique was not available to me and hence I could only compute
moment Floer homology under some monotonicity assumption.
 
In section~\ref{mofoho} we introduce moment Floer homology. The chain
complex of moment Floer homology is generated by the intersection points
of some Lagrangian submanifold of a Marsden-Weinstein quotient
which is the fixpoint set of an antisymplectic involution with its
image under a generic Hamiltonian isotopy. The boundary operator is
defined by counting solutions of the symplectic vortex equations on the
strip. These equations contain an
equivariant version of Floer's equation and a condition which relates
the curvature to the moment map. We prove under some topological 
restrictions on the envelopping manifold a compactness theorem which
allows us to prove that the boundary operator is well-defined. To 
compute the moment Floer homology we set the Hamiltonian equal to
zero. This is the infinite dimensional analogon of a Morse-Bott
situation. Combining the
techniques developed in section~\ref{transi} together with the 
approach to Morse-Bott theory described in the
appendix by defining a boundary operator by
counting flow lines with cascades we prove that moment Floer homology
is isomorphic to the singular homology of the Lagrangian in the
Marsden-Weinstein quotient with coefficients in some Novikov ring. 

In the appendix we develop some approach to define Morse-Bott homology
by counting flow-lines with cascades. We prove that for 
a given Morse-Bott function and for generic choice of a Riemannian metric on
the manifold and of a Morse-function and a Riemannian
metric on the critical manifold transversality for flow-lines with
cascades can be achieved and hence the boundary operator is
well-defined. We show that Morse-Bott homology is independent of
the Morse-Bott function and hence isomorphic to the ordinary Morse homology.
 
\subsection[Acknowledgements]{Acknowledgements}

I would like to express my deep gratitude to the supervisor of my
thesis, Prof. D. Salamon, for pointing
my attention to moment Floer homology, for his encouragement and a lot of
lively discussions. I cordially thank the two coexaminers Prof. K. Ono and
Prof. E. Zehnder for a great number of valuable suggestions. In
particular, I would like to thank Prof. Ono for enabling me to come to
Hokkaido University and for explaining to me the theory of Kuranishi
structures.

\section[Arnold-Givental conjecture]{Arnold-Givental conjecture}\label{agi}

We will give in this section  a heuristic argument 
based on Floer theory why the Arnold-Givental
conjecture should hold. 
We refer the reader to Floer's original papers
\cite{floer1,floer2, floer3, floer4, floer5}
and to \cite{salamon1} for an introduction into the topic
of Floer theory. 
We will here just introduce the main objects of the theory.
Moreover, we
completely ignore questions of bubbling and of transversality. We will
address the question of transversality in 
section~\ref{transi}. 
We hope that these techniques combined with the techniques 
developed in \cite{fukaya-ono2} and \cite{fukaya-oh-ohta-ono} to
overcome the bubbling phenomenon should lead to a proof of the
Arnold-Givental conjecture in general. However, in this paper we
will not persecute such an approach but consider instead 
in section~\ref{mofoho} moment
Floer homology where no bubbling occurs.  

We abbreviate by $\Gamma$ the group
$$\Gamma=\frac{\pi_2(M,L)}{\mathrm{ker}I_\mu \cap
\mathrm{ker}I_\omega}$$
where $I_\omega \colon \pi_2(M,L) \to \mathbb{R}$ and 
$I_\mu \colon \pi_2(M,L) \to \mathbb{Z}$ are the homomorphisms induced
from the symplectic structure $\omega$ respectively the Maslov number
$\mu$. We refer the reader to \cite{robbin-salamon1, robbin-salamon2,
salamon-zehnder} for a discussion of the Maslov index.
Following \cite{hofer-salamon} 
we denote the Novikov ring $\Lambda=\Lambda_\Gamma$ as the ring 
consisting of formal sums
$$r=\sum_{\gamma \in \Gamma} r_\gamma \gamma, \quad r_\gamma 
\in \mathbb{Z}_2$$
which satisfy the finiteness condition
$$\#\{\gamma \in \Gamma: r_\gamma \neq 0,\,\,I_\omega(\gamma) \geq
\kappa\}
<\infty, \quad \forall \,\,\kappa \in \mathbb{R}.$$
Note that since the coefficients are taken in the field $\mathbb{Z}_2$
the Novikov ring is actually a field. It is naturally graded by
$I_\mu$.

To define the Floer chain complex we consider pairs 
$(\bar{x},x) \in C^\infty((\Omega, \Omega \cap \mathbb{R}),(M,L)) 
\times C^\infty(([0,1],\{0,1\}),(M,L))$
\footnote{We abbreviate for manifolds $A_2 \subset A_1$ and $B_2 \subset B_1$
by $C^\infty((A_1,A_2),(B_1,B_2))$ the space of smooth maps from
$A_1$ to $B_1$ which map $A_2$ to $B_2$.} for the 
half-disk
$\Omega=\{z \in \mathbb{C}: |z| \leq 1,\,\,\mathrm{Im}(z) \geq 0\}$ 
which satisfy the following conditions
$$\dot{x}(t)=X_{H_t}(x(t)), \quad \bar{x}(e^{i\pi t})=x(t),
\quad t \in [0,1].$$
We introduce an equivalence relation on these pairs by
$$(\bar{x},x) \cong (\bar{y},y) \Longleftrightarrow x=y,\,\,
\omega(\bar{x})=\omega(\bar{y}),\,\,\mu(\bar{x})=\mu(\bar{y})$$
and denote the set of equivalence classes by $\mathscr{C}$, recalling
the fact that this set may be interpreted as the critical set of
an action functional on a covering of the space of paths in $M$
which connect two points of the Lagrangian $L$. The Floer chain complex
$CF_*(M,L;H)$ can now be defined as the graded $\mathbb{Z}_2$ vector space
consisting of formal sums
$$\xi=\sum_{c \in \mathscr{C}}\xi_c c, \quad \xi_c \in \mathbb{Z}_2$$
satisfying the finiteness condition
$$\#\{c=[\bar{x},x] \in \mathscr{C}: \xi_c \neq 0,\,\,I_\omega(\bar{x}) \geq
\kappa\}
<\infty, \quad \forall \,\,\kappa \in \mathbb{R}.$$
The grading of $CF_*$ is induced from the Maslov index. The natural
action of $\Gamma$ on $CF_*$ by cocatenation induces on $CF_*$ the
structure of a graded vector space over the graded field $\Lambda$.
\\
To define the boundary operator we choose a smooth family
of $\omega$-compatible almost complex structures $J_t$ for $t \in [0,1]$
and count for two critical points
$[\bar{x},x], [\bar{y},y] \in \mathscr{C}$
solutions $u \in C^\infty([0,1]\times \mathbb{R},M)$
of the following problem
\begin{eqnarray}\label{floer}
\partial_s u+J_t(u)(\partial_t u-X_{H_t}(u)) &=&0\\ \nonumber
u(s,j) \in L, \,\,& &j \in \{0,1\}\\ \nonumber
\lim_{s \to -\infty}u(s,t) &=& x(t)\\ \nonumber
\lim_{s \to \infty}u(s,t) &=& y(t)\\ \nonumber
\bar{x} \#[u]\# \bar{y} &=& 0.
\end{eqnarray}
Here the limites are uniform in the $t$-variable with respect to the
$C^1$-topology and $\#$ denotes cocatenation. For generic choice of
the almost complex structures the moduli space 
$\tilde{\mathcal{M}}([\bar{x},x],[\bar{y},y])$ of solutions of
(\ref{floer}) is a smooth manifold of dimension
$$\mathrm{dim}(\tilde{\mathcal{M}}([\bar{x},x],[\bar{y},y]))
=\mu(\bar{x})-\mu(\bar{y}).$$
If $[\bar{x},x]$ is different from $[\bar{y},y]$ the group 
$\mathbb{R}$ acts freely on the solutions of (\ref{floer}) by timeshift
$$u(s,t) \mapsto u(s+r,t), \quad r \in \mathbb{R}.$$
We denote the quotient by
$$\mathcal{M}([\bar{x},x],[\bar{y},y])=\tilde{\mathcal{M}}([\bar{x},x],
[\bar{y},y])/\mathbb{R}.$$
If we ignore the question of bubbling the manifold
$\mathcal{M}([\bar{x},x],[\bar{y},y])$ for critical points
$[\bar{x},x],[\bar{y},y] \in \mathscr{C}$ satisfying
$\mu(\bar{x})-\mu(\bar{y})=1$ is compact and we may define the
$\mathbb{Z}_2$ numbers
$$n([\bar{x},x],[\bar{y},y]):=\# \mathcal{M}([\bar{x},x],[\bar{y},y])
\,\,\mathrm{mod} \,\,2.$$
The Floer boundary operator is now defined as the linear extension of
$$\partial[\bar{x},x]=\sum_{\substack{[\bar{y},y] \in \mathscr{C},\\
\mu(\bar{y})=\mu(\bar{x})-1}}n([\bar{x},x],[\bar{y},y])[\bar{y},y].$$
Ignoring again the bubbling problem one can ``prove''
$$\partial^2=0$$
and hence one can define
$$HF_*(M,L;H,J)=\frac{\mathrm{ker}\partial}{\mathrm{Im}\partial}.$$
One can show - always ignoring the bubbling problem - that for different
choices of $H$ and $J$ there are canonical isomorphisms between the
graded $\Lambda$ vector spaces and hence Floer homology 
$HF_*(M,L)$ is defined independent of $H$ and $J$. It follows from its
definition that the dimension of $HF_*$ as $\Lambda$ vector space gives
a lower bound on the number of intersection points of $L$ and
$\phi_H(L)$.
\\
To ``compute'' the Floer homology we consider the case where $H=0$.
Actually, in this case $L$ and $\phi_H(L)=L$ do not intersect
transversally but still cleanly
\footnote{Two submanifolds $L_1,L_2 \subset M$ are said to intersect
cleanly if their intersection is a manifold such that for each 
$x \in L_1 \cap L_2$ it holds that
$T_x (L_1 \cap L_2)=T_x L_1 \cap T_x L_2$.}. The finite dimensional analogon of
clean intersections in Floer theory are Morse-Bott functions. In our
case the critical manifold consists of different copies of the
Lagrangian $L$ indexed by the group $\Gamma$. 
Following the approach developed in the
appendix one can still define Floer homology in the
case of clean intersections. To do that one chooses a Morse function on
the critical manifold. The critical points of this Morse function give
a basis for the chain complex. The boundary operator is defined by
counting flow lines with cascades. In our case, a cascade is just a
nonconstant solution $u \in C^\infty([0,1]\times \mathbb{R},M)$
of the unperturbed Floer equation
\begin{eqnarray*}
\partial_s u+J_t(u)\partial_t u&=&0,\\
u(s,j) \in L,& & j \in \{0,1\}
\end{eqnarray*}
which converges uniformly in the $t$-variable to points on $L$
as $s$ goes to $\pm \infty$.
A flow line with zero cascades is just a Morse flow line. A flow line
with one cascade consists of a piece of a Morse flow line a cascade
which converges on the left to the endpoint of this piece and which
converges on the right to the initial point of a second piece of a Morse
flow line, see appendix~\ref{morsebott} for details. The boundary
operator now splits naturally
$$\partial=\partial^0+\partial^1$$
where $\partial^0$ consists of the flow lines with zero cascades and
$\partial^1$ consists of the flow lines with at most one cascade. Note
that $\partial^0$ is precisely the Morse boundary operator. If one can
show that $\partial^1$ is zero it follows that 
$$HF_*(M,L)=HM_*(L;\mathbb{Z}_2)\times_{\mathbb{Z}_2}\Lambda$$
where $HM_*(L;\mathbb{Z}_2)$ is the Morse homology of $L$ with
$\mathbb{Z}_2$ coefficients, which is isomorphic to the singular
homology of $L$ with $\mathbb{Z}_2$ coefficients, see \cite{schwarz1}.
\\
It now remains to explain why one has to expect that $\partial^1$
vanishes. Since the Lagrangian $L$ is the fixpoint set of the
antisymplectic involution $R$ there is an induced involution on the
space of cascades. More precisely, since we ignore questions of transversality
we may now choose the family of $\omega$ compatible
almost complex structures $J_t$
independent of the
$t$-variable and impose furthermore the condition that 
$J$ is
antiinvariant under $R$, i.e.
$$R^*J=-J.$$
When $J$ is independent of the $t$-variable a cascade corresponds
to a nonconstant $J$-holomorphic disk satisfying Lagrangian boundary
conditions, i.e.
\begin{eqnarray}\label{agf}
\partial_s u+J(u)\partial_t u &=&0\\ \nonumber
u(s,j) \in L, \,\,& &j \in \{0,1\}.
\end{eqnarray}
The antisymplectic involution $R$ induces a natural involution on the
solutions of the problem above defined by
$$I_1 u(s,t)=R(u(s,1-t)).$$
This involution does not act freely in general, but it was observed by
K.Ono that on the "lanterns", the fixpoints of the involution $I_1$,
one can define a new involution defined by
$$I_2 u(s,t)=\left\{\begin{array}{cc}
u(s,t+1/2) & 0 \leq t \leq 1/2\\
u(s,t-1/2) & 1/2 \leq t \leq 1
\end{array}\right.$$
Finally, on the ``double lanterns'', the fixpoints of the involution $I_2$, one
can define a third involution given by
$$I_3 u(s,t)=\left\{\begin{array}{cc}
u(s,t+1/4) & 0 \leq t \leq 3/4\\
u(s,t-3/4) & 3/4 \leq t \leq 1
\end{array}\right.$$
and so on. Since we consider only nonconstant solution, there exists
$m \in \mathbb{N}$ such that $I_m$ acts freely. Hence one may expect
that there is an even number of flow lines with at most one
cascade. Since we are working with $\mathbb{Z}_2$ coefficients this
implies that $\partial^1$ is zero.

\section[Transversality]{Transversality}\label{transi}

For general symplectic manifolds one cannot expect to find an
$R$-invariant $\omega$-compatible
almost complex structure which is regular so that the
relevant moduli spaces have the structure of finite dimensional
manifolds. Hence to make precise the heuristic argument at the end of
the last section one has to use abstract perturbation. To do that
one considers solutions of (\ref{agf}) modulo the natural
$\mathbb{R}$ action as the
zero set of a section of an infinite dimensional Banach manifold 
$\mathcal{B}$ into a
Banach bundle $\mathcal{E}$ over it. 
There are natural extensions of the involutions 
$I_k$ to involutions of the Banach manifold, so that again $I_1$ is
defined on the whole space, $I_2$ on the fixpoint set of $I_1$ and so on.
The idea of our abstract perturbation is now to perturb our section 
in such a way that the perturbed section intersects the zero section 
transversally but that its zero set is still invariant under the involutions.  
There are two problems we have to overcome. The first one is that one 
cannot define invariants from the 
zero set of an arbitrary section from in infinite dimensional space and
hence one should define the section in some ``finite dimensional''
neighbourhood of the zero set of the original section. This problem was
solved in \cite{fukaya-ono2} by using Kuranishi structures. The second
problem is that for each involution one should find
extensions $I_k^{T\mathcal{B}}$ 
and $I_k^{\mathcal{E}}$ to the tangent space $T\mathcal{B}$ or the bundle 
$\mathcal{E}$ respectively  
restricted to the domain of definition of $I_k$ in order to
achieve that the zero set of the perturbed section is invariant under
the involutions. 

\subsection[Banach spaces]{Banach spaces}

We interpret solutions of (\ref{agf}) as
the zero-set of a smooth section from a Banach manifold 
$\tilde{\mathcal{B}}$
to a Banachbundle $\tilde{\mathcal{E}}$ over 
$\tilde{\mathcal{B}}$. To define 
$\tilde{\mathcal{B}}$ we first have to introduce some weighted Sobolev norms. 
Choose a smooth cutoff function $\beta \in C^\infty(\mathbb{R})$
such that $\beta(s)=0$ for $s<0$ and $\beta(s)=1$ for $s>1$. Choose
$\delta>0$ and define $\gamma_\delta \in C^\infty(\mathbb{R})$
by
$$\gamma_\delta(s):=e^{\delta \beta(s) s}.$$
Let $\Omega$ be an open subset of 
the strip $\mathbb{R}\times [0,1]$. 
For $1 \leq p \leq \infty$ and $k \in \mathbb{N}_0$ 
we define the $||\,||_{k,p,\delta}$-norm 
for $v \in W^{k,p}(\Omega)$ by
$$||v||_{k,p,\delta}:=\sum_{i+j \leq k}||\gamma_\delta \cdot
\partial^i_s\partial^j_t v||_p.$$
We introduce the following weighted Sobolev spaces
$$W^{k,p}_\delta(\Omega):=\{v \in W^{k,p}(\Omega):||v||_{k,p,\delta}<\infty\}
=\{v \in W^{k,p}(\Omega):\gamma_\delta v \in W^{k,p}(\Omega)\}.$$
We abbreviate
$$L^p_\delta(\Omega):=W^{0,p}_\delta(\Omega).$$
Fix a real number $p>2$ and a Riemannian metric $g$ on $TM$. 
The Banach manifold $\tilde{\mathcal{B}}=
\tilde{\mathcal{B}}^{1,p}_\delta(M,L)$
consists of $W^{1,p}_{loc}$-maps $u$ from the strip 
$\mathbb{R}\times [0,1]$
to $M$ which map the boundary of the strip to the
Lagrangian $L$ and satisfy in addition the following conditions.
\begin{description}
 \item[B1:] There exists a point
  $x^- \in L$,
  a real number $T_1$, and an element
  $v_1 \in W^{1,p}_\delta((-\infty,-T_1)\times [0,1],T_{x^-}M)$
  such that
  $$u(s,t)=\exp_{x^-}(v_1(s,t)), \quad s < -T_1.$$
  Here the exponential map is taken with respect to the metric $g$. 
 \item[B2:] There exists a point $x^+ \in L$, a real number $T_2$, and
  an element
  $v_2 \in W^{1,p}_\delta((T_2,\infty)\times [0,1],T_{x^+}M)$ such that
  $$u(s,t)=\exp_p (v_2(s,t)), \quad s >T_2.$$
\end{description}
We introduce the Banach bundle $\tilde{\mathcal{E}}$ over the Banach-manifold
$\tilde{\mathcal{B}}$ whose fiber over $u \in \tilde{\mathcal{B}}$ is given by
$$\tilde{\mathcal{E}}_u:=L^p_\delta(u^*TM).$$
Define the section
$\tilde{\mathcal{F}} \colon \tilde{\mathcal{B}} \to \tilde{\mathcal{E}}$ by
$$\tilde{\mathcal{F}}(u):=\partial_s u+J(u)\partial_t u$$
for $u \in \tilde{\mathcal{B}}$. Note that the zero set
$\tilde{\mathcal{F}}^{-1}(0)$ consists of solutions of (\ref{agf}).
The vertical differential of $\tilde{\mathcal{F}}$ at $u \in
\mathcal{F}^{-1}(0)$
is given by
$$\tilde{D}_u \xi:=D\tilde{\mathcal{F}}(u)\xi=\partial_s
\xi+J(u)\partial_t
\xi+\nabla_\xi J(u)\partial_t u, \quad \xi \in T_u
\tilde{\mathcal{B}},$$
where $\nabla$ denotes the Levi-Civita connection of the metric
$g(\cdot,\cdot)=\omega(\cdot,J\cdot)$.
There is a natural action of the
group $\mathbb{R}$ on $\tilde{\mathcal{B}}$ and $\tilde{\mathcal{E}}$
given by timeshift
$$u(s,t) \mapsto u(s+r,t), \quad r \in \mathbb{R}.$$
We denote the quotient by
$$\mathcal{B}:=\tilde{\mathcal{B}}/\mathbb{R}, \quad 
\mathcal{E}:=\tilde{\mathcal{E}}/\mathbb{R},$$ 
by 
$$\mathcal{F} \colon \mathcal{B} \to \mathcal{E}$$
we denote the section induced from $\tilde{\mathcal{F}}$ and
by $D_u$ we denote the vertical differential of $\mathcal{F}$ for
$u \in \mathcal{F}^{-1}(0)$. 
By induction we define in the obvious way for every $k \in \mathbb{N}$
smooth involutions $I_k$ on $\mathcal{B}_k=\mathrm{Fix}(I_{k-1})$
where we set $\mathcal{B}_1=\mathcal{B}$. Our aim is to construct
natural extensions of
these involutions to $T\mathcal{B}$ and $\mathcal{E}$
restricted to their domain of definition. We prove the following
theorem.
\begin{thm}\label{invo}
For each $k \in \mathbb{N}$ there exist smooth involutative bundle
maps
$$I^{T\mathcal{B}}_k \colon T\mathcal{B}|_{\mathcal{B}_k} \to
T\mathcal{B}|_{\mathcal{B}_k}, \quad
I^{\mathcal{E}}_k \colon \mathcal{E}|_{\mathcal{B}_k} \to
\mathcal{E}|_{\mathcal{B}_k}$$
which extend $I_k$ and which have the following properties.
\begin{description}
 \item[(i)] The bundle maps commute on their common domain of
	    definition, i.e.
  $$I_k^{T\mathcal{B}} \circ I_{\ell}^{T\mathcal{B}}\xi=
  I_{\ell}^{T\mathcal{B}}\circ I_k^{T\mathcal{B}}\xi, \quad
  I_k^{\mathcal{E}} \circ I_{\ell}^{\mathcal{E}}\eta=
  I_{\ell}^{\mathcal{E}}\circ I_k^{\mathcal{E}}\eta, \quad
  \xi \in T\mathcal{B}|_{\mathcal{B}_{\ell}},\,\,
  \eta \in \mathcal{E}|_{\mathcal{B}_{\ell}}$$
  for $\ell \geq k$.
 \item[(ii)] For $k \in \mathbb{N}$ and
  $u \in \mathcal{F}^{-1}(0) \cap \mathcal{B}_k$ the vertical
	    differential
  of $\mathcal{F}$ commutes with the involutions modulo a compact
	    operator.
  More precisely, there exists a compact operator
  $Q_u \colon T_u \mathcal{B} \to \mathcal{E}_u$ such that
  $$I^{\mathcal{E}}_k \circ D_u-D_u \circ I^{T\mathcal{B}}_k
  =I_k^{\mathcal{E}} \circ Q_u-Q_u \circ I_k^{T\mathcal{B}}.$$
  Moreover, 
  $I_k^{\mathcal{E}} \circ Q_u-Q_u \circ I_k^{T\mathcal{B}}$ 
  vanishes on $T_u\mathcal{B}_k$. 
 \item[(iii)] The restriction of of $I_k^{T\mathcal{B}}$ to the tangent
  space of $\mathcal{B}_k$ equals the differential of $I_k$, i.e.
  $$I_k^{T\mathcal{B}}|_{T\mathcal{B}_k}=d I_k.$$  
\end{description}
\end{thm} 

\subsection[A sequence of involutions]{A sequence of involutions}

The aim of this subsection is to prove Theorem~\ref{invo}. Since the 
involution we want to construct are independent of the $s$-variable
it is most convenient to define them on the space of paths  
whose endpoints lie on the Lagrangian $L$. We define for a real number
$p>1$ the path space 
$$\mathscr{P}:=\mathscr{P}^{1,p}(M,L):=W^{1,p}(([0,1],\{0,1\}),(M,L)).$$
For $x \in \mathscr{P}$ the tangent space of $\mathscr{P}$ at $x$ is
given by
$$T_x\mathscr{P}=W^{1,p}(([0,1],\{0,1\}),(T^*_x M,T^*_x L)).$$
We define the bundle $\mathscr{E}$ over $\mathscr{P}$ by
$$\mathscr{E}=L^p([0,1],T_x^*M).$$
Note that $T \mathscr{P}$ is a subbundle of $\mathscr{E}$. As before we
set
$$\mathscr{P}_1:=\mathscr{P}$$
and define the first involution 
$\mathscr{I}_1 \in \mathrm{Diff}(\mathscr{P}_1)$ by
$$\mathscr{I}_1 x(t):=R(x(1-t)), \quad x \in \mathscr{P}_1.$$
By induction on $k \in \mathbb{N}$ we define for $k \geq 2$
$$\mathscr{P}_k:=\mathrm{Fix}(\mathscr{I}_{k-1}), \quad
\mathscr{I}_k x(t):=x(t+\frac{1}{2^{k-1}}-\lfloor t+\frac{1}{2^{k-1}}
\rfloor), \quad x \in \mathscr{P}_k.$$
Here $\lfloor \,\, \rfloor$ denote the Gauss brackets, i.e. the largest
integer which is smaller then a given real number
$$\lfloor r \rfloor:=\max\{n \in \mathbb{Z}: n \leq r\}, \quad r \in
\mathbb{R}.$$ 
Observe that if the index of integrability $p$ is greater than two and
$u \in \tilde{\mathcal{B}}$ has the property that
$u_s(t):=u(s,t) \in \mathscr{P}_k$ for every $s \in \mathbb{R}$ and
some $k \in \mathbb{N}$, then $\mathscr{I}_k$ induces an involution on
$u$ which commutes with the $\mathbb{R}$ action on $\tilde{\mathscr{B}}$
so that the induced map in the quotient
$\mathcal{B}=\tilde{\mathcal{B}}/\mathbb{R}$ coincides with the
involution
$I_k$.

We will find in this subsection an extension of these
involutions to smooth involutative bundle maps
$\mathscr{I}_k \colon 
\mathscr{E}|_{\mathscr{P}_k} \to \mathscr{E}|_{\mathscr{P}_k}$ 
such that the following properties are satisfied.
\begin{description}
 \item[(i)] The tangent bundle of the path space is invariant under the
  involutions
  $$\mathscr{I}_k(T\mathscr{P}|_{\mathscr{P}_k})=T\mathscr{P}|_{\mathscr{P}_k}.
  $$
  Moreover,
  $$\mathscr{I}_k|_{T\mathscr{P}_k}=d \mathscr{I}_k|_{\mathscr{P}_k}.$$
 \item[(ii)] The involutions commute on their common domain of
	    definition
  $$\mathscr{I}_k \circ \mathscr{I}_\ell \xi=\mathscr{I}_\ell \circ 
  \mathscr{I}_k \xi, \quad \xi \in \mathscr{E}|_{\mathscr{P}_\ell}$$
  for $\ell \geq k$. 
\end{description}
We will then define the involution $I^{T\mathcal{B}}_k$ and
$I^{\mathcal{E}}_k$ of Theorem~\ref{invo} as the induced maps of
$\mathscr{I}_k$ on $T\mathcal{B}|_{\mathcal{B}_k}$ respectively
$\mathcal{E}|_{\mathcal{B}_k}$. Property (i) of the involutions
$\mathscr{I}_k$ guarantees that $I^{T\mathcal{B}}_k$ is well
defined and guarantees assertion (iii) of Theorem~\ref{invo}. 
Assertion (i) of Theorem~\ref{invo} follows from property
(ii) of the involutions $\mathscr{I}_k$ and assertion (ii) of
Theorem~\ref{invo} will follow from our construction.

For $x \in \mathscr{P}$ the 
extension of the first involution to $\mathscr{E}$ is defined in the
following way
$$\mathscr{I}_1 \xi(t):=R^* \xi(1-t) \in \mathscr{E}_{\mathscr{I}_1 x}, 
\quad \xi \in \mathscr{E}_x. $$
If $x \in \mathscr{P}_2$, then $\mathscr{I}_1$ is a bounded
linear involutative map of the Banach space $\mathscr{E}_x$ and leads to a 
decomposition 
$$\mathscr{E}_x=\mathscr{E}_{x,-1} \oplus \mathscr{E}_{x,1},$$
where $\mathscr{E}_{x,-1}$ is the eigenspace of 
$\mathscr{I}_1|_{\mathscr{E}_x}$ to the eigenvalue $-1$ and 
$\mathscr{E}_{x,1}$ is the eigenspace to the eigenvalue $1$. The two
projections to the eigenspaces are given by
$$\pi_{x,1}=
\frac{1}{2}\bigg(\mathrm{id}|_{\mathscr{E}_x}+\mathscr{I}_1|_{\mathscr{E}_x}
\bigg)
\colon \mathscr{E}_x \to \mathscr{E}_{x,1}, \quad
\pi_{x,-1}=
\frac{1}{2}\bigg(\mathrm{id}|_{\mathscr{E}_x}-\mathscr{I}_1|_{\mathscr{E}_x}
\bigg)
\colon \mathscr{E}_x \to \mathscr{E}_{x,-1}.$$
We first define the involutions on the subspace $\mathscr{E}_{x,-1}$. 
For $\xi \in \mathscr{E}_{x,-1}$ the second involution $\mathscr{I}_2$
is defined by the formula
$$\mathscr{I}_2 \xi(t)=(-1)^{\lfloor t+1/2 \rfloor}
J(\mathscr{I}_2 x(t))\xi(t+1/2-\lfloor t+1/2 \rfloor))
\in \mathscr{E}_{\mathscr{I}_2 x}.$$
To see that $\mathscr{I}_2 \xi$ satisfies condition (i), i.e. maps the
tangent space of the path space to itself, we have to check that if
$\xi$ is in $W^{1,p}$ and satisfies the Lagrangian boundary conditions,
then $\mathscr{I}_2 \xi$ satisfies these conditions, too. Since
$\xi$ lies in the eigenspace of the eigenvalue $-1$ of 
the first involution $\mathscr{I}_1$, it follows that 
$$\xi(1/2)=-dR(x(1/2))\xi(1/2)$$
and using anticommutativity of the almost complex structure $J$ with the
antisymplectic involution $R$ together with the fact that the Lagrangian
$L$ equals the fixpoint set of $R$ it follows that
$$J(x(1/2))\xi(1/2) \in T_{x(1/2)}L$$
and hence $\mathscr{I}_2 \xi$ satisfies the required boundary
condition. To prove that $\mathscr{I}_2 \xi \in W^{1,p}$ one has to
check continuity at the point $t=1/2$. This is done similarly as above
by noting that
$$\xi(0)=-\xi(1) \in T_{x(0)}L=T_{x(1)}L.$$
Moreover, a straightforward calculation shows that 
$$\mathscr{I}_1(\mathscr{I}_2 \xi)=-\mathscr{I}_2 \xi$$
and hence
$$\mathscr{I}_2 \xi \in \mathscr{E}_{\mathscr{I}_2 x, -1}.$$ 
Now assume that $x \in \mathscr{P}_m$ for some integer $m>2$. To
define the involutions $\mathscr{I}_3, \ldots, \mathscr{I}_m$ on
$\mathscr{E}_{x,-1}$ we first consider the following linear maps 
$$\mathscr{L}_k \colon \mathscr{E}_{x,-1} \to \mathscr{E}_
{\mathscr{I}_{k+2}x,-1}, \quad k \in \{1,\ldots,m-2\}$$
defined by
$$\mathscr{L}_k \xi(t)=\sum_{j=0}^{2^k-1}
(-1)^{\lfloor \frac{j}{2^{k-1}}\rfloor+\lfloor t+\frac{1}{2^{k+1}}+
\frac{j}{2^k}\rfloor}\xi\bigg(t+\frac{1}{2^{k+1}}+\frac{j}{2^k}
-\lfloor t+\frac{1}{2^{k+1}}+\frac{j}{2^k}\rfloor\bigg).$$
\begin{lemma} The maps $\mathscr{L}_k$ for $1 \leq k \leq m-2$ are
well-defined, i.e. $\mathscr{I}_1(\mathscr{L}_k \xi)=-\mathscr{L}_k \xi$.
They have the following properties.
\begin{description}
 \item[(i)] They leave the tangent bundle of the
  path space invariant, i.e. for $T_{x,-1}\mathscr{P}:=T_x \mathscr{P}
  \cap \mathscr{E}_{x,-1}$ we have
  $$\mathscr{L}_k(T_{x,-1}\mathscr{P})=T_{\mathscr{I}_{k+2}x,-1}\mathscr{P}.$$
 \item[(ii)] They commute with each other and
  with $\mathscr{I}_2|_{\mathscr{E}_{x,-1}}$.
 \item[(iii)] Setting $\mathscr{L}_0:=\mathrm{id}$
  their squares can be determined recursively from the formula
  \begin{equation}\label{square}  
  \mathscr{L}_{k+1}^2=2\Bigg(\mathscr{L}_k^2+\mathscr{L}_k
  \bigg(\sum_{i=0}^{k-1}
  \mathscr{L}_i\bigg)\Bigg).
  \end{equation}
 \item[(iv)] If $p=2$ then $\mathscr{L}_k$ is selfadjoint with respect
  to the $L^2$-inner product 
  $$\langle \xi, \eta \rangle=\int_0^1 \xi(t)\eta(t)dt.$$
 \item[(v)] The maps $\mathscr{L}_k$ are injective. 
 \end{description}
\end{lemma}
\textbf{Proof: }To prove that the maps $\mathscr{L}_k$ are well defined
we calculate
\begin{eqnarray*}
& &R^*\mathscr{L}_k \xi(1-t)\\
&=&\sum_{j=0}^{2^k-1}(-1)^
{\lfloor \frac{j}{2^{k-1}}\rfloor+\lfloor 1-t+\frac{1}{2^{k+1}}+\frac{j}{2^k}
\rfloor}\\
& &\quad \cdot R^*\xi\bigg(1-t+\frac{1}{2^{k+1}}+\frac{j}{2^k}-
\lfloor 1-t+\frac{1}{2^{k+1}}+\frac{j}{2^k}\rfloor\bigg)\\
&=&-\sum_{j=0}^{2^k-1}(-1)^
{\lfloor \frac{j}{2^{k-1}}\rfloor+\lfloor 1-t+\frac{1}{2^{k+1}}+\frac{j}{2^k}
\rfloor}\\
& &\quad \cdot \xi\bigg(t-\frac{1}{2^{k+1}}-\frac{j}{2^k}
+\lfloor 1-t+\frac{1}{2^{k+1}}+\frac{j}{2^k}\rfloor\bigg)\\
&=&-\sum_{i=0}^{2^k-1}(-1)^
{\lfloor 2-\frac{i}{2^{k-1}}-\frac{1}{2^{k-1}}\rfloor
+\lfloor 1-t+\frac{1}{2^{k+1}}+1-\frac{i}{2^k}-\frac{1}{2^k}\rfloor}\\
& &\quad \cdot \xi\bigg(t-\frac{1}{2^{k+1}}-1+\frac{i}{2^k}+\frac{1}{2^k}+
\lfloor 1-t+\frac{1}{2^{k+1}}+1-\frac{i}{2^k}-\frac{1}{2^k}\rfloor\bigg)\\
&=&-\sum_{i=0}^{2^k-1}(-1)^{\lfloor -\frac{i}{2^{k-1}}-\frac{1}{2^{k-1}}
\rfloor+\lfloor -t-\frac{1}{2^{k+1}}-\frac{i}{2^k}\rfloor}\\
& &\quad \cdot \xi\bigg(t+\frac{1}{2^{k+1}}+\frac{i}{2^k}+1+
\lfloor -t-\frac{1}{2^{k+1}}-\frac{i}{2^k}\rfloor\bigg)\\
&=&-\sum_{i=0}^{2^k-1}(-1)^{-\lfloor \frac{i}{2^{k-1}}-\frac{1}{2^{k-1}}
\rfloor-1-\lfloor t+\frac{1}{2^{k+1}}+\frac{i}{2^k}\rfloor-1}\\
& &\quad \cdot \xi\bigg(t+\frac{1}{2^{k+1}}+\frac{i}{2^k}-
\lfloor t+\frac{1}{2^{k+1}}+\frac{i}{2^k}\rfloor\bigg)\\
&=&-\mathscr{L}_k \xi(t).
\end{eqnarray*}
The second last equality only holds in the case, when
$t \neq \frac{j}{2^k}+\frac{1}{2^{k+1}}$ for $j \in \{0,\ldots,
2^k-1\}$. However, the equality above still holds for $\xi$ in the 
$L^p$-sense.

To prove assertion (i), one has to show that if $\xi \in
T_{x,-1}\mathscr{P}$, then $\mathscr{L}_k \xi$ satisfies the Lagrangian
boundary condition and is continuous at the points 
$t=\frac{1}{2^{k+1}}+\frac{j}{2^k}$ for $j \in \{0, \ldots 2^k-1\}$.
To prove the boundary conditions one checks that 
$$dR \mathscr{L}_k \xi(0)=\mathscr{L}_k \xi(0), \quad
dR \mathscr{L}_k \xi(1)=\mathscr{L}_k \xi(1).$$
To prove continuity one uses 
$$\xi(0)=-\xi(1)$$
which follows from the assumption that $\xi \in
T_{x,-1}\mathscr{P}$.

To prove assertion (ii) one checks using the formula
$$\lfloor r-\lfloor s \rfloor \rfloor=\lfloor r \rfloor -\lfloor s
\rfloor$$
that
\begin{eqnarray*}
& &\mathscr{L}_k \mathscr{L}_\ell \xi(t)\\
&=&\sum_{j=0}^{2^k-1}\sum_{i=0}^{2^\ell-1}
(-1)^{\lfloor \frac{j}{2^{k-1}}\rfloor+\lfloor \frac{i}{2^{\ell-1}}\rfloor
+\lfloor t+\frac{1}{2^{k+1}}+\frac{1}{2^{\ell+1}}+
\frac{j}{2^k}+\frac{i}{2^\ell}\rfloor}\\
& &\quad \cdot \xi\bigg(t+\frac{1}{2^{k+1}}+\frac{1}{2^{\ell+1}}+
\frac{j}{2^k}+\frac{i}{2^\ell}-
\lfloor t+\frac{1}{2^{k+1}}+\frac{1}{2^{\ell+1}}+
\frac{j}{2^k}+\frac{i}{2^\ell}\rfloor\bigg).
\end{eqnarray*}
This formula is symmetric in $k$ and $\ell$ and hence commutativity
follows. In a similar way one proves commutativity with 
$\mathscr{I}_2|_{\mathscr{E}_{x,-1}}$.

It is straightforward to check that (\ref{square}) holds
if $k=0$. We now assume $k \geq 1$. Using the formula
$$\sum_{i=0}^{k-1}\mathscr{L}_i \xi(t)=\sum^{2^k-1}_{\substack{\ell=0\\
\ell \neq 2^{k-1}}}(-1)^{\lfloor \frac{\ell}{2^{k-1}}\rfloor+
\lfloor t+\frac{\ell}{2^k}\rfloor}\xi\bigg(t+\frac{\ell}{2^k}-
\lfloor t+\frac{\ell}{2^k}\rfloor\bigg)$$
we deduce that
\begin{eqnarray*}
\mathscr{L}_k\bigg( \sum_{i=0}^{k-1}\mathscr{L}_i\bigg)\xi(t) &=&
\sum_{\substack{\ell=0\\
\ell \neq 2^{k-1}}}^{2^k-1}\sum_{j=0}^{2^k-1}
(-1)^{\lfloor \frac{\ell}{2^{k-1}}\rfloor+\lfloor
\frac{j}{2^{k-1}}\rfloor
+\lfloor t+\frac{\ell}{2^k}+\frac{j}{2^k}+\frac{1}{2^{k+1}}\rfloor}\\
& &\quad\cdot \xi\bigg(t+\frac{\ell}{2^k}+\frac{j}{2^k}+\frac{1}{2^{k+1}}
-\lfloor t+\frac{\ell}{2^k}+\frac{j}{2^k}+\frac{1}{2^{k+1}}\rfloor\bigg).
\end{eqnarray*}
We calculate
\begin{eqnarray*}
& &\Bigg(\mathscr{L}_{k+1}^2-2\mathscr{L}_k^2-2\mathscr{L}_k
\bigg(\sum_{i=0}^{k-1}\mathscr{L}_i \bigg)\Bigg)\xi(t)\\
&=&\sum_{j=0}^{2^{k+1}-1}\sum_{i=0}^{2^{k+1}-1}
(-1)^{\lfloor \frac{j}{2^k}\rfloor+\lfloor \frac{i}{2^k}\rfloor+
\lfloor t+\frac{1+j+i}{2^{k+1}}\rfloor}\\
& &\quad \cdot
\xi\bigg(t+\frac{1+j+i}{2^{k+1}}-\lfloor t+\frac{1+j+i}{2^{k+1}}\rfloor \bigg)
\\
& &-2\sum_{j=0}^{2^k-1}\sum_{i=0}^{2^k-1}
(-1)^{\lfloor \frac{j}{2^{k-1}}\rfloor+\lfloor \frac{i}{2^{k-1}}\rfloor+
\lfloor t+\frac{1+j+i}{2^k}\rfloor}\\
& &\quad \cdot
\xi\bigg(t+\frac{1+j+i}{2^k}-\lfloor t+\frac{1+j+i}{2^k}\rfloor \bigg)\\
& &-2\sum_{\substack{j=0\\
j \neq 2^{k-1}}}^{2^k-1}\sum_{i=0}^{2^k-1}
(-1)^{\lfloor \frac{j}{2^{k-1}}\rfloor+\lfloor \frac{i}{2^{k-1}}\rfloor+
\lfloor t+\frac{1+2j+2i}{2^{k+1}}\rfloor}\\
& &\quad \cdot
\xi\bigg(t+\frac{1+2j+2i}{2^{k+1}}-\lfloor t+\frac{1+2j+2i}{2^{k+1}}
\rfloor \bigg)
\end{eqnarray*}
\begin{eqnarray*}
&=&\sum_{j=0}^{2^{k+1}-1}\sum_{i=0}^{2^{k+1}-1}
(-1)^{\lfloor \frac{j}{2^k}\rfloor+\lfloor \frac{i}{2^k}\rfloor+
\lfloor t+\frac{1+j+i}{2^{k+1}}\rfloor}\\
& &\quad \cdot
\xi\bigg(t+\frac{1+j+i}{2^{k+1}}-\lfloor t+\frac{1+j+i}{2^{k+1}}\rfloor \bigg)
\\
& &-\sum_{j=0}^{2^k-1}\sum_{i=0}^{2^k-1}
(-1)^{\lfloor \frac{2j+1}{2^k}\rfloor+\lfloor \frac{2i}{2^k}\rfloor+
\lfloor t+\frac{1+2j+1+2i}{2^{k+1}}\rfloor}\\
& &\quad \cdot
\xi\bigg(t+\frac{1+2j+1+2i}{2^{k+1}}-
\lfloor t+\frac{1+2j+1+2i}{2^{k+1}}\rfloor \bigg)\\
& &-\sum_{j=0}^{2^k-1}\sum_{i=0}^{2^k-1}
(-1)^{\lfloor \frac{2j}{2^k}\rfloor+\lfloor \frac{2i+1}{2^k}\rfloor+
\lfloor t+\frac{1+2j+2i+1}{2^{k+1}}\rfloor}\\
& &\quad \cdot
\xi\bigg(t+\frac{1+2j+2i+1}{2^{k+1}}-
\lfloor t+\frac{1+2j+2i+1}{2^{k+1}}\rfloor \bigg)\\
& &-\sum_{j=0}^{2^k-1}\sum_{i=0}^{2^k-1}
(-1)^{\lfloor \frac{2j}{2^k}\rfloor+\lfloor \frac{2i}{2^k}\rfloor+
\lfloor t+\frac{1+2j+2i}{2^{k+1}}\rfloor}\\
& &\quad \cdot
\xi\bigg(t+\frac{1+2j+2i}{2^{k+1}}-\lfloor t+\frac{1+2j+2i}{2^{k+1}}
\rfloor \bigg)\\
& &-\sum_{j=1}^{2^k}\sum_{i=0}^{2^k-1}
(-1)^{\lfloor \frac{2j-1}{2^k}\rfloor+\lfloor \frac{2i+1}{2^k}\rfloor+
\lfloor t+\frac{1+2j-1+2i+1}{2^{k+1}}\rfloor}\\
& &\quad \cdot
\xi\bigg(t+\frac{1+2j-1+2i+1}{2^{k+1}}-\lfloor t+\frac{1+2j-1+2i+1}{2^{k+1}}
\rfloor \bigg)\\
&=&0
\end{eqnarray*}
This proves (\ref{square}).

To prove selfadjointness we calculate
\begin{eqnarray*}
& &\int_0^1 \mathscr{L}_k \xi(t) \eta(t)dt\\
&=&\sum_{j=0}^{2^k-1}\int_0^1
(-1)^{\lfloor \frac{j}{2^{k-1}}\rfloor+
\lfloor t+\frac{1}{2^{k+1}}+\frac{j}{2^k}\rfloor}\\
& &\quad \cdot
\xi\bigg(t+\frac{1}{2^{k+1}}+\frac{j}{2^k}-\lfloor t+\frac{1}{2^{k+1}}+
\frac{j}{2^k}\rfloor\bigg)\eta(t)dt\\
&=&\sum_{j=0}^{2^k-1}\int_0^1
(-1)^{\lfloor \frac{j}{2^{k-1}}\rfloor+
\lfloor s-\frac{1}{2^{k+1}}-\frac{j}{2^k}-\lfloor s-\frac{1}{2^{k+1}}-
\frac{j}{2^k}\rfloor+\frac{1}{2^{k+1}}+\frac{j}{2^k}\rfloor }\\
& &\quad \cdot
\xi(s)\eta\bigg(s-\frac{1}{2^{k+1}}-\frac{j}{2^k}-
\rfloor s-\frac{1}{2^{k+1}}-\frac{j}{2^k}\rfloor\bigg)ds\\
&=&\sum_{j=0}^{2^k-1}\int_0^1
(-1)^{\lfloor \frac{j}{2^{k-1}}\rfloor+
\lfloor s-\lfloor s-\frac{1}{2^{k+1}}-
\frac{j}{2^k}\rfloor \rfloor }\\
& &\quad \cdot
\xi(s)\eta\bigg(s-\frac{1}{2^{k+1}}-\frac{j}{2^k}-
\rfloor s-\frac{1}{2^{k+1}}-\frac{j}{2^k}\rfloor\bigg)ds\\
&=&\sum_{i=0}^{2^k-1}\int_0^1
(-1)^{\lfloor 2-\frac{1}{2^{k-1}}-\frac{i}{2^{k-1}}\rfloor
+\lfloor s-\lfloor s-\frac{1}{2^{k+1}}-1+\frac{1}{2^k}+\frac{i}{2^k}\rfloor
\rfloor}\\
& &\quad \cdot \xi(s)\eta\bigg(s-\frac{1}{2^{k+1}}-1+\frac{1}{2^k}
+\frac{i}{2^k}-\lfloor s-\frac{1}{2^{k+1}}-1+\frac{1}{2^k}
+\frac{i}{2^k}\rfloor \bigg)\\
&=&\sum_{i=0}^{2^k-1}
\int_0^1 (-1)^{\lfloor-\frac{1}{2^{k-1}}-\frac{i}{2^{k-1}}\rfloor
+\lfloor s \rfloor -\lfloor s+\frac{1}{2^{k+1}}+\frac{i}{2^k}\rfloor+1}\\
& &\quad \cdot \xi(s)\eta\bigg(s+\frac{1}{2^{k+1}}+\frac{i}{2^k}
-\lfloor s+\frac{1}{2^{k+1}}+\frac{i}{2^k}\rfloor\bigg)\\
&=&\sum_{i=0}^{2^k-1}
\int_0^1 (-1)^{-\lfloor \frac{1}{2^{k-1}}+\frac{i}{2^{k-1}}\rfloor
-\lfloor s+\frac{1}{2^{k+1}}+\frac{i}{2^k}\rfloor}\\
& &\quad \cdot \xi(s)\eta\bigg(s+\frac{1}{2^{k+1}}+\frac{i}{2^k}
-\lfloor s+\frac{1}{2^{k+1}}+\frac{i}{2^k}\rfloor\bigg)\\
&=&\int_0^1 \xi(s)\mathscr{L}_k \eta(s)ds
\end{eqnarray*}

To prove injectivity we observe that using the formula
\begin{eqnarray*}
\mathscr{L}_k \xi(t+\frac{i}{2^k})&=&
\sum_{\ell=0}^{2^k-1}(-1)^{\lfloor \frac{\ell-i}{2^{k-1}}\rfloor
+\lfloor \frac{\ell-i}{2^k}\rfloor+\lfloor t+\frac{1}{2^{k+1}}+
\frac{\ell}{2^k}\rfloor}\\
& &\quad \cdot \xi\bigg(t+\frac{1}{2^{k+1}}+\frac{\ell}{2^k}
-\lfloor t+\frac{1}{2^{k+1}}+\frac{\ell}{2^k}\rfloor\bigg)
\end{eqnarray*}
for $i \in \{0,\ldots,2^k-1\}$ and $0 \leq t < 1/2^k$ it suffices to
show that the matrix $A(t)\in \mathbb{R}^{2^k \times 2^k}$ whose entries
are given by
$$A_{i,\ell}(t)=(-1)^{\lfloor \frac{\ell-i}{2^{k-1}}\rfloor
+\lfloor \frac{\ell-i}{2^k}\rfloor+\lfloor t+\frac{1}{2^{k+1}}+
\frac{\ell}{2^k}\rfloor}$$
is nondegenerate. Using the elementary fact that the determinant of a
matrix is invariant if one subtracts from a row a different row together
with the formulas
$$|A_{i+1,\ell}(t)-A_{i,\ell}(t)|=2 \delta_{i+2^{k-1},\ell},
\quad i \in \{0,\ldots, 2^{k-1}-1\}$$
$$|A_{i+1,\ell}(t)-A_{i,\ell}(t)|=2 \delta_{i-2^{k-1},\ell},
\quad i \in \{2^{k-1}-1, \ldots, 2^k-2\}$$
$$|A_{0,\ell}(t)-A_{2^k-1,\ell}(t)|=2\delta_{2^{k-1}-1,\ell}$$
one concludes that 
$$|\mathrm{det}(A)|=2^{(2^k)}.$$
This proves injectivity and hence the lemma follows. \hfill $\square$
\\ \\
It follows from (\ref{square}) and from (ii) in the preceding lemma
that for each $k \in \{1,\ldots,m-2\}$
the spectrum of $\mathscr{L}_k^2 \colon \mathscr{E}_{x,-1} \to
\mathscr{E}_{x,-1}$ consists of a finite number of eigenvalues. These 
eigenvalues do not depend on $x$ and can be computed
recursively from (\ref{square}). For $\mathscr{L}_1^2$ the unique
eigenvalue is given by
$$\lambda_1^1=2$$
for $\mathscr{L}_2^2$ the two eigenvalues are given by
$$\lambda_1^2=4+2 \cdot \sqrt{2}, \quad \lambda_2^2=4-2 \cdot \sqrt{2}$$
and for $\mathscr{L}_3^2$ the four eigenvalues are given by
\begin{eqnarray*}
\lambda_1^3&=&2\Bigg
(4+2\cdot \sqrt{2}+\sqrt{4+2 \cdot \sqrt{2}}\bigg(1+\sqrt{2}\bigg)\Bigg)\\
\lambda_2^3&=&2\Bigg
(4-2\cdot \sqrt{2}+\sqrt{4-2 \cdot \sqrt{2}}\bigg(1-\sqrt{2}\bigg)\Bigg)\\
\lambda_3^3&=&2\Bigg
(4+2\cdot \sqrt{2}-\sqrt{4+2 \cdot \sqrt{2}}\bigg(1+\sqrt{2}\bigg)\Bigg)\\
\lambda_4^3&=&2\Bigg
(4-2\cdot \sqrt{2}-\sqrt{4-2 \cdot \sqrt{2}}\bigg(1-\sqrt{2}\bigg)\Bigg).
\end{eqnarray*}
We claim that all the eigenvalues of all the maps $\mathscr{L}_k$ for
$1 \leq k \leq m-2$ are positive. It follows from assertion (v)
in the preceding lemma that no eigenvalue of $\mathscr{L}_k$ is
zero. To see that the eigenvalues are nonnegative we use the fact that
they are independent of the path $x$ and the index of integrability $p$.
For $k<m$ the map $\mathscr{L}_k$ maps $\mathscr{E}_{x,-1}$ to
itself and moreover $\mathscr{L}_k$ is selfadjoint with
respect to the $L^2$-inner product 
$L^2([0,1],T^* M)$ by assertion (iv). It follows that
all the eigenvalues of
$\mathscr{L}_k$ are real and hence all eigenvalues of $\mathscr{L}_k^2$
are nonnegative. Using independence of the eigenvalues from
$x$ and $p$ the claim follows. 

For $1 \leq k \leq m-2$ 
denote by $\Pi_{k,\lambda}$ for $\lambda \in \sigma(\mathscr{L}^2_k)$ the
projection to the eigenspace of $\lambda$ in $\mathscr{E}_{x,-1}$. For
$3 \leq k \leq m$ we
now extend the involution $\mathscr{I}_k$ to $\mathscr{E}_{x,-1}$ by
the formula
$$\mathscr{I}_k \xi:=\sum_{\lambda \in \sigma(\mathscr{L}^2_{k-2})}
\frac{1}{\sqrt{\lambda}}\mathscr{L}_{k-2} \circ \Pi_{k-2,\lambda}\xi,
\quad \xi \in \mathscr{E}_{x,-1}.$$

It remains to extend the involutions also the $\mathcal{E}_{x,1}$ the
eigenspace to the eigenvalue $1$ of the first involution. To do that
we introduce the maps 
$$\mathscr{H}_k \colon \mathscr{P}_{k+1} \to \mathscr{P}_k, \quad
\mathscr{D}_k \colon \mathscr{P}_k \to \mathscr{P}_{k+1}, \quad
k \in \mathbb{N}$$
by 
$$\mathscr{H}_k x(t):=x(t/2), \quad 0 \leq t \leq 1,\,\,x \in 
\mathscr{P}_{k+1}$$
and
$$\mathscr{D}_k x(t)=\left\{\begin{array}{cc}
x(2t) & 0 \leq t \leq 1/2\\
R(x(2-2t)) & 1/2 \leq t \leq 1.
\end{array}\right.$$
We extend these maps in the obvious way to bundle maps
$$\mathscr{H}_k \colon \mathscr{E}_x \to \mathscr{E}_{\mathscr{H}_k(x)},
\quad x \in \mathscr{P}_{k+1}$$
and
$$\mathscr{D}_k \colon \mathscr{E}_x \to \mathscr{E}_{\mathscr{D}_k(x)},
\quad x \in \mathscr{P}_k$$
by setting
$$\mathscr{H}_k\xi(t):=\xi(t/2), \quad 0 \leq t \leq 1,\,\,\xi \in
\mathscr{E}_x$$
and
$$\mathscr{D}_k \xi(t)=\left\{\begin{array}{cc}
\xi(2t) & 0 \leq t \leq 1/2\\
R^*\xi(2-2t) & 1/2 \leq t \leq 1.
\end{array}\right.$$
Note that
$$\mathscr{H}_k \circ \mathscr{D}_k=\mathrm{id}|_{\mathscr{E}_x}.$$
We now define recursively for $\xi \in \mathscr{E}_{x,1}$ 
$$\mathscr{I}_{k+1} \xi:=\mathscr{D}_{k}\circ \mathscr{I}_k \circ
\mathscr{H}_k \xi, \quad k \in \{1,\ldots, m-1\}.$$
\begin{rem}\label{commute}
Recall that the vertical differential of the section $\mathcal{F}$ is
given by
$$D_u\xi=\partial_s \xi+J(u)\partial_t \xi+\nabla_\xi J(u)\partial_t u$$ 
where $u \in \mathcal{F}^{-1}(0)$, $\xi \in T_u \mathcal{B}$, and 
$\nabla$ denotes the Levi-Civita connection of the metric
$g(\cdot,\cdot)=\omega(\cdot,J\cdot)$.
Observe that the last term does anticommute with the almost complex
structure $J$. The compact
operator in assertion (ii) of Theorem~\ref{invo} is given by 
$$Q_u\xi=\nabla_\xi J(u)\partial_t u, \quad \xi \in T_u\mathcal{B}.$$
If the almost
complex structure is integrable, i.e. $\nabla J=0$, then
$Q_u$ vanishes and $D_u$
commutes with $J$ and hence interchanges the
two involutions. 
\end{rem}

\subsection[Kuranishi structures]{Kuranishi structures} 

We recall here the definition of Kuranishi structure as defined in
\cite{fukaya-ono2}, see also \cite{fukaya-oh-ohta-ono}. Our definition
will be less general than the one in \cite{fukaya-ono2} since our 
Kuranishi neighbourhoods consist of manifolds instead of orbifolds,
which is sufficient for our purposes. 
Let $X$ be a compact topological Hausdorff
space. A Kuranishi structure assigns
$(V_p, E_p, \psi_p,s_p)$ to each $p \in X$ and 
$(V_{pq},\hat{\phi}_{pq}, \phi_{pq} )$ to points
$p,q \in X$ which are close to each other. They are required to
satisfy the following properties:
\begin{description}
 \item[K1:] $V_p$ is a smooth manifold, and
  $E_p$ is a smooth vector bundle on it.
 \item[K2:] $s_p$ is a continuous section of $E_p$.
 \item[K3:] $\psi_p$ is a homeomorphism from $s_p^{-1}(0)$ to
  a neighbourhood of $p \in X$.
 \item[K4:] $V_{pq}, \hat{\phi}_{pq},\phi_{pq}$ are defined
  if $q \in \psi_p(s_p^{-1}(0))$.
 \item[K5:] $V_{pq}$ is an open subset of $V_q$ containing
  $\psi_q^{-1}(q)$.
 \item[K6:] $(\hat{\phi}_{pq},\phi_{pq})$ is a map of vector
  bundles $E_q|_{V_{pq}} \to E_p$.
 \item[K7:] $\hat{\phi}_{pq}s_q=s_p\phi_{pq}$.
 \item[K8:] $\psi_q=\psi_p \phi_{pq}$.
 \item[K9:] If $r \in \psi_q(s_q^{-1}(0) \cap V_{pq})$, then
  $$\hat{\phi}_{pq} \circ \hat{\phi}_{qr}=\hat{\phi}_{pr}$$
  in a neighbourhood of $\psi_r^{-1}(r)$.
 \item[K10:] $\mathrm{dim}V_p-\mathrm{rank}E_p$ does only depend on
  the connected component of $X$ in which $p$ lies and is called the 
  \emph{virtual dimension} of the 
  Kuranishi structure of the connected component of $p$.
\end{description}

Following \cite{fukaya-ono2} we will say that our Kuranishi structure
has a \emph{tangent bundle} if there exists a family of isomorphisms
$$\Phi_{pq}\colon N_{V_p}V_q \cong E_p/E_q$$
satisfying the usual compatibility conditions. Here $N_{V_p}V_q$ denotes
the normal bundle of $V_q$ in $V_p$.

\begin{fed}
We say that a compact topological Hausdorff
space $X$ and a sequence of continuous 
involutions $\{I_k\}_{1 \leq k \leq m}$ defined on closed subspaces 
$X_m \subset X_{m-1} \cdots \subset X_1=X$ of
$X$ are of \textbf{\emph{Arnold-Givental type}}, if the following holds.
\begin{description}
 \item[(i)] The domain of the first involution $X_1$ is the whole space
  $X$, and the domain $X_k$ of $I_k$ for $2 \leq k \leq m$ is the
  fixpoint set of the previous involution, i.e. 
  $$X_k=\mathrm{Fix}(I_{k-1}).$$
 \item[(ii)] The last involution acts freely, i.e.
  $$\mathrm{Fix}(I_m)=\emptyset.$$ 
\end{description}
\end{fed}

\begin{fed}
We say that a space of Arnold-Givental type $(X,\{I_k\}_{1 \leq k \leq m})$
has a \emph{\textbf{Kuranishi structure}} if $X$ admits a Kuranishi
structure in the sense of Fukaya-Ono, such that in addition for each
$p \in X_k$ there exist involutions $I_{p,j}$ for $1 \leq j \leq k$
where $I_{p,1}$ is defined on $V_{p,1}:=V_p$ and $I_{p,j}$ is defined on
$V_{p,j}:=\mathrm{Fix}(I_{p,j-1})$ for $2 \leq j \leq k$, and
extensions of $I_{p,j}$ to smooth bundle involutions 
$$\hat{I}_{p,j}\colon E_p|_{V_{p,j}} \to E_p|_{V_{p,j}},$$ 
where the following conditions are satisfied.
\begin{description}
 \item[(i)] If 
  $p \in X_k \setminus X_{k+1}$ and $q \in \psi_p(s^{-1}_p(0)) \cap X_j$,
  then $j \leq k$, $I_{p,k}$ acts freely on $V_{p,k}$, and 
  $$I_j(q)=\psi_p \circ I_{p,j} \circ \psi_p^{-1}(q).$$
 \item[(ii)] If $p,q \in X$ are close enough, $x \in V_{pq} \cap V_{q,j}$,
  and $\xi \in (E_q)_x$
  then
  $$I_{p,j}(\phi_{pq}(x))=\phi_{pq}(I_{q,j}(x))$$
  and
  $$\hat{I}_{p,j} \circ \hat{\phi}_{pq} \xi=\hat{\phi}_{pq}\circ 
  \hat{I}_{q,j} \xi.$$
 \item[(iii)] The bundle involutions commute on their common domain of
  definition, i.e. if $x \in V_{p,j}$ and $\xi \in (E_p)_x$, then
  $$\hat{I}_{p,\ell} \circ \hat{I}_{p,j}\xi=\hat{I}_{p,j} 
  \circ \hat{I}_{p,\ell}\xi$$
  for $1 \leq \ell \leq j$.
\end{description}
\end{fed}

\begin{rem}\label{restriction}
Assume that $(X,\{I_k\}_{1 \leq k \leq m})$ has Arnold-Givental type. 
Then also
$(X_j,\{I_k\}_{j \leq k \leq m})$ for
$1 \leq j \leq m$ has Arnold-Givental type. Moreover, if
$X$ has a Kuranishi structure, then also 
$X_j$ has a Kuranishi structure. A 
Kuranishi neighbourhood is constructed in the following way. 
For $p \in X_k$
with $k \geq j$ take 
$$V_{p,j}=\mathrm{Fix}(I_{p,j-1}|_{V_p})$$
and as obstruction bundle take the intersection of the eigenspaces to
the eigenvalue $1$ of the previous involutions, i.e. 
$$E_{p,j}:=\bigcap_{1 \leq i \leq j-1}\mathrm{ker}(\hat{I}_{p,i}|_{E_p|V_p^j}
-\mathrm{id}|_{E_p|V_p^j}).$$
The other ingredients of the Kuranishi structure are then given by the
obvious restrictions. Note that since the involutions commute on their
common domain of definition $E_{p,j}$ is invariant under $\hat{I}_{p,i}$
for $j \leq i \leq k$.  
\end{rem}

\begin{fed}
We say that a space $(X,\{I_k\}_{1 \leq k \leq m})$ of Arnold-Givental
which admits a Kuranishi structure has a \emph{\textbf{tangent bundle}}
if the Kuranishi spaces $X_j$ admit a tangent bundle in the sense of
Fukaya-Ono and the isomorphisms $\Phi_{pq,j}$ for $p, q \in X_j$ close
enough are
obtained by restriction of $\Phi_{pq,1}=\Phi_{pq}$ for $1 \leq j \leq m$.
\end{fed}

In order to do useful perturbation theory we have to extend the
involutions to the tangent bundle. Since by assertion (ii) of 
Theorem~\ref{invo} the vertical differential of the section
$\mathcal{F}$ commutes with the involutions only modulo a compact
operator the tangent bundle will in general only admit
a ``stable'' Arnold-Givental structure. A similar phenomenon appeared in
\cite{fukaya-ono2} where the Kuranishi structure in general only
admitted a stable almost complex structure. 

We first have to recall the following terminology from
\cite{fukaya-ono2}. 
A tuple 
$((F_{1,p},F_{2,p}),(\Phi_{1,pq},\Phi_{2,pq},\Phi_{pq}))$ is a bundle
system over the Kuranishi space $X=(X,(V_p,E_p,\psi_p,s_p))$ if 
$F_{1,p}$ and $F_{2,p}$ are two vector bundles over $V_p$ for
every $p \in X$, $\Phi_{1,pq}\colon F_{1,q} \to F_{1,p}|_{V_q}$ and
$\Phi_{2,pq} \colon F_{2,q} \to F_{2,p}|_{V_q}$ are embeddings for
$q$ sufficiently close to $p$, and 
$$\Phi_{pq} \colon \frac{F_{1,p}|_{V_q}}{F_{1,q}} \to 
\frac{F_{2,p}|_{V_q}}{F_{2,q}}$$
are isomorphisms of vector bundles. These maps are required to satisfy
some compatibility conditions. Moreover, there
is some obvious notion of isomorphism, Whitney sum,
tensor product etc. for bundle systems. 
We refer the reader to \cite{fukaya-ono2} for details. 

If a Kuranishi structure has a tangent bundle, i.e. a family of 
isomorphisms $\Phi_{pq} \colon N_{V_p}V_q \cong E_p/E_q$, 
then one can define a
bundle system 
$$TX=(TV_p, E_p, \hat{\phi}_{pq},d \phi_{pq}, \Phi_{pq})$$
over X.  
If the space $X$ is of Arnold-Givental type and admits a tangent bundle,
then we define the normal bundle
\begin{eqnarray*}
NX&=&\{TV_{p,j}|_{V_{p,j+1}}/TV_{p,j+1},\\
& &E_{p,j}|_{V_{p,j+1}}/E_{p,j+1},\\
& &d\phi_{pq,j}|_{V_{pq,j+1}}/d\phi_{pq,j+1},\\
& &\hat{\phi}_{pq,j}|_{V_{pq,j+1}}/\hat{\phi}_{pq,j+1},\\
& &\Phi_{pq,j}|_{V_{pq,j+1}}/\Phi_{pq,j+1}\}_{1 \leq j \leq m-1}
\end{eqnarray*}
as a bundle system over $\bigcup_{j=1}^{m-1}X_{j+1}$.
The real $K$-group $KO(X)$ of
a space with Kuranishi structure $X$ was defined
in \cite{fukaya-ono2} as the quotient of the
free abelian group generated by the set of all isomorphism classes of
bundle systems modulo the relations
$$
[((F_{1,p},F_{2,p}),(\Phi_{1,pq},\Phi_{2,pq},\Phi_{pq}))
\oplus((F'_{1,p},F'_{2,p}),(\Phi'_{1,pq},\Phi'_{2,pq},\Phi'_{pq}))]=$$
$$[((F_{1,p},F_{2,p}),(\Phi_{1,pq},\Phi_{2,pq},\Phi_{pq}))]+
[((F'_{1,p},F'_{2,p}),(\Phi'_{1,pq},\Phi'_{2,pq},\Phi'_{pq}))]$$
$$[((F_{1,p},F_{2,p}),(\Phi_{1,pq},\Phi_{2,pq},\Phi_{pq}))]=0$$
$$\mathrm{if}\,\,((F_{1,p},F_{2,p}),(\Phi_{1,pq},\Phi_{2,pq},\Phi_{pq}))
\,\,\textrm{is trivial.}$$
If the Kuranishi structure has a tangent bundle
$[TX]$ denotes the class of $TX$ in $KO(X)$. If the Kuranishi structure
is of Arnold-Givental type and admits a tangent bundle
$[NX]$ denotes the class of $NX$ in $\bigoplus_{j=2}^{m}KO(X_j)$.

\begin{fed}
Assume that $(X,\{I_k\}_{1 \leq k \leq m})$ is a space of Arnold-Givental type
which admits a Kuranishi structure
 $(X,(V_p,E_p,s_p,\psi_p,\phi_{pq},
\hat{\phi}_{pq}))$ and assume further that 
$((F_{1,p},F_{2,p}),(\Phi_{1;pq},
\Phi_{2;pq},\Phi_{pq}))$ is a bundle system on
it. We say that 
$((F_{1,p},F_{2,p}),(\Phi_{1;pq},\Phi_{2;pq},\Phi_{pq}))$ 
is a \emph{
\textbf{Bundle System of Arnold-Givental type}} if for every $p \in X_k$ 
there exist extensions of the
involutions $I_{p,j}$ for $1 \leq j \leq k$ to
smooth involutative bundle maps
$\hat{I}_{1,p,j} \colon F_{1,p}|_{V_{p,j}} \to F_{1,p}|_{V_{p,j}}$ and 
$\hat{I}_{2,p,j} \colon F_{2,p}|_{V_{p,j}} \to F_{2,p}|_{V_{p,j}}$ such 
that the following conditions are satisfied.
\begin{description}
 \item[Compatibility:] The transition maps $\Phi_{pq}$, $\Phi_{1;pq}$,
  and $\Phi_{2;pq}$ restricted to the domain of definition of the
  involutions interchange them.
 \item[Commutativity:] The involutions commute on their common domain of
  definition.
\end{description}
\end{fed}
The Whitney sum of two bundle systems of Arnold-Givental type still has 
Arnold-Givental type and hence one can consider the $K$-group
$K_{AG}(X)$ of bundle systems of Arnold-Givental type over $X$. 
There is an obvious map
$$K_{AG}(X) \to KO(X).$$
If $X$ admits a tangent bundle then there is an obvious extension of 
$I_{p,j}$ to the bundle $E_p$ given by 
$\hat{I}_{p,j}$. However there is no obvious extension of 
$I_{p,j}$ to $T V_p$ in general. This motivates the following definition.
\begin{fed}
Assume that $(X,\{I_k\}_{1 \leq k \leq m})$ is a space of
Arnold-Givental type which has a Kuranishi structure with tangent
bundle. We say that $X$ has a \emph{\textbf{Kuranishi structure of
Arnold-Givental type}} if the normal bundle $NX$ 
is a bundle system of Arnold-Givental
type in $\bigcup_{j=2}^{m}X_j$. 
We say that $X$ has a \emph{\textbf{Kuranishi structure of
stable Arnold-Givental type}} if $[NX]$ is in the image of 
$\bigoplus_{j=2}^m K_{AG}(X_j) \to \bigoplus_{j=2}^m KO(X_j)$. 
\end{fed}
Note that a Kuranishi structure of Arnold-Givental type is also a 
Kuranishi structure of stable Arnold-Givental type.

\begin{thm}\label{kurtra}
We assume that $(X,\{I_k\}_{1 \leq k \leq m})$ has a Kuranishi
structure of Arnold-Givental type whose virtual dimension is zero. 
Then for each $p \in X$, there
exist smooth sections $\tilde{s}_p$ such that the following
holds.
\begin{description}
 \item[(i)] $\tilde{s}_p\circ \phi_{pq}=\hat{\phi}_{pq}\circ 
 \tilde{s}_q$,
 \item[(ii)] Let $x \in V_{pq}$. Then the restriction of the
  differential of the composition of $\tilde{s}_p$ and the projection 
  $E_p \to E_p/E_q$ coincides with the isomorphism 
  $\Phi_{pq}\colon N_{V_p}V_q \cong E_p/E_q$.
 \item[(iii)] The sections $\tilde{s}_p$ are transversal to $0$, and if
  $p \in X_k$ then $\tilde{s}_p^{-1}(0)$ is invariant under 
  $I_{p,j}$ for $1 \leq j \leq k$.  
\end{description}
\end{thm}
\begin{rem}
The sections $\tilde{s}_p$ will not necessarily be
invariant under the involutions, only their zero set will be.
\end{rem}
We need the following lemma.

\begin{lemma}\label{chnorz}
Assume that $Y$ is a compact manifold, $N$ and $E$ are two vector bundle
over $Y$ and $U \subset N$ is an open neighbourhood of $Y \subset N$. 
Denote by $\pi_N \colon N \to Y$ the canonical projection. 
Then there exists a section $s \colon N \to \pi^*_N E$ which satisfies
the following conditions.
\begin{description}
 \item[(i)] The zero section $s^{-1}(0)$ is invariant under the 
  involutative bundle map on $N$ defined by
  $(n,y) \mapsto (-n,y)$ for $y \in Y$ and $n \in N_y$.
 \item[(ii)] Outside $U$ the section $s$ is invariant under the involutative
  bundle map of $\pi_N^* E$ given by $(e,n,y) \mapsto (-e,-n,y)$
  for $y \in Y$, $n \in N_y$, and $e \in \pi_N^*E_{(n,y)}$.
 \item[(iii)] The boundary of the set of transversal points has
  codimension at least one. More precisely, there exists a manifold 
  $\Omega$ of dimension 
  $\mathrm{dim}(\Omega)=\mathrm{dim}(Y)-\mathrm{rk}(E)+\mathrm{rk}(N)-1$
  and a smooth map $f \colon \Omega \to N$ such that the $\omega$-limit
  set of the $\mathrm{dim}(Y)-\mathrm{rk}(E)+\mathrm{rk}(N)$-dimensional
  manifold of points of transversal intersection
  $$\mathscr{T}:=\{n \in s^{-1}(0): Ds(n)\,\,onto\}$$
  is contained in the image of $\Omega$, i.e.
  $$\bigcap_{K \subset \mathscr{T}\,\,compact}\mathrm{cl}
  (\mathscr{T}\setminus K) \subset f(\Omega).$$
\end{description}
\end{lemma}
\textbf{Proof:} Choose a bundle map $\Phi \colon N \to E$ and define
$$Q:=\{y \in Y: \mathrm{dim}(\mathrm{ker}\Phi(y))>0\}.$$
Let $\Psi \colon Y \to E$ be a section such that 
$$\Psi|_Q=0, \quad \Psi(y) \cap \Phi(N_y)=\{0\},\quad \forall\,\,
y \in Y.$$
Choose further a smooth cutoff function $\beta \colon N \to \mathbb{R}$
such that
$$\beta|_Y=1, \quad \mathrm{supp}(\beta) \subset U.$$
Define a section $s \colon N \to \pi_N^* E$ by
$$s(n):=\beta(n)\pi_N^*\Psi(\pi_N(n))+\pi_N^*\Phi(n), \quad n \in N.$$
Then $s$ satisfies conditions (i) and (ii) and for generic choice of
$\Psi$ and $\Phi$ also condition (iii). This proves the lemma. 
\hfill $\square$
\\ \\
\textbf{Proof of Theorem~\ref{kurtra}:} We continue the notation of
Remark~\ref{restriction}. For $j \in \{1,\cdots, m\}$ denote
by $s_p^j$ for $p \in X_j$ the section from $V_{p,j}$ to $E_{p,j}$ which is
induced from $s_p$. By induction on $j$ from $m$ to $1$ we find 
sections $\tilde{s}_p^j$ of $s_p^j$ which satisfy conditions
(i) and (ii) of Theorem~\ref{kurtra} as well invariance of the zero set
under the involutions, but instead of the transversality condition we
impose the condition that the boundary of the manifold of points of
transversal intersection has codimension at least one.
\begin{description}
 \item[(iii`)] For each $j \in \{1, \cdots, m\}$ and for each $p \in X_j$ 
  there exists a manifold $\Omega_p^j$ of dimension
  $\mathrm{dim}(\Omega_p^j)=d^j_p-1$ where $d^j_p$ is the virtual
  dimension of the connected component of $p$ of the
  Kuranishi structure of $\mathrm{X}_j$ and a smooth map 
  $f_p^j \colon \Omega_p^j \to V_p^j$ such that the 
  $\omega$-limit set of the set of points of transversal intersection
  $$\mathscr{T}^j_p:=\{q \in (\tilde{s}^j_p)^{-1}(0): D\tilde{s}_p^j(q) 
  \,\,onto\}$$
  is contained in the image of $\Omega_p^j$, i.e.
  $$\bigcap_{K \subset \mathscr{T}^j_p\,\,compact} \mathrm{cl}(
  \mathscr{T}^j_p \setminus K) \subset f_p^j(\Omega_p^j).$$
\end{description}
To prove the induction step we observe that by the assumption that the
Kuranishi structure is of Arnold-Givental type it follows that
the two vector bundles $TV_{p,j}|_{(s_p^{j+1})^{-1}(0)}/
TV_{p,j+1}|_{(s_p^{j+1})^{-1}(0)}$ and
$E_{p,j}|_{(s_p^{j+1})^{-1}(0)}/E_{p,j+1}|_{(s_p^{j+1})^{-1}(0)}$
induce bundles on the quotient $(s_p^{j+1})^{-1}(0)/I_{p,j+1}$. The
induction step can now be concluded from Lemma~\ref{chnorz}.
Since the Kuranishi structure of $X=X_1$ has virtual dimension zero it
follows that for every $p \in X_1$ the set $\Omega^1_p$ is empty and
hence condition (iii) holds. This proves the theorem. \hfill $\square$
\\ \\
We recall from \cite{fukaya-ono2} that if $Y$ is a topological space and
$X$ is a space with Kuranishi structure, a \emph{strongly continuous
(smooth)
map} $f: X \to Y$ is a family of continuous (smooth) 
maps $f_p \colon V_p \to Y$
for each $p \in X$ such that $f_p \circ \phi_{pq}=f_q$. If $n$ is the
virtual dimension of the Kuranishi structure we define the
homology class
$$f([X]) \in H_n(Y;\mathbb{Z}_2)$$
as in \cite{fukaya-ono2}. We are now able to draw the following
Corollary from Theorem~\ref{kurtra}.

\begin{cor}
Assume that $(X,\{I_k\}_{1 \leq k \leq m})$ has a Kuranishi
neighbourhood of stable Arnold-Givental type whose virtual dimension is
zero. Suppose further that $Y$ is a topological space and
$f \colon X \to Y$ is a strongly continuous map. Then
$$f([X])=0 \in H_0(Y;\mathbb{Z}_2).$$
\end{cor}

If the space $\mathcal{M}$ consisting of $J$-holomorphic disks whose
boundary is mapped to the Lagrangian $L$ is compact with respect to
the Gromov topology, then
$\mathcal{M}$ together with the involutions 
described in section~\ref{agi} is of Arnold-Givental type. We next prove
that under the compactness assumption $\mathcal{M}$ has a Kuranishi 
neighbourhood of Arnold-Givental type. In our proof we mainly follow
\cite{fukaya-ono2}. The new ingredient is to choose the
obstruction bundle in such a way that it is invariant under the involutions. 

In general one cannot expect that $\mathcal{M}$ is compact due to the
bubbling phenomenon. We hope that the approach pursued in this article
together with the techniques developped in \cite{fukaya-ono2} and
\cite{fukaya-oh-ohta-ono} will allow us to show that the
compactification of $\mathcal{M}$ by bubble trees has a Kuranishi
neighbourhood of Arnold-Givental type. 

\begin{thm}\label{kusta}
Assume that $\mathcal{M}$ is compact with respect to the Gromov topology.
Then $(\mathcal{M},\{I_k\}_{1 \leq k \leq m})$ admits
a Kuranishi structure of stable Arnold-Givental type. 
\end{thm}
\textbf{Proof:} We first choose local trivialisations of the vector
bundle $\mathcal{E}$. More precisely, for each  $q \in \mathcal{B}$ we
choose an open neighbourhood  $U_q \subset \mathcal{B}$ of $q$
and a smooth family of Banach space isomorphisms
$$\mathcal{P}^p_q \colon \mathcal{E}_q \to \mathcal{E}_p, \quad
p \in U_q$$
such that the following conditions are satisfied.
\begin{description}
 \item[(T1)] If $q \in \mathcal{B}_k \setminus \mathcal{B}_{k+1}$ for
  $k \in \mathbb{N}$, then $U_q$ is invariant under $I_j$ for
  $1 \leq j \leq k$ and $I_k$ acts freely on $U_q$.
 \item[(T2)] For $1 \leq j \leq k$ the trivialisations commute with the
  involutions, i.e.
  $$\mathcal{P}_q^p \circ I^{\mathcal{E}}_j=I^{\mathcal{E}}_j \circ
  \mathcal{P}_{I_j q}^{I_j p}, \quad p \in U_q.$$
\end{description}

Now choose for every $p \in \mathcal{M}$ an open neighbourhood 
$\hat{V}_p \subset U_p$ of $p$ in $\mathcal{B}$,
choose a finite set $Q \subset \mathcal{M}$, for each $q \in Q$ a closed
neighbourhood $\hat{U}_q \subset U_q$ of $q$ in $\mathcal{B}$, 
and a finite dimensional
subspace $\hat{E}_q \subset \mathcal{E}_q$ consisting
of smooth sections with the following properties:
\begin{description}
 \item[(i)] For every $k \in \mathbb{N}$ the sets $\hat{V}_p$ for every 
  $p \in \mathcal{M} \cap (\mathcal{B}_k \setminus \mathcal{B}_{k+1})$
  and the set
  $\bigcup_{q \in Q \cap (\mathcal{B}_k\setminus \mathcal{B}_{k+1})}
  \hat{U}_q$ 
  are invariant under $I_j$ for $1 \leq j \leq k$ and $I_k$ acts freely
  on them.
 \item[(ii)] For $1 \leq k \leq m$ the family of vectorspaces
  $\bigcup_{q \in \mathcal{M}_k} \hat{E}_q$ is
  invariant under $I_k^{\mathcal{E}}$.
 \item[(iii)] For $p \in \mathcal{B}$ let $Q_p:=\{q \in Q: p \in \hat{U}_q\}$.
  Assume that $p \in \mathcal{M}$ and $p' \in \hat{V}_p$.
  Then the sum $\bigoplus_{q \in Q_p}\mathcal{P}_q^{p'}\hat{E}_q$ is
  direct, i.e. for every
  $q_0 \in Q_p$ we have $\mathcal{P}_{q_0}^{p'}\hat{E}_{q_0} 
  \cap \sum_{q \in Q_p \setminus
  \{q_0\}}\mathcal{P}_q^{p'}\hat{E}_q=\emptyset.$ 
 \item[(iv)] For $p \in \mathcal{M}$ and $p' \in \hat{V}_p$
  the operator
  $$\Pi^{p'}_p \circ D_{p'} \colon T_{p'}\mathcal{B} \to \mathcal{E}_{p'}
  /\bigoplus_{q \in Q_p} \mathcal{P}^{p'}_q\hat{E}_q$$
  is surjective, where 
  $$\Pi^{p'}_p \colon \mathcal{E}_{p'} \to \mathcal{E}_{p'}
  / \bigoplus_{q \in Q_p} \mathcal{P}_q^{p'}\hat{E}_q$$
  denotes the canonical projection.
\end{description}

Our first aim is to define the manifold $V_p$ of $p$ which occurs
in the definition of Kuranishi structure. In our construction this 
manifold will be a small neighbourhood of zero in the kernel of the map
$\Pi^p_p \circ D_p$. The size of this neighbourhood will
depend on the domain of definition of the smooth injective evaluation maps
$$\mathrm{ev}_p \colon V_p \to \hat{V}_p$$
which we have to define first. 

By (iv) we can choose a smooth family of uniformly 
bounded right inverses
$$R^{p'}_p \colon \mathcal{E}_{p'}/\bigoplus_{q \in Q_p} 
\mathcal{P}_q^{p'}\hat{E}_q \to T_{p'}\mathcal{B}$$
of $\Pi^{p'}_p \circ D_{p'}$, i.e.
$$\Pi^{p'}_p \circ D_{p'} \circ R^{p'}_p =\mathrm{id}
|_{\mathcal{E}_{p'}/\bigoplus_{q \in Q_p}\mathcal{P}_q^{p'}\hat{E}_p}.$$
We need in addition some compatibility of the right inverses with the
involutions. To state it we observe that
it follows from condition (T2) on the local trivialisations together
with conditions (i) and (ii) that for $p \in \mathcal{M}_k$ and
$p' \in \hat{V}_p \cap \mathcal{B}_j$ for $j \leq k$ the family
of vector spaces
$\bigoplus_{q \in Q_p} \mathcal{P}_q^{p'}\hat{E}_q \cup
\bigoplus_{q \in Q_{I_i p}}\mathcal{P}^{I_i p'}_{I_i q}\hat{E}_{I_i q}$
is invariant
under $I_i$ for $1 \leq i \leq j$. Hence the
involutions $I_j^{\mathcal{E}}$ induce involutions 
$$I_j^{\mathcal{E},p} \colon \bigcup_{p' \in \hat{V}_p \cap
\mathcal{B}_j}\mathcal{E}_{p'}/\bigoplus_{q \in Q_p}\mathcal{P}^{p'}_q 
\hat{E}_q \to
\bigcup_{p' \in \hat{V}_p \cap
\mathcal{B}_j}\mathcal{E}_{p'}/\bigoplus_{q \in Q_p}\mathcal{P}^{p'}_q 
\hat{E}_q.
$$
Using the fact that by (ii) of Theorem~\ref{invo} the operator
$I^{\mathcal{E}}_k \circ D_u-D_u\circ I^{T\mathcal{B}}_k$
vanishes on $T_u\mathcal{B}_k$ together with (iii) of Theorem~\ref{invo},
we can impose the following
compatibility condition of the right inverse $R^{p'}_p$ and the 
involutions
\begin{equation}\label{inv1}
R_p^{p'} \circ I_k^{\mathcal{E},p}\big|
_{\bigcap_{j=1}^{k-1}\mathrm{ker}(I_j^{\mathcal{E},p}-\mathrm{id})}
=I_k^{T\mathcal{B}}\circ 
R_p^{p'}\big|
_{\bigcap_{j=1}^{k-1}\mathrm{ker}(I_j^{\mathcal{E},p}-\mathrm{id})}
, \quad \forall\,\,p' \in \hat{V}_p \cap \mathcal{B}_k.
\end{equation}

Now we are able to define the manifold $V_p$ which occurs in the
definition of the Kuranishi structure. 
For $p \in \mathcal{B}$ and for $\xi \in T_p\mathcal{B}$ small enough
we define 
$$\exp_p \xi \in \mathcal{B}$$
as the pointwise exponential map with respect to the metric
$g(\cdot,\cdot)=\omega(\cdot,J\cdot)$. Note that if $p \in
\mathcal{B}_k$ and $\xi \in T_p \mathcal{B}_k$ then 
$\exp_p \xi \in \mathcal{B}_k$. Now we choose $V_p$ as an open
neighbourhood of zero in 
$\mathrm{ker}(\Pi_p\circ D_p) \subset T_p \mathcal{B}$ 
which is invariant under $I^{T\mathcal{B}}_k$ if $p \in \mathcal{B}_k$
and which is so small that  we are able to define
$$\mathrm{ev}^0_p \colon V_p \to \hat{V}_p, \quad
\mathrm{ev}^0_p:=\exp_p|_{V_p}$$
and be recursion for $\nu \in \mathbb{N}$
$$\mathrm{ev}^\nu_p\colon V_p \to \hat{V}_p,\quad \xi \mapsto
\exp_{\mathrm{ev}^{\nu-1}_p \xi}\bigg(R^{\mathrm{ev}^{\nu-1}_p \xi}_p \circ 
\Pi^{\mathrm{ev}^{\nu-1}_p \xi}_p 
\circ \mathcal{F} \circ \mathrm{ev}^{\nu-1}_p \xi\bigg),
$$
and finally
$$\mathrm{ev}_p \colon V_p \to \hat{V}_p, \quad \mathrm{ev}_p:=
\lim_{\nu \to \infty} \mathrm{ev}^\nu_p.$$
We now define for each $p \in \mathcal{M}$ the obstruction bundle
$E_p \to V_p$ by
$$E_p:=(\mathrm{ev}_p)^*\bigoplus_{q \in
Q_p}\mathcal{P}_q^{\mathrm{ev}_p}\hat{E}_q,$$
and the section $s_p \colon V_p \to E_p$ by
$$s_p:=(\mathrm{ev}_p)^*\mathcal{F}.$$
The homeomorphisms $\psi_p$ from $s_p^{-1}(0)$ to a neighbourhood of
$p \in \mathcal{M}$ are defined  by
$$\psi_p:=\mathrm{ev}_p|_{s_p^{-1}(0)}.$$
Perhaps after shrinking the manifolds $V_p$ we may assume using the
assumption that the set $\hat{U}_q$ are closed for every $q \in Q$ that 
\begin{equation}\label{Vpq}
Q_{\psi_p(x)} \subset Q_p, \quad \forall \,\,p \in \mathcal{M},\,\,\forall
\,\,x \in s_p^{-1}(0).
\end{equation}
This implies that for $p \in \mathcal{M}$ and $x \in s_p^{-1}(0)$ the
set 
$$V_{p
\psi_p(x)}:=\mathrm{ev}_{\psi_p(x)}^{-1}
\bigg(\mathrm{ev}_{\psi_p(x)}(V_{\psi_p(x)})
\cap \mathrm{ev}_p(V_p)\bigg)$$
is open in $V_{\psi_p(x)}$.
We now define 
$$\phi_{p\psi_p(x)}\colon V_{p\psi_p(x)} \to V_p, \quad 
\phi_{p\psi_p(x)}:=\mathrm{ev}_p^{-1}\circ \mathrm{ev}_{\psi_p(x)}.$$
Using again (\ref{Vpq}) we observe that for every $y \in V_{p\psi_p(x)}$
we have
$$\bigoplus_{q \in Q_{\psi_p(x)}}\mathcal{P}_q^{\mathrm{ev}_{\psi_p(x)}(y)}
\hat{E}_q \subset 
\bigoplus_{q \in Q_p}\mathcal{P}_q^{\mathrm{ev}_p \circ \phi_{p\psi_p(x)}(y)}
\hat{E}_q.$$
We now define the bundle maps $\hat{\phi}_{p\psi_p(x)} \colon
E_{\psi_p(x)}|_{V_{p\psi_p(x)}} \to E_p$ as the map induced by the above
inclusion. To define the isomorphisms $\Phi_{p\psi_p(x)} \colon N_{V_p}
V_{\psi_p(x)} \to E_p/E_{\psi_p(x)}$ which are required for the tangent
bundle of the Kuranishi structure, we observe that there is
a natural identification of the  normal bundle of
$\mathrm{ev}_{\psi_p(x)}(V_{p\psi_p(x)})$ in $\mathrm{ev}_p(V_p)$ 
$$N_{\mathrm{ev}_p(V_p)} \mathrm{ev}_{\psi_p(x)}(V_{p\psi_p(x)})
\cong \bigoplus_{q \in Q_p \setminus Q_{\psi_p(x)}}\mathcal{P}_q^
{\mathrm{ev}_{\psi_p(x)}}\hat{E}_q$$
due to the fact that 
$$\Pi_{\psi_p(x)}^{\mathrm{ev}_{\psi_p(x)}(y)} \circ 
D_{\mathrm{ev}_{\psi_p(x)}(y)}
\colon T_{\mathrm{ev}_{\psi_p(x)}(y)}\mathcal{B} \to
\mathcal{E}_{\mathrm{ev}_{\psi_p(x)}(y)}/\bigoplus_{q \in Q_{\psi_p(x)}}
\mathcal{P}_q^{\mathrm{ev}_{\psi_p(x)}(y)}\hat{E}_q$$
is already surjective. The isomorphisms $\Phi_{p\psi_p(x)}$ are defined
to be the induced isomorphisms of the identification above. Finally,
using (\ref{inv1}), the involutions $I^{\mathcal{E}}_k$ induce involutions
$\hat{I}_{p,k}$ on the obstruction bundle $E_p$ for $p \in
\mathcal{M}_k$. Hence we have proved that $\mathcal{M}$ admits a
Kuranishi structure.

We next show in the integrable case the Kuranishi structure is of 
Arnold-Givental type. As it was explained in 
Remark~\ref{commute}, if the almost complex structure is integrable
the operator $D_u$ will interchange the involutions, 
i.e. for $k \in \mathbb{N}$ and
$u \in \mathcal{F}^{-1}(0) \cap \mathcal{B}_k$ it holds that
\begin{equation}\label{comm}
I_k^{\mathcal{E}}\circ D_u=D_u \circ I_k^{T\mathcal{B}}.
\end{equation}
If $p \in \mathcal{M}_k$ and $x \in V_{p,j}$ for $j \leq k$ then 
$$T_x V_{p,j}=\mathrm{ev}_p^*\bigg(\mathrm{ker}\big(\Pi^{\mathrm{ev}_p(x)}_p
\circ D_{\mathrm{ev}_p(x)}\big) \cap \bigcap_{i=1}^{j-1}
\big(\mathrm{ker}I_i^{T\mathcal{B}}-\mathrm{id}|_{T\mathcal{B}}\big)\bigg),$$
and it follows from (\ref{comm}) and the invariance of the obstruction
bundle under the involutions, that $I_i^{T\mathcal{B}}$ for $1 \leq i
\leq j$ induce involutions on $TV_{p,j}$. Using these involutions one
can endow the normal bundle with the structure of a
bundle system of Arnold-Givental type.

It remains to treat the non-integrable case. In general (\ref{comm})
does not hold but by assertion (ii) of Theorem~\ref{invo} we can homotop
$D_u$ through Fredholm operators to a Fredholm operator $D_u^1$ which
interchanges the involutions by setting
$$D_u^\lambda:=D_u-\lambda Q_u, \quad \lambda \in [0,1].$$
Now choose $\hat{V}^1_p \subset U_p$ for every $p \in \mathcal{B}$, a
finite set $Q^1 \in \mathcal{B}$, and for every
$q \in Q^1$ a finite dimensional subspace
$\hat{E}_q^1 \subset \mathcal{E}_q$ consisting of smooth sections, which
satisfy again assertions (i) to (iii) but assertion (iv) replaces by
\begin{description}
 \item[(iv)'] For $p \in \mathcal{M}$ and $p' \in \hat{V}_p$ the
  operator 
  $$\Pi^{p',1}_p \circ D^1_{p'} \colon T_{p'}\mathcal{B}
  \to \mathcal{E}_{p'}/\bigoplus_{q \in Q^1_p}\mathcal{P}^{p'}_q
  \hat{E}^1_q$$
  is surjective. 
\end{description}
Define for $p \in \mathcal{M}_k$ and $x \in V_{p,j}$ for $j \leq k$
$$W_{p,j}:=\mathrm{ev}_p^*\bigg(\mathrm{ker}\big(\Pi^{\mathrm{ev}_p(x),1}_p
\circ D^1_{\mathrm{ev}_p(x)}\big) \cap \bigcap_{i=1}^{j-1}
\big(\mathrm{ker}I_i^{T\mathcal{B}}-\mathrm{id}|_{T\mathcal{B}}\big)\bigg).$$
Note that $W_{p,j}$ is invariant under the involutions $I_i^{T\mathcal{B}}$ 
for $1 \leq i \leq j$. One can now define a bundle system of
Arnold-Givental type $N^1 X$ over $\bigcup_{j=2}^{m}X_j$ where $F_{1,p}$
is given by $\bigcup_{j=2}^m W_{p,j}|_{V_{p,j+1}}/W_{p,j+1}$,
$F_{2,p}=\bigcup_{j=2}^m
E_{p,j}|_{V_{p,j+1}}/E_{p,j+1}$, and the transition functions
are defined in a similar manner as in the case of the normal
bundle. Using the homotopy between $D_u$ and $D^1_u$ one shows that 
$$[N^1 X]=[NX] \in \bigoplus_{j=2}^m KO(X_j).$$
It follows that the Kuranishi structure is of stable Arnold-Givental
type. This proves the Theorem. \hfill $\square$

\section[Moment Floer homology]{Moment Floer homology}\label{mofoho}

Moment Floer homology was introduces in \cite{frauenfelder2}. Moment
Floer homology is a tool to count intersection points of some
Lagrangians in Marsden-Weinstein quotients which are
fixpoint set of some antisymplectic involution. In general, due to the
bubbling phenomenon, the ordinary Floer homology for Lagrangians in
Marsden-Weinstein quotients cannot be defined by standard means, see
\cite{cho,cho-oh} for a computation of the Floer homology of Lagrangian
torus fibers of Fano toric manifolds. To overcome the bubbling problem
one replaces Floer's equations by the symplectic vortex equations to
define the boundary operator. Under some topological assumptions on the
enveloping manifold one can prove compactness of the relevant moduli
spaces of the symplectic vortex equations. In the special case where the
two Lagrangians are hamiltonian isotopic to each other one can use the
antisymplectic involution to prove that moment Floer homology is equal 
to the singular homology of the Lagrangian with coefficients in some
Novikov ring. This leads to a prove of the Arnold-Givental conjecture
for some class of Lagrangians in Marsden-Weinstein quotients which are
fixpoint sets of some antisymplectic involution. We will give in this
section proofs of the main properties of the symplectic vortex equations
to define moment Floer homology and refer to \cite{frauenfelder2} for
complete details. To compute it we will need the techniques of 
section~\ref{transi}. These techniques were not available in 
\cite{frauenfelder2} and hence moment Floer homology could there only be
computed under some additional monotonicity assumption which is removed
here.

\subsection[The set-up]{The set-up}

In this subsection we introduce the notation to define 
the symplectic vortex equations 
and formulate the hypotheses under which compactness of
the relevant moduli spaces can be proven. 

Let $G$ be a Lie group with Lie algebra $\mathfrak{g}$ which acts 
covariantly on a manifold $M$, i.e. there exists a smooth homomorphism
$\psi: G \to \mathrm{Diff}(M)$. We will often drop $\psi$ and
identify $g$ with $\psi(g)$.  

For $\xi \in \mathfrak{g}$ we denote
by $X_\xi$ the vector field $M \to TM$ which is generated by the
one-parameter subgroup generated by $\xi$, i.e.
$$X_\xi(x):=\frac{d}{dt}\bigg|_{t=0}\exp(t \xi)(x), \quad
\forall \,\,x \in M.$$
We shall use the linear mapping $L_x: \mathfrak{g} \to T_x M$ defined by
$$L_x \xi:=X_\xi(x) \in T_x M.$$
We will denote the adjoint action of $G$ on $\mathfrak{g}$ by
$$g \xi g^{-1}=\mathrm{Ad}(g)\xi
=\frac{d}{dt}\bigg|_{t=0}g \exp(t\xi) g^{-1}.$$
If $I$ is an open intervall, $t_0 \in I$,  
and $g:I \to M$ is a smooth path, we will write
$$(g^{-1}\partial_t g)(t_0):=
d\mathcal{L}_{g(t_0)}^{-1}(g(t_0))\partial_t g(t_0) \in 
T_{\mathrm{id}} G=\mathfrak{g}$$
where $\mathcal{L}_g \in \mathrm{Diff}(G)$ is the left-multiplication by $g$.

Assume that the Lie algebra is endowed with an inner product
$\langle \cdot,\cdot \rangle$ which is invariant under the adjoint action
of the Lie group. 
If $(M,\omega)$ is symplectic, we say that the action of $G$ is Hamiltonian,
if there exists a moment map for the action, i.e. 
an equivariant function $\mu: M \to \mathfrak{g}$ 
\footnote{Some authors use the convention that the moment map takes values
in the dual of the Lie algebra. Since we have an inner product we can 
identify the Lie algebra with its dual.}
where the action of $G$ on $\mathfrak{g}$ is the adjoint action, such that
for every $\xi \in \mathfrak{g}$
$$d \langle \mu (\cdot),\xi \rangle=\iota(X_\xi)\omega.$$
Note that the function
$\langle \mu( \cdot),\xi \rangle$ is a Hamiltonian function for the
vector field $X_\xi$. Observe that if  $\xi \in Z(\mathfrak{g})$
\footnote{The centraliser $Z(\mathfrak{g})$ consists of all $\xi \in \mathfrak{g}$
such that $[\xi,\eta]=0$ for every $\eta \in \mathfrak{g}$},
then $\mu_\xi(\cdot):=\mu(\cdot)+\xi$ is also a moment map for the action of 
$G$. Hence the moment map is determined by the action 
up to addition of a central element in each 
connected component of $M$.

We assume now
that  $(M,\omega)$ is a symplectic (not necessarily compact) 
connected manifold, and 
$G$ a compact connected Lie group that acts on $M$ by 
Hamiltonian symplectomorphisms as above. We assume that the action is
effective, i.e. the homomorphism $\psi$ is injective.  
Let $L_0$ and $L_1$ be two closed Lagrangian submanifolds of $M$. 
We do not require that the Lagrangians are $G$-invariant but we assume
throughout this section the following compatibility condition with $G$.
\begin{description}
 \item[(H1)] \emph{
  For $j \in \{0,1\}$ there exist antisymplectic involutions 
  $R_j \in \mathrm{Diff}(M)$, i.e.
  $$R_j^*\omega=-\omega, \quad R_j^2=\mathrm{id},$$
  which commute which $G$, i.e. for every $g \in G$ the symplectomorphism
  $R_j \psi(g) R_j$ lies in the image of $\psi$, such that
  $$L_j=\mathrm{Fix}(R_j)=\{x \in M: R_j(x)=x\}.$$}
\end{description}
The maps $R_j$ lead to Lie group Automorphisms $S_j:G \to G$
defined by
\begin{equation}\label{S}
S_j(g):=\psi^{-1}(R_j \psi(g) R_j), \quad \forall \,\, g \in G.
\end{equation}
Note that $S_j$ are involutative, i.e.
$$S^2_j=\mathrm{id}.$$
We assume that the inner product in the Lie algebra is also invariant 
under the differentials of $S_j$ at the identity. These are determined by
the formula
$$X_{dS_j(\mathrm{id})(\xi)}(x)=dR_j(R_j(x))^{-1}X_\xi(R_j x), \quad \forall\,\,
x \in M.$$
In the following we will write $\dot{S}_j$ for $dS_j(\mathrm{id})$.
If one identifies $G$ with $\psi(G)$, then formally
$$\dot{S}_j=\mathrm{Ad}(R_j).$$
Let $\mu$ be a moment map for the action of $G$ on $M$.
We further impose the following hypothesis throughout this section.
\begin{description}
\item[(H2)] \emph{The moment map $\mu$ is proper, zero is a regular value of
$\mu$, and $G$ acts freely on $\mu^{-1}(0)$, i.e. 
$\psi(g)p=p$ for $p \in \mu^{-1}(0)$ implies that $g=\mathrm{id}$.}
\end{description}
The Marsden-Weinstein quotient is defined to be the set of $G$-orbits in
$\mu^{-1}(0)$
$$\bar{M}:=M//G:=\mu^{-1}(0)/G,$$
i.e. $x, y \in \mu^{-1}(0)$ are equivalent if there exists $g \in G$ such that
$\psi(g)x=y$. It follows from hypothesis (H2) that $\bar{M}$ 
is a compact manifold of dimension
$$\mathrm{dim}(\bar{M})=\mathrm{dim}(M)-2\mathrm{dim}(G).$$
The Marsden-Weinstein quotient carries a natural symplectic structure induced from
the symplectic structure on $M$, see \cite[Proposition 5.40]{mcduff-salamon1}.

We denote by 
$$G_{L_j}:=\{g \in G:gL_j=L_j\}$$
for $j \in \{0,1\}$ the isotropy subgroup of the Lagrangian $L_j$.
It follows directly from the definitions that 
$G_{S_j}:=\{g \in G: S_j g=g\}$ is a subgroup of $G_{L_j}$. If 
$\mu^{-1}(0) \cap L_j \neq \emptyset$, than the two groups agree. 
To see that, note that by (H1) there
exists $p \in L_j$ whose isotropy subgroup is trivial, i.e.
$G_p:=\{g \in G :gp=p\}=\{\mathrm{id}\}.$
If $g \in G_{L_j}$ then 
$g^{-1} S_j(g) p=p$
and hence $g \in G_{S_j}$. 

We denote by $\mathfrak{g}_{L_j}$
the Lie-algebra of $G_{L_j}$. Note that if $\mu^{-1}(0) \cap L_j \neq \emptyset$ 
$$\mathfrak{g}_{L_j}=\{\xi \in \mathfrak{g}: \dot{S}_j(\xi)=
\xi\}, \quad \mathfrak{g}_{L_j}^\perp=
\{\xi \in \mathfrak{g}: \dot{S}_j(\xi)=-\xi\},$$
where $\perp$ stands for the invariant inner product defined above. The
following proposition says that the two Lagrangians induce Lagrangian submanifolds
in the Marsden-Weinstein quotient. 
\begin{prop} \label{indlag}
Assume $(H1)$ and $(H2)$. Then the subsets of $\bar{M}$
$$\bar{L}_j:=G(L_j \cap \mu^{-1}(0))/G$$
are Lagrangian submanifolds of $\bar{M}$ and they are naturally diffeomorphic to 
$$(L_j \cap \mu^{-1}(0))/G_{L_j}.$$
\end{prop}
The following example shows that Proposition~\ref{indlag} will in general be wrong
if we do not assume hypothesis $(H1)$.

\begin{ex}
Consider the standard action of $S^1$ on $\mathbb{C}^2$ given by
$$(z_1,z_2) \mapsto (e^{i\theta}z_1 e^{i\theta} z_2).$$
A moment map for this action is given by
$$\mu(z):=\frac{i}{2}(|z|^2-1)$$
and the Marsden-Weinstein quotient is the two sphere. Consider now tthe family of
Lagrangian submanifolds of $\mathbb{C}^2$ given by 
$$L_a:=\{(x_1+ia,x_2): x_1,x_2 \in \mathbb{R}\}$$
where $a \in \mathbb{R}$. Note that if $a=0$ then $L_a$ equals the fixpoint
set of the antisymplectic involution on $\mathbb{C}^2$ which is given by
complex conjugation. This involution commutes with the $S^1$-action. For $a \neq 0$
the Lagrangians $L_a$ do not satisfy hypothesis $(H1)$. Consider the chart 
$$\bigg\{\frac{z_2}{z_1}:|z_1|^2+|z_2|^2=1,\,\,z_2 \neq 0\bigg\} \cong \mathbb{C}$$
of $\mu^{-1}/S^1 \cong S^2$. Then the images of $\bar{L}_a$ are given by
$$\bigg\{\frac{t \sqrt{1-t^2-a^2}}{\sqrt{1-t^2}} \pm\frac{iat}{\sqrt{1-t^2}}:
|t| \leq 1-a^2\bigg\}.$$
For $a \neq 0$ these are figure eights with nodal point $(0,0)$.
\end{ex}
To prove Proposition~\ref{indlag} we need two lemmas.

\begin{lemma}\label{clean} Assume $(H1)$ and $(H2)$.
The Lagrangians $L_j$ intersect cleanly with $\mu^{-1}(0)$, i.e.
$\mu^{-1}(0) \cap L_j$ is a submanifold of $\mu^{-1}(0)$ and for every
$p \in \mu^{-1}(0) \cap L_j$ we have 
$T_p \mu^{-1}(0) \cap T_p L_j=T_p( \mu^{-1}(0) \cap L_j)$.
\end{lemma}
\textbf{Proof: }We may assume without loss
of generality that $\mu^{-1}(0) \cap L_j \neq \emptyset$.
For $p \in \mu^{-1}(0) \cap L_j$ we claim that
\begin{equation}\label{tracl}
d \mu(p) T_p L_j=\mathfrak{g}_{L_j}^\perp.
\end{equation}
We first calculate for $v \in  T_p L_j$ and $\xi \in \mathfrak{g}_{L_j}$
\begin{eqnarray*}
\langle d \mu(p) v,\xi \rangle=d\langle \mu,\xi \rangle(p) v
=\omega(X_\xi(p),v)=0
\end{eqnarray*}
and hence
$$d \mu(p) T_pL_j \subset \mathfrak{g}^\perp_{L_j}.$$
To prove that equality holds in (\ref{tracl}) it suffices to show the
following implication
\begin{equation}\label{tracl2}
\xi \in \mathfrak{g}_{L_j}^\perp,\,\,\langle \xi, d\mu(p)v \rangle=0\,\,
\forall v \in T_p L_j \quad \Longrightarrow \quad \xi=0.
\end{equation}
Assume that $\xi$ satisfies the assumption in (\ref{tracl2}). 
Let $w \in T_p M$. Then $w=w_1+w_2$, where
$w_1 \in T_p L_j$, i.e. $dR_j w_1=w_1$, and $dR_j w_2=-w_2$. We calculate
\begin{eqnarray*}
\langle \xi,d\mu(p) w \rangle &=& \langle \xi d \mu(p) w_2 \rangle\\
&=& \omega(X_\xi(p),w_2)\\
&=& \omega(X_\xi(p),-dR_j w_2)\\
&=& \omega(dR_j X_\xi(p),w_2)\\
&=& \omega(X_{\dot{S}_j \xi}(p),w_2)\\
&=&-\omega(X_\xi(p),w_2)\\
&=&-\langle \xi, d\mu(p) w \rangle
\end{eqnarray*}
and hence
$$\langle \xi, d \mu(p) w \rangle=0 \,\, \forall w \in T_p M.$$
Since $0=\mu(p)$ is a regular value of $\mu$ it follows that
$d\mu(p)$ is surjective and hence
$\xi=0$. This proves (\ref{tracl2}) and hence (\ref{tracl}). If $U$ is
a sufficiently small open neighbourhood of $0$ in $\mathfrak{g}_{L_j}^\perp$, then
$\mu^{-1}(U)$ is a submanifold of $M$. 
It follows from
(\ref{tracl}), that $\mu^{-1}(0)$ and $L_j \cap \mu^{-1}(U)$ intersect transversally
in $\mu^{-1}(U)$. Hence $\mu^{-1}(0)$ and
$L_j$ intersect cleanly. \hfill $\square$

\begin{lemma}\label{base}
Assume $(H1)$ and $(H2)$. 
If $p \in L_j$, $G_p=\{g \in G: gp=p\}=\{\mathrm{id}\}$ and 
$\psi(g)p \in L_j$ for some $g \in G$, then $g \in G_{L_j}.$
\end{lemma}
\textbf{Proof: } Because $L_j=\mathrm{Fix}(R_j)$,
$$R_jgR_jp=gp.$$
Since $R_jG=GR_j$ there exists $\tilde{g} \in G$ such that
$$\tilde{g}x=R_jgR_jx \quad \forall \,\,x \in M.$$
Hence 
$$(\tilde{g})^{-1} g p=p$$
and because $G_p=\{\mathrm{id}\}$
$$\tilde{g}=g.$$
Hence $g=R_jgR_j$ and for every $q \in L_j$
$$gq=R_jgq.$$ 
This implies that $gq \in L_j$ and hence $g \in G_{L_j}$. \hfill $\square$
\\ \\
\textbf{Proof of Proposition~\ref{indlag}:} It follows from Lemma~\ref{clean}
and the fact that $G_{L_j}$ acts freely on $L_j \cap \mu^{-1}(0)$ that
$(L_j \cap \mu^{-1}(0))/G_{L_j}$ is a manifold. There is an obvious 
surjective map from
$(L_j \cap \mu^{-1}(0))/G_{L_j}$ to $\bar{L}_j$ which assigns to a representative 
$x \in L_j \cap \mu^{-1}(0)$ of an equivalence class in
$(L_j \cap \mu^{-1}(0))/G_{L_j}$ the equivalence class of $x$ in $\bar{L}_j$. 
It follows from Lemma~\ref{base} that this map is an injection. \hfill $\square$
\\ \\ 
In addition we make the following topological assumptions.
\begin{description}
 \item[(H3)] \emph{$\pi_2(M)$, $\pi_1(M)$, $\pi_1(L_j)$, and 
  $\pi_0(L_j)$
  for $j \in \{0,1\}$ are trivial.}\footnote{Most of the results of this
  paper could be generalized to the case, where we replace (H3) by the following
  weaker assumption $(H3')$\\ \\
  \textbf{(H3')}\,\, 
   \emph{For every smooth map $v:(B,\partial B) \to (M,L_j)$, we have
    $$\int_{B} v^* \omega=0.$$}
  However, if we only assume (H3') instead of (H3), then our path space
  will in general neither be connected nor simply connected. In particular,
  there will be no well defined action functional.}
\end{description}
Convex structures for Hamiltonian group actions on symplectic manifolds were
introduced in \cite{cieliebak-mundet-salamon}. We give a similar definition
which takes care of the Lagrangian submanifolds. Recall that an almost complex
structure is called $\omega$-compatible if
$$\langle \cdot, \cdot \rangle=\omega(\cdot, J\cdot )$$
is a Riemannian metric on $TM$. We say that an almost complex structure is
$G$-invariant, if
$$J(z)=g_* J(z):=d\psi(g)^{-1}(gz)J(gz) d\psi(g)z, \quad \forall \,\,g \in G,
\,\,\forall \,\,z \in M.$$
It follows from \cite[Proposition 2.50]{mcduff-salamon1} that the space of
$G$-invariant compatible almost complex structures is nonempty and
contractible. 
\begin{fed} \label{convexstr}
A \emph{\textbf{convex structure}} on $(M,\omega,\mu,L_0,L_1)$
is a pair $(f,J)$ where $J$ is a $G$-invariant
$\omega$-compatible almost complex structure on $M$ 
which satisfies 
\begin{equation}\label{Rinv}
dR_j(R_jz)J(R_jz)dR_j(z)=-J(z), \quad \forall z \in M.
\end{equation}
for $j \in \{0,1\}$ and $f:M \to [0,\infty)$
is a smooth function satisfying the following conditions.
\begin{description}
 \item[(C1)] $f$ is $G$-invariant and proper.
 \item[(C2)] There exists a constant $c_0>0$ such that
  $$f(x) \geq c_0 \quad \Longrightarrow \quad 
  \langle \nabla_\xi \nabla f(x),\xi\rangle \geq 0, \,\,
  df(x)J(x)X_{\mu(x)}(x) \geq 0,\,\, \mu(x) \neq 0$$
  for every $x \in M$ and every $\xi \in T_xM$. Here $\nabla$ denotes the
  Levi-Civita connection of the metric 
  $\langle \cdot, \cdot \rangle=\omega(\cdot,J\cdot)$.
 \item[(C3)] For every $p \in L_j$ it holds that $\nabla f(p) \in T_p L_j$.
\end{description}
\end{fed}
As our fourth hypothesis we assume that a convex structure exists
\begin{description}
 \item[(H4)] \emph{There exists a convex structure $(f,J_0)$ on 
  $(M,\omega,\mu,L_0,L_1)$.}
\end{description}
A convex structure will guarantee, that solutions of our gradient equation
will remain in a compact domain.
\\ 

The main examples we have in mind, are of the following form. The symplectic
manifold $(M,\omega)$ equals the complex vector space $\mathbb{C}^n$ endowed
with its canonical symplectic structure $\omega_0$, the Lagrangians equal some
linear Lagrangian subspace of $\mathbb{C}^n$, and the group action 
$\psi$ is given by some injective linear representation
of a connected compact
Lie group $G$ to $U(n)$. 

For a linear Lagrangian subspace $L$ of $\mathbb{C}^n$
there is a $\mathbb{R}$-linear splitting
$$\mathbb{C}^n=L \oplus J_0 L$$
where $J_0$ is the standard complex structure given by multiplication with 
$i$. Let $R=R_L$ be the canonical antisymplectic involution given by
$$R(x+J_0y)=x-J_0y$$
for $x,y \in L$. Hypothesis (H1) means
$$\rho(G)R=R\rho(G).$$
In the special case where $\rho=\textrm{id}$ and $L=\mathbb{R}^n$, the
induced involution $S$ on $G$ is given by complex conjugation. To see that,
choose $A \in U(n)$ and $z \in \mathbb{C}$. We calculate
$$S(A)z=RAR(z)=RA\bar{z}=\bar{A}z$$
and hence $S(A)=\bar{A}$. 

The Lie algebra 
$u(n)$ of $U(n)$ carries a natural invariant inner product given by
$$\langle A,B \rangle:=\textrm{trace}(A^*B),$$
where $A^*$ is the complex conjugated transposed of $A$.
Let $\dot{\rho}: \mathfrak{g} \to u(n)$ be the induced representation of $\rho$. 
Endow $\mathfrak{g}$ with the invariant inner product 
$$\langle \xi_1,\xi_2  \rangle_\rho=
\langle \dot{\rho}(\xi_1),\dot{\rho}(\xi_2)\rangle, \quad \forall \,\,
\xi_1,\xi_2 \in 
\mathfrak{g}.$$
Let $\dot{\rho}^*: u(n) \to \mathfrak{g}$ be the adjoint of $\dot{\rho}$, i.e.
$$\langle \dot{\rho}(\xi),\eta \rangle=\langle \xi, \dot{\rho}^*(\eta)\rangle_\rho
\quad \forall \,\,\xi \in \mathfrak{g}, \eta \in u(n),$$ 
then a moment map $\mu$ for the action of 
$G$ is given by
$$\mu(z)=-\frac{1}{2} \dot{\rho}^*(i z z^*)-\tau$$
where $\tau$ is a central element of $\mathfrak{g}$.
Since $\dot{\rho}$ is isometric with respect to our inner products, 
$\langle \cdot, \cdot \rangle$ is $\dot{S}$ invariant. 

A convex structure on 
$(\mathbb{C}^n,\omega_0, \mu,L)$ is defined for example by 
$$(f,J)=(\frac{1}{2}|z|^2,J_0).$$

\begin{ex}[Toric manifolds]
Let $A$ be a $k \times n$-matrix of rank $k$ 
whose entries are positive integers.
Let the $k$-torus $T^k$ act on $\mathbb{C}^n$ by
$$z \mapsto \exp(2 \pi i A \theta^T)z$$
where $\theta=(\theta_1, \ldots,\theta_k) \in \mathbb{R} /
\mathbb{Z} \times \ldots \times \mathbb{R} / \mathbb{Z}=T^k$
and $z=(z_1, \ldots ,z_n) \in \mathbb{C}^n$. For some
$\tau \in i\mathbb{R}^k$ which is equal to the Lie algebra
of the torus a moment map for the torus
action above is defined by 
$$\mu(z)=\frac{1}{2i}(A^T A)^{-1}A^T
\left(\begin{array}{c}
|z_1|^2\\
\vdots\\ 
|z_n|^2
\end{array}\right)-\tau.$$
If $T^k$ acts freely on $\mu^{-1}(0)$, 
then the Marsden-Weinstein quotient
$\mathbb{C}^n //T^k=\mu^{-1}(0)/T^k$ is called a toric manifold.
\end{ex} 

\begin{ex}[Grassmannians]
There is a natural action of the unitary group $U(k)$ on
$\mathbb{C}^{n\times k}$, the space of $k$-frames in $\mathbb{C}^n$,
having the moment map
$$\mu(B)=\frac{1}{2i}(B^* B-\textrm{id}).$$
The symplectic quotient is the complex Grassmannian 
manifold $G_\mathbb{C}(n,k)$.
Let $L=\mathbb{R}^{n \times k}$ be the space of real $k$-frames. 
Its isotropy subgroup is the orthogonal group $O(k)$ and the induced 
Lagrangian in the symplectic quotient equals
the real Grassmannian $G_\mathbb{R}(n,k)$.
\end{ex}

\begin{rem}[Naturality]\label{natur}
For $U \in U(n)$ let the representation $\rho_U$ be given by 
$$\rho_U(g)=U\rho(g) U^{-1},\quad \forall g \in G.$$ 
Then a moment map for this action is given by
$$\mu_U(z)=\mu(U^{-1}z)$$
and there is a natural induced isomorphism from $\bar{M}$ to 
$\bar{M}_U:=\mu_U^{-1}(0)/G$
given by
$$[z] \mapsto [Uz].$$
Defining the linear Lagrangian subspace
$$L_U:=U(L)$$
the image of $\bar{L}$ under the above isomorphism equals
$$\bar{L}_U:=G(\mu^{-1}_U(0) \cap L_U)/G.$$
Because the group $U(n)$ acts transitively on the set of Lagrangian subspaces
of $\mathbb{C}^n$ one can always assume after applying some
$U$ as above, that $L=\mathbb{R}^n$. 
\end{rem}

\subsection[The symplectic vortex equations on the strip]
{The symplectic vortex equations on the strip}

In this subsection we show how one can derive the symplectic vortex
equations from an action functional.

We define the path space $\mathscr{P}$ by 
\begin{equation}\label{pathspa}
\mathscr{P}:=\{(x,\eta) \in C^\infty([0,1],M
\times \mathfrak{g}): x(j) \in L_j,\,\, \eta(j) \in \mathfrak{g}_{L_j}^\perp,
\,\, j \in \{0,1\}\},
\end{equation}
The assumption (H3) implies that $\mathscr{P}$ is connected and
simply connected.
The gauge group $\mathcal{H}$ is defined by
$$\mathcal{H}:=\{g \in C^\infty([0,1],G):g(j) \in G_{L_j},\,\,
g(j)^{-1}\partial_t g(j) \in \mathfrak{g}_{L_j}^\perp,\,\, j \in \{0,1\}\}.$$
The group structure is the pointwise multiplication of $G$.
The gauge group $\mathcal{H}$  acts on $\mathscr{P}$ as follows
$$g_*(x,\eta)=(g x,g \eta g^{-1}-\partial_t g g^{-1}),
\quad g \in \mathcal{H}.$$ 
Choose a path $x_0: [0,1] \to M$ with $x_0(j) \in L_j$ for $j \in \{0,1\}$. 
For a smooth family of $G$-invariant functions $H_t:M \to \mathbb{R}$ for
$t \in [0,1]$, we define the action functional 
$$\mathcal{A}_{\mu,H}:\mathscr{P} \to \mathbb{R}$$
by
$$\mathcal{A}_{\mu,H}(x,\eta)=-\int_{[0,1]\times [0,1]}\bar{x}^*\omega
+\int_0^1\big(\langle \mu(x(t)),\eta(t)\rangle-H_t(x(t))\big)dt,$$
where $\bar{x}:[0,1] \times [0,1] \to M$ is a smooth map, which satisfies 
$$\bar{x}(t,1)=x(t),\quad \bar{x}(t,0)=x_0(t),\quad \bar{x}(0,s) \in L_0,\quad
\bar{x}(1,s) \in L_1.$$
Since $\omega$ is closed and vanishes on the Lagrangians the assumption that
$\pi_2(M)=0$ and $\pi_1(L_j)=0$ for $j \in \{0,1\}$
together with Stokes theorem implies that 
the value of $\mathcal{A}_{\mu,H}(x,\eta)$
does not depend on the choice of $\bar{x}$. Moreover,
$\mathcal{A}_{\mu,H}$ is invariant under the action of
$\mathcal{H}_0$, the path-connected component of the identity of 
$\mathcal{H}$. To see this, let $g \in \mathcal{H}_0$, then there exists
$h: [0,1]\times [0,1] \to G$ with 
$$h(t,1)=g(t), \quad h(t,0)=\mathrm{id}, \quad h(j,s) \in G_{L_j},
\quad h^{-1}\partial_t h(j,s) \in \mathfrak{g}_{L_j}^\perp,\,\,j \in \{0,1\}.$$
The claim follows with $\overline{g^*x}=h \bar{x}$. 

The tangent space $T_{(x,\eta)}\mathscr{P}$ 
of the path space $\mathscr{P}$ 
at $(x,\eta) \in \mathscr{P}$ is defined as the vector space
$$\{(\hat{x},\hat{\eta}) \in 
C^\infty([0,1],x^*TM \times \mathfrak{g}):
\hat{x}(j) \in T_{x(j)}L_j,\,\,\hat{\eta}(j) \in \mathfrak{g}^\perp_{L_j},\,\,
j \in \{0,1\}\}.$$
A family of $G$-invariant, $\omega$-compatible, almost complex
structures $J_t$ determines an inner product on 
$\mathscr{P}$ by 
\begin{equation}\label{metric}
\langle (\hat{x}_1,\hat{\eta}_1),(\hat{x}_2,\hat{\eta}_2)\rangle
=\int_0^1(\langle \hat{x}_1(t),\hat{x}_2(t)\rangle_t+
\langle \hat{\eta}_1(t),\hat{\eta}_2(t)\rangle)
dt
\end{equation}
for $(\hat{x}_1,\hat{\eta}_1),(\hat{x}_2,\hat{\eta}_2) \in T_{(x,\eta)}\mathscr{P}$,
where
$$\langle \cdot,\cdot \rangle_t=\langle \cdot, \cdot \rangle_{J_t}
=\omega( \cdot, J_t \cdot).$$
The gradient of $\mathcal{A}_{\mu,H}$ with
respect to the above inner product as usual defined by
$$d \mathcal{A}_{\mu,H}(x,\eta)[\hat{x},\hat{\eta}]=\langle
\mathrm{grad} \mathcal{A}_{\mu,H}(x,\eta),[\hat{x},\hat{\eta}]\rangle$$
is given by
\begin{displaymath}
\textrm{grad}\mathcal{A}_{\mu,H}(x,\eta)=
\left(\begin{array}{c}
J_t(\dot{x}+X_\eta(x)-X_{H_t}(x))\\
\mu(x)
\end{array}\right).
\end{displaymath}
The set $\mathrm{crit}(\mathcal{A}) \subset \mathscr{P}$ of
critical points of $\mathcal{A}_{\mu,H}$ consists of paths 
$(x,\eta):[0,1] \to M \times \mathfrak{g}$ which satisfy
$$\dot{x}+X_\eta(x)=X_{H_t}(x), \quad \mu(x)=0, \quad x(j) \in L_j,
\quad \eta(j) \in \mathfrak{g}_{L_j}^\perp,\,\,j \in \{0,1\}.$$
Since $G$ acts freely on $\mu^{-1}(0)$ the group $\mathcal{H}$ acts freely 
on $\mathrm{crit}(\mathcal{A})$. 
If $\bar{H}$ is the induced Hamiltonian function of $H$ in 
the Marsden-Weinstein quotient $\bar{M}$
and $\phi^t_{\bar{H}}$ its flow, i.e.
$$\frac{d}{dt}\phi^t_{\bar{H}}=X_{\bar{H}_t} \circ \phi^t_{\bar{H}},
\quad \phi^0_{\bar{H}}=\mathrm{id},$$
then we will prove in Lemma~\ref{zehnder} below that there is a natural bijection
$$\mathrm{crit}(\mathcal{A})/\mathcal{H} \cong
\phi^1_{\bar{H}} (\bar{L}_0) \cap \bar{L}_1.$$
Let
$$\Theta=\{z=s+it \in \mathbb{C}:0 \leq t \leq 1\}$$
be the strip. The flow lines of the vector field
$\textrm{grad}\mathcal{A}_{\mu,H}$ are pairs 
$(u,\Psi) \in C^\infty_{loc}(\Theta,M \times \mathfrak{g})$, 
which satisfy the following partial differential equation
\begin{equation}\label{eq1}
\begin{array}{c}
\partial_s u+J_t(u)(\partial_tu+X_\Psi(u)-X_{H_t}(u))=0\\
\partial_s \Psi+\mu(u)=0\\
u(s,j)\in L_j, \quad\eta(s,j) \in \mathfrak{g}_{L_j}^\perp,\,\,\, j\in\{0,1\}.
\end{array}
\end{equation}
We define further the gauge group
$$\mathcal{G}_{loc}=
\{g \in C^\infty_{loc}(\Theta,G):g(s,j) \in G_{L_j},\,
g^{-1}\partial_t g(s,j) \in \mathfrak{g}_{L_j}^\perp, \,
j \in \{0,1\}\}.$$
Solutions of the problem (\ref{eq1}) are invariant under the action of 
$\mathcal{H}$ but not of $\mathcal{G}_{loc}$. 
To make the problem invariant under the gauge
group $\mathcal{G}_{loc}$, we introduce an additional variable $\Phi$.  
Given a solution
$(u_0,\Psi_0)$ of (\ref{eq1}) and $g \in \mathcal{G}_{loc}$ then
$(u,\Psi,\Phi)=(g u_0,g \Psi_0 g^{-1}-g^{-1}\partial_t g,-g^{-1}\partial_s g)$ 
is a solution of the 
so called \textbf{symplectic vortex equations on the strip}
\begin{equation}\label{eq}
\begin{array}{c}
\partial_s u+X_\Phi(u)+J_t(u)(\partial_tu+X_\Psi(u)-X_{H_t}(u))=0\\
\partial_s \Psi-\partial_t\Phi+[\Phi,\Psi]+\mu(u)=0\\
u(s,j)\in L_j, \quad \Phi(s,j)\in \mathfrak{g}_{L_j}, 
\quad \Psi(s,j) \in \mathfrak{g}_{L_j}^\perp \quad j\in\{0,1\}.
\end{array}
\end{equation}
Moreover, (\ref{eq}) is invariant under the action
of $g \in \mathcal{G}_{loc}$ given by
$$g_*(u,\Psi,\Phi)=(g u,
g \Psi g^{-1}-\partial_t g g^{-1},g \Phi g^{-1}-\partial_s g g^{-1}).$$
On the other hand each solution of (\ref{eq}) is gauge 
equivalent to a solution of (\ref{eq1}). To see that, 
let $(u,\Psi,\Phi)$ be a solution of (\ref{eq}) and
take the solution $g:\Theta \to G$ of the following ordinary differential equation
on the strip
$$\partial_s g=g \Phi, \quad g(0,t)=\textrm{id}.$$
Then $g \in \mathcal{G}_{loc}$ and $g_*\Phi=0$. In the
terminology of gauge theory, this means that solutions of (\ref{eq1})
are solutions of (\ref{eq}) in so called \textbf{radial gauge}.  
\begin{rem}\label{connection}
If one introduces the connection $A=\Phi ds+\Psi dt$ on the trivial 
$G$-bundle over the strip, then the first two equations of (\ref{eq})
can be written as
$$\bar{\partial}_{J,H,A}(u)=0, \qquad *F_A+\mu(u)=0.$$
These equations were discovered independently by D.Salamon and I.Mundet
(see \cite{cieliebak-gaio-salamon}, \cite{cieliebak-mundet-salamon}, and 
\cite{mundet}). In the physics literature they are known as
gauged sigma models.
\end{rem}

\begin{rem}[Naturality]\label{naturality}
Solutions of the problem (\ref{eq}) have the following properties. Let
$K_t$ be some smooth family of $G$-invariant functions on $M$ and
let $\psi_K^t:M \to M$ be the Hamiltonian symplectomorphism defined
by
$$\frac{d}{dt}\psi^t_K=X_{K_t} \circ \psi^t_K, \quad \psi_K^0=\mathrm{id}.$$
If $(u,\Psi,\Phi)$ is a solution of (\ref{eq}), then
$$(\tilde{u},\tilde{\Psi},\tilde{\Phi})(s,t):=(\psi_K^{-t} \circ u,
\Psi,\Phi)(s,t)$$
is also a solution of (\ref{eq}) with $H,J,L_0,L_1$ replaced by
$$\tilde{H}_t:=(H_t-K_t) \circ \psi^t_K, \quad
\tilde{J}_t:=(\psi_K^t)^*J_t, \quad \tilde{L}_0=L_0, \quad
\tilde{L}_1=\psi^{-1}_K L_1.$$
In particular, by choosing $H_t=K_t$ one can always assume
that $H \equiv 0$.
\end{rem}

\begin{prop}\label{zehnder}
There is a natural bijection between $\mathrm{crit}(\mathcal{A})/\mathcal{H}$
and $\phi^1_{\bar{H}}(\bar{L}_0) \cap \bar{L}_1$. 
\end{prop}
\textbf{Proof: } By Remark~\ref{naturality} we may assume without loss of
generality that $H=0$. Denote by $\pi$ the canonical projection from
$\mu^{-1}(0)$ to $\bar{M}=\mu^{-1}(0)/G$. If $q \in \bar{L}_0 \cap \bar{L}_1$,
then there exists $x_0 \in L_0$, $x_1 \in L_1$, and $h \in G$ such that
$$\pi(x_0)=\pi(x_1)=q, \quad x_1=hx_0.$$
Choose a smooth path $g \in C^\infty([0,1],G)$ such that
$$g(0)=\mathrm{id}, \quad g(1)=h, \quad (\partial_t g) g^{-1}(0) \in 
\mathfrak{g}_{L_0}^\perp, \quad (\partial_t g) g^{-1}(1) 
\in \mathfrak{g}_{L_1}^\perp.$$
Such a path exists, since $G$ is connected. Now define
$$I: \bar{L}_0 \cap \bar{L}_1 \to \mathrm{crit}(\mathcal{A})/\mathcal{H},
\quad q \mapsto [(g(t)x_0,-(\partial_t g) g^{-1}(t))].$$
Here $[ \cdot,\cdot]$ denotes the equivalence class in 
$\mathrm{crit}(\mathcal{A})/\mathcal{H}$. 
\\
We have to show that $I$ is well defined. To see that choose another
quadruple $(\tilde{x}_0,\tilde{x}_1,\tilde{h},\tilde{g})$ which satisfies
the relations above. It follows from Lemma~\ref{base} that there
exists $h_j \in G_{L_j}$ for $j \in \{0,1\}$
such that
$$\tilde{x}_j=h_j x_j.$$
Define 
$$\gamma(t):=g(t) h_0^{-1}\tilde{g}(t)^{-1}.$$
Then 
$$\gamma_*(g(t)x_0,-(\partial_t g) g^{-1})=
(\tilde{g}(t)\tilde{x}_0,-(\partial_t \tilde{g}) \tilde{g}^{-1})$$
and
$$\gamma \in \mathcal{H}.$$
This shows that $I$ is well defined. The verification that $I$ is a bijection is
easy, namely to construct the inverse of $I$ 
map $(x,\eta) \in \mathrm{crit}(\mathcal{A})$
to $\pi(x)(0)$. \hfill $\square$

\begin{rem}[Extension]\label{extens}
Every solution can be extended to the whole complex plane. To see that
let $(u,\Psi,\Phi)$ be a solution of (\ref{eq}) and assume that
\begin{equation}\label{R-inv}
J_j(z)=-dR_j(R_jz) J_j(R_jz)dR_j(z), \quad z \in M,\,\, j \in \{0,1\}.
\end{equation}
For simplicity, assume also that $H \equiv 0$. 
Let $\hat{J}_t$
for $t \in \mathbb{R}$ be the unique $G$-invariant extension of
$\omega$-compatible almost complex structures on $M$ defined by the
following conditions
\begin{eqnarray*}
\hat{J}|_{[0,1] \times M}&=&J\\
\hat{J}_{2n-t}(z)&=&-dR_0(R_0 z) \hat{J}_{2n+t}(R_0 z) dR_0(z),  \quad
n \in \mathbb{Z},\,\,\,t \in (0,1]\\
\hat{J}_{2n+1-t}(z)&=&-dR_1(R_1 z) \hat{J}_{2n+1+t}(R_1 z) dR_1(z), \quad
n \in \mathbb{Z},\,\,\,t \in (0,1].
\end{eqnarray*}
For $n \in \mathbb{Z}$ and $t \in (0,1]$ let 
$(\hat{u}, \hat{\Psi},\hat{\Phi}) \in W^{1,p}_{loc}(\mathbb{C},M \times
\mathfrak{g} \times \mathfrak{g})$ be defined by the conditions
\begin{eqnarray*}
(\hat{u},\hat{\Psi},\hat{\Phi})|_{\Theta}&=&(u,\Psi,\Phi),\\
(\hat{u},\hat{\Psi},\hat{\Phi})(s,2n-t)&=&
(R_0 \hat{u},-\dot{S}_0(\hat{\Psi}), \dot{S}_0(\hat{\Phi}))(s,2n+t),\\
(\hat{u},\hat{\Psi},\hat{\Phi})(s,2n+1-t)&=&
(R_1 \hat{u},-\dot{S}_1(\hat{\Psi}), \dot{S}_1(\hat{\Phi}))(s,2n+1+t).
\end{eqnarray*}
Here $S_j$ for $j \in \{0,1\}$ was defined in (\ref{S}).
The map $(\hat{u},\hat{\Psi},\hat{\Phi})$ solves
\begin{equation}\label{eqext}
\begin{array}{c}
\partial_s \hat{u}+X_{\hat{\Phi}}(\hat{u})+J_t(\hat{u})
(\partial_t\hat{u}+X_{\hat{\Psi}}(\hat{u}))=0\\
\partial_s \hat{\Psi}-\partial_t\hat{\Phi}+
[\hat{\Phi},\hat{\Psi}]+\mu(\hat{u})=0\\
(\hat{u},\hat{\Psi},\hat{\Phi})(s+2,t)=[(R_0) \circ (R_1) \hat{u},
(\dot{S}_0) \circ (\dot{S}_1)\hat{\Psi}, (\dot{S}_0) \circ 
(\dot{S}_1) \hat{\Phi})(s,t).
\end{array}
\end{equation}
Solutions of (\ref{eqext}) are invariant under the action of the
gauge group
$$\hat{\mathcal{G}}_{loc}:=\{g \in C^\infty_{loc}(\mathbb{C},G):
g(s,t+2)=(S_0) \circ (S_1) g(s,t)\}.$$
\end{rem}

\subsection[Compactness]{Compactness}

The energy of a solution of (\ref{eq}) is defined by
$$E(u,\Psi,\Phi):=\int_\Theta 
 \bigg(|\partial_s u+X_\Phi(u)|^2
+|\mu(u)|^2 \bigg)dsdt.$$
The aim of this subsection is to prove that every sequence of finite energy
solutions of (\ref{eq}) has a convergent subsequence modulo gauge 
invariance. The main ingredient in the proof is Uhlenbeck's compactness
theorem, which states that a connection with an
$L^p$-bound on the curvature is gauge equivalent to a connection
which satisfies an $L^p$-bound on all its first derivatives. 

Compactness fails if solutions of (\ref{eq}) can escape to
infinity. To make sure that this cannot happen we have to choose our
almost complex structure and the Hamiltonian function appropriately. 
Fix some convex structure $\mathcal{K}=(f,\tilde{J})$ on 
$(M,\omega,\mu,L_1,L_2)$.
Let 
$\mathcal{J}(M,\omega,\mu,\mathcal{K})$ be
the space of all $G$-invariant $\omega$-compatible almost complex 
structures 
$J$ on $(M,\omega)$ which 
equal $\tilde{J}$ outside of a compact set in $M$. It is proven in 
Proposition 2.50 in \cite{mcduff-salamon1} that the space
$\mathcal{J}(M,\omega,\mu,\mathcal{K})$ is nonempty and contractible. 
We define the \textbf{space of admissible families of almost complex
structures}
$$\mathcal{J}:=\mathcal{J}([0,1],M,\omega,\mu,\mathcal{K})
\subset C^\infty([0,1],\mathcal{J}(M,\omega,\mu,\mathcal{K}))$$
as the space consisting of smooth families of 
$J_t \in \mathcal{J}(M,\omega,\mu,\mathcal{K})$ which satisfy
(\ref{R-inv}).
Let $C^\infty_{0,G}(M)$ be the space of smooth 
$G$-invariant functions on 
$M$ with compact support, and
$$\mathrm{Ham}:=\mathrm{Ham}(M,G):=
\{H \in C_0^\infty([0,1]\times M): H_t \in C^\infty_{0,G}(M)\}$$
the space of $G$-invariant functions parametrised by $t \in [0,1]$. 

\begin{thm}[Compactness]\label{compactn}
Let $(u_\nu,\Psi_\nu,\Phi_\nu) \in
C^\infty_{loc}(\Theta,M\times \mathfrak{g} \times \mathfrak{g})$  be a
sequence of solutions of (\ref{eq}) with respect to a smooth family
of  almost complex
structures $J_t \in \mathcal{J}$ and to a smooth family of Hamiltonian 
functions
$H_t \in \mathrm{Ham}$. If the energies are uniformly bounded,
then there exists a sequence of gauge transformations 
$g_\nu \in \mathcal{G}_{loc}$ such that a subsequence of
$(g_\nu)_*(u_\nu,\Psi_\nu,\Phi_\nu)$ converges in the 
$C^{\infty}_{loc}$-topology to a smooth solution 
$(u,\Psi,\Phi)$ of the vortex problem(\ref{eq}). 
\end{thm}
Instead of Theorem~\ref{compactn} we prove the following stronger theorem.
\begin{thm}\label{uh}
Let $(u_\nu,\Psi_\nu) 
\in C^\infty_{loc}(\Theta,M\times\mathfrak{g})$
be a sequence of solutions of (\ref{eq1}) with respect to a smooth family of
almost complex structures $J_t \in \mathcal{J}$ and to a smooth family of
Hamiltonian functions $H_t \in \mathrm{Ham}$. If the energies are 
uniformly bounded,
then there exists a sequence of gauge transformations
$g_\nu \in \mathcal{H}$ such that a subsequence of
$(g_\nu)_*(u_\nu,\Psi_\nu)$ converges in the
$C_{loc}^{\infty}$-topology to a smooth
solution $(u,\Psi)$ of the gradient equation (\ref{eq1}).
\end{thm}
\textbf{Proof: } Let 
$(\hat{u}_\nu,\hat{\Psi}_\nu) 
\in C^\infty_{loc}(\mathbb{C},M \times \mathfrak{g})$ be the extension of
$(u_\nu,\Psi_\nu)$ as in Remark~\ref{extens}. We will prove the
theorem in four steps.
\\ \\
\textbf{Step 1:} \emph{For every compact subset $K$ of $\mathbb{C}$
there exists a sequence of gauge transformations 
$g_\nu \in C^\infty(K,G)$ such that 
$(g_\nu)_*(\hat{u}_\nu,\hat{\Psi}_\nu,0)|_K$
converges in the $C^\infty$-topology.}
\\ \\
Step 1 was proved in \cite[Theorem 3.4.]{cieliebak-mundet-salamon}. The main
ideas are the following. 
Let $\hat{A}_\nu=\hat{\Psi}_\nu dt$ be the connection on the trivial
$G$-bundle over $\mathbb{C}$. By convexity
$\hat{u}_\nu(K)$ is contained in a compact subset of $M$. The curvature
of $\hat{A}_\nu$ is given by 
$F_{\hat{A}_\nu}=\partial_s \hat{\Psi}_\nu ds \wedge dt$.
Hence the equation
$\partial_s \hat{\Psi}_\nu+\mu(\hat{u}_\nu)=0$ 
implies that the curvature is uniformly bounded.
Moreover, because of (H4) there is no bubbling and hence
$$\sup_\nu||\partial_s \hat{u}_\nu|_K||_\infty<\infty.$$
Step 1 follows now from a combination of Uhlenbecks compactness theorem, 
see \cite{uhlenbeck, wehrheim}, and the Compactness Theorem for the
Cauchy-Riemann operator (see for example \cite[Appendix B]{mcduff-salamon2}).
\\ \\
\textbf{Step 2:} \emph{There exists a sequence of gauge transformations
$g_\nu \in C^\infty_{loc}(\mathbb{C},G)$ such that a subsequence of 
$(g_\nu)_*(\hat{u}_\nu,\hat{\Psi}_\nu,0)$ converges in the 
$C^\infty_{loc}$-topology.}
\\ \\
We use the fact that $\mathbb{C}$ can be exhausted by compact sets,
$\mathbb{C}=\bigcup_{n \in \mathbb{N}}B_n$ where
$B_n:=\{z \in \mathbb{C}:|z|=n\}$.
Let $g_\nu^n$ be the sequence of gauge transformations on $B_n$
for $n \in \mathbb{N}$ obtained in step 1. We show that we may
assume that
\begin{equation}\label{induc}
g^n_\nu|_{B_{n-1}}=g^{n+1}_\nu|_{B_{n-1}}, 
\quad \forall \,\,\nu \in \mathbb{N}.
\end{equation}
To prove (\ref{induc}) we first observe that there exists
$h^n \in C^\infty(B_n,G)$ such that
$h^n_\nu:=(g_\nu^{n+1})^{-1} \circ g_\nu^n|_{B_n}$ has a subsequence  which
converges in the $C^\infty$-topology to $h^n$. Hence there exists a 
sequence $\tilde{h}^n_\nu$ which has a converging subsequence and satisfies
$$\tilde{h}^n_\nu(z):=\left\{\begin{array}{cc}
h^n_\nu(z) & z \in B_{n-1}\\
\mathrm{id} & z \in B_{n+1} \setminus B_n. \end{array}\right.$$
Now replace $g^{n+1}_\nu$ by $g^{n+1}_\nu \circ \tilde{h}^n_\nu$, which
satisfies (\ref{induc}). Now define $g_\nu$ by
$$g_\nu|_{B_n}:=g_\nu^{n+1}|_{B_n}, \quad \forall \,\,n \in \mathbb{N}.$$
\textbf{Step 3:} \emph{The sequence of gauge transformations
$g_\nu(s+it) \in C^\infty_{loc}(\mathbb{C},G)$ in
step 2 may be chosen independent of $s$.} 
\\ \\
We use an idea from \cite{jones-rawnsley-salamon}. Denote the
limit of the sequence
$(g_\nu)_*(\hat{u}_\nu,\hat{\Psi}_\nu,0)$ as $\nu$ goes to infinity by
$(\hat{u},\hat{\Psi},\hat{\Phi})$.
Choose 
$h\in C^{\infty}_{loc}(\mathbb{C},G)$ such that 
$$h_*\hat{\Phi}=0.$$
Observe that
$$\lim_{\nu \to \infty}(\partial_s(h \circ g_\nu))(h \circ g_\nu)^{-1})
=\lim_{\nu \to \infty}(h \circ g_\nu)_*(0)=h_*\hat{\Phi}=0.$$ 
It follows that  $\partial_s (h \circ g_\nu)$
converges to zero in the $C^{\infty}_{loc}$-topology. Now set
$$\tilde{g}_\nu(s,t):=h \circ g_\nu(0,t).$$
Obviously, $\tilde{g}_\nu$ is independent of $s$. Moreover, since
$\tilde{g}_\nu \circ (h \circ g_\nu)^{-1}$ converges to the identity in
the $C^\infty_{loc}$-topology, $\tilde{g}_\nu$ satisfies the assumptions of
step 2. We denote $\tilde{g}_\nu$ by $g_\nu$ as before. 
\\ \\
\textbf{Step 4:} \emph{We prove the theorem.}
\\ \\
We have to modify further our sequence of gauge transformations $g_\nu$,
such that they satisfy the boundary conditions.
Because the energy of the sequence $(u_\nu,\Psi_\nu)$ is bounded, the energy
of $(\hat{u},\hat{\Psi})|_{\Theta}$ is bounded. Using the fact that
finite energy solutions converge uniformly at the two ends of the strip,
see~\cite{frauenfelder2}, we conclude that  
$\mu(\hat{u})(s,t)$ converges to zero as $s$
goes to infinity uniformly in the $t$-variable. 
Since $G$ acts freely on $\mu^{-1}(0)$, there exists
$s_0 \in \mathbb{R}$ such that
$$G_{u(s_0,0)}=G_{u(s_0,1)}=\{\textrm{id}\}.$$
Since by assumption $\hat{u}_\nu(s_0,j) \in L_j$ for $j \in \{0,1\}$ it
follows that $\hat{u}(s_0,j) \in GL_j$. 
Now choose $h \in C^\infty([0,1],G)$ such that
$h(0)\hat{u}(s_0,0) \in L_0$ and $h(1)\hat{u}(s_0,1) \in L_1$. 
This is
possible, because $G$ is connected. Choose a sequence 
$h_\nu \in C^\infty([0,1],G)$ converging
to the identity satisfying
$$(h_\nu(j) \circ h(j) \circ g_\nu(j))u_\nu(s_0,j) 
\in L_j , \quad j \in \{0,1\}.$$
We can assume without loss of generality that 
$$G_{u_\nu(s_0,0)}=G_{u_\nu(s_0,1)}=\{\textrm{id}\}$$
for every $\nu$. By Lemma~\ref{base}, 
$h_\nu(j) \circ h(j) \circ g_\nu(j) \in G_{L_j}$ for 
$j \in \{0,1\}$. In particular,
$$(h_\nu \circ h \circ g_\nu)u_\nu(s,j) 
\in L_j, \quad j \in \{0,1\}.$$
Hence $h_*\hat{u}(s,j) \in L_j$ for $j \in \{0,1\}$. By a similar
procedure as above we may find a further 
sequence of gauge transformations
$\tilde{h}_\nu \in C^\infty([0,1],G)$ which satisfy the following 
conditions.
\begin{description}
 \item[(i)] $\tilde{h}_\nu(j) \in G_{L_j}$ for $j \in \{0,1\}$, 
 \item[(ii)] $(\tilde{h}_\nu)_* (g_\nu)_* h_*\hat{\Psi}_\nu(s,j) 
  \in \mathfrak{g}_{L_j}^\perp$
  for every $\nu$ and $j \in \{0,1\}$,
 \item[(iii)] There exists $\tilde{h} \in C^\infty([0,1],G)$ such that
  $\tilde{h}_\nu$ converges as $\nu$ goes to infinity in the
  $C^\infty$-topology to $\tilde{h}$.
\end{description} 
Now set
$$\tilde{g}_\nu:=\tilde{h}_\nu \circ h_\nu \circ h \circ g_\nu|_{[0,1]} \in 
C^\infty([0,1],G).$$
Moreover, using the assumption that 
$\hat{\Psi}_\nu(j) \in \mathfrak{g}_{L_j}^\perp$
for $j \in \{0,1\}$, one calculates
$$\tilde{g}_\nu(j) \in G_{L_j}, \quad \tilde{g}_\nu^{-1}\partial_t \tilde{g}_\nu(j)
\in \mathfrak{g}_{L_j}^\perp, \quad j \in \{0,1\},$$
and hence $\tilde{g}_\nu \in\mathcal{H}$. Set
$$(u,\Psi):=\tilde{h}_* h_*(\hat{u},\hat{\Psi})|_\Theta.$$
Now the theorem follows with $g_\nu$ replaced by $\tilde{g}_\nu$.
\hfill $\square$

\subsection[Moduli spaces]{Moduli spaces}

We assume that the Hamiltonian $H \in \mathrm{Ham}$ has the property
that the Lagrangians $\phi^1_{\bar{H}}(\bar{L}_0)$ and $\bar{L}_1$ in the
Marsden-Weinstein quotient intersect transversally. Under this
assumption it can be shown, see \cite{frauenfelder2}, that 
finite energy solutions of
the symplectic vortex equation (\ref{eq}) are gauge equivalent to
solutions which decay exponentially fast at the two ends of the strip.
More precisely, define for some small number $\delta>0$ and some smooth
cutoff function $\beta$ satisfying $\beta(s)=-1$ if $s<0$ and
$\beta(s)=1$ if $s>1$ 
$$\gamma_\delta \in C^\infty(\mathbb{R}), \quad s \mapsto
e^{\delta \beta(s)s}.$$
For an open subset $S \subset \Theta$ we define the
$||\,\,||_{C^k_\delta}$-norm for some smooth function
$f \colon S \to \mathbb{R}$ by
$$||f||_{C^k_\delta}:=||\gamma_\delta \cdot f||_{C^k}$$
and denote
\begin{equation}\label{Cdelta}
C^\infty_\delta(S):=\{f \in C^\infty(S) \colon
||f||_{C^k_\delta}<\infty,\,\,\forall \,\,k \in \mathbb{N}\}.
\end{equation}
We now introduce the Fr\'echet manifold 
$\mathcal{B}=\mathcal{B}_\delta$ as the set consisting of
all $w=(u,\Psi,\Phi) \in 
C^\infty_{loc}(\Theta,M \times \mathfrak{g}) \times 
C^\infty_{\delta}(\Theta,\mathfrak{g})$
which satisfy the following conditions:
\begin{description}
 \item[(i)] $w$ maps $(s,j)$ to $L_j \times \mathfrak{g}_{L_j} \times
  \mathfrak{g}_{L_j}^{\perp}$ for $j \in \{0,1\}$ and $s \in \mathbb{R}$.
 \item[(ii)] There exists a
  critical point of the action functional
  $(x_1,\eta_1) \in \mathrm{crit}(\mathcal{A})$, a real number $T_1 \in
  \mathbb{R}$, and $(\xi_1,\psi_1) \in C^\infty_\delta((-\infty,T_1]
  \times [0,1],x^*_1 TM \times \mathfrak{g})$ such that
  $$(u,\Psi)(s,t)=(\mathrm{exp}_{x_1(t)}(\xi_1(s,t)),\eta_1(t)
  +\psi_1(s,t)), \quad s \leq -T_1.$$
 \item[(iii)] There exists a critical point of the action functional
  $(x_2,\eta_2) \in \mathrm{crit}(\mathcal{A})$, a real number
  $T_2 \in \mathbb{R}$, and $(\xi_2,\psi_2) \in C^\infty_\delta
  ([T_2,\infty) \times [0,1],x_2^*TM\times \mathfrak{g})$ such that
  $$(u,\Psi)(s,t)=(\mathrm{exp}_{x_2(t)}(\xi_2(s,t)),\eta_2(t)+
  \psi_2(s,t)).$$
\end{description}
The theorem about exponential decay proved in \cite{frauenfelder2} now
tells us, that for $\delta>0$ chosen small enough
every finite energy solution of the symplectic vortex
equations is gauge equivalent to an element in $\mathcal{B}_\delta$.
Moreover, one can prove that every solution of (\ref{eq}) which lies
in $\mathcal{B}_{\delta}$ has finite energy.

In a similar vein we define the gauge group \label{banachgauge}
$\mathcal{G}=\mathcal{G}_\delta$ consisting of gauge transformations
$g \in \mathcal{G}_{loc}$ which decay exponentially fast at the two ends
of the strip to elements of $\mathcal{H}_0$ the connected component of
the identity of the gauge group $\mathcal{H}$. Note that there are
natural evaluation maps
$$\mathrm{ev}_1, \mathrm{ev}_2 \colon \frac{\{(\ref{eq})\} \cap 
\mathcal{B}}{\mathcal{G}} \to \pi_0(\mathrm{crit}(\mathcal{A}))
\cong (\phi^1_{\bar{H}}(\bar{L}_0) \cap \bar{L}_1) \times 
\pi_0(\mathcal{H}),$$
induced by the maps $w \mapsto (x_1,\eta_1)$ and $w \mapsto (x_2,\eta_2)$.
One can check that the energy of an element $w \in (\{(\ref{eq})\}\cap
\mathcal{B})/\mathcal{G}$ is given by the difference of the actions
$$E(w)=\mathcal{A}(\mathrm{ev}_1(w))-\mathcal{A}(\mathrm{ev}_2(w))
\footnote{Since the action functional is invariant under the action
of $\mathcal{H}_0$ we denote by abuse of notation the function induced
from the action functional 
on $\pi_0(\mathrm{crit}(\mathcal{A}))$ also by $\mathcal{A}$.}.$$
It can be shown, see \cite{frauenfelder2} that the 
linearization of the symplectic vortex
equations considered as an operator between suitable Banach spaces 
is a Fredholm operator. Moreover, using hypothesis (H3) it follows that
the path space $\mathscr{P}$ defined in (\ref{pathspa}) is connected
and simply connected. This implies that there
exists a function
$$I \colon \pi_0(\mathrm{crit}(\mathcal{A})) \to \mathbb{Z}$$
such that the 
Fredholm index of the
linearized symplectic vortex equations at a point 
$w \in \mathcal{B}/\mathcal{G}$ is given by the difference
$I(\mathrm{ev}_1(w))-I(\mathrm{ev}_2(w))$. 

We can now introduce on $\pi_0(\mathrm{crit}(\mathcal{A}))$ the
following equivalence relation. We say that two connected components
of $\mathrm{crit}(\mathcal{A})$ are equivalent if they project to the
same intersection point of $\phi^1_{\bar{H}}(\bar{L}_0) \cap
\bar{L}_1$ and the action functional $\mathcal{A}$ and the
index $I$ agree on them. We denote the set of such equivalence classes
by $\mathscr{C}=\mathscr{C}(\mathcal{A})$. For 
$c_1, c_2 \in \mathscr{C}$ we define the moduli space
\begin{equation}\label{modulispa}
\tilde{\mathcal{M}}(c_1,c_2):=\{w \in (\{\ref{eq}\}\cap \mathcal{B})
/\mathcal{G}: \mathrm{ev}_1(w) \in c_1,\,\,\mathrm{ev}_2(w) \in c_2\}.
\end{equation}
Note that the moduli spaces depend on the choice of the almost complex
structure $J \in \mathcal{J}$, i.e. 
$\tilde{\mathcal{M}}(c_1,c_2)=\tilde{\mathcal{M}}_J(c_1,c_2)$.
It is proved in \cite{frauenfelder2} that for generic choice of the
almost complex structure the Fredholm operators obtained
by linearizing the symplectic vortex equations are surjectiv. 
In particular, the moduli spaces
$\tilde{\mathcal{M}}_J(c_1,c_2)$ are smooth manifolds whose dimension is
given by the difference of the Fredholm indices of $c_1$ and $c_2$.
\begin{thm} There exists a subset $\mathcal{J}_{reg} \subset \mathcal{J}$
of second category such that $\tilde{\mathcal{M}}_J(c_1,c_2)$ for
every $c_1,c_2 \in \mathscr{C}$ are smooth finite dimensional
manifolds whose dimension is given by
$$\mathrm{dim}(\tilde{\mathcal{M}}_J(c_1,c_2))=I(c_1)-I(c_2).$$
\end{thm}
The group $\mathbb{R}$ acts on $\tilde{\mathcal{M}}(c_1,c_2)$ by
timeshift
We define the path space $\mathscr{P}$ by 
$$w(s,t) \mapsto w(s+r,t), \quad r \in \mathbb{R}.$$
If $c_1 \neq c_2$ then this action is free and the quotient 
$\tilde{\mathcal{M}}(c_1,c_2)/\mathbb{R}$ is again a manifold. Using
the compactness result in Theorem~\ref{compactn} one can show as in
\cite{salamon1} that
the only obstruction to compactness of the spaces
$\tilde{\mathcal{M}}(c_1,c_2)$ is breaking off of flow lines, see
\cite{frauenfelder2}.

\subsection[A Novikov ring]{A Novikov ring}

If $c_1, c_2 \in \pi_0(\mathrm{crit}(\mathcal{A}))$
and $h \in \mathcal{H}$ then
$$\mathcal{A}(c_1)-\mathcal{A}(c_2)=\mathcal{A}(h c_1)-\mathcal{A}(h
c_2), \quad
I(c_1)-I(c_2)=I(hc_1)-I(hc_2)$$
where $I$ is the Fredholm-index introduced in the previous subsection.
It follows that $I(hc)-I(c)$ and $\mathcal{A}(hc)-\mathcal{A}(c)$ is 
independent of the choice of $c \in \pi_0(\mathrm{crit}(\mathcal{A}))$.
Hence we may define the maps
$$I_{\mathcal{H}} \colon \mathcal{H} \to \mathbb{Z}, \quad
E_{\mathcal{H}} \colon \mathcal{H} \to \mathbb{R}, \quad
h \mapsto I(hc)-I(c), \quad h \mapsto \mathcal{A}(hc)-\mathcal{A}(c)$$
for some arbitrary $c \in \pi_0(\mathrm{crit}(\mathcal{A}))$. One easily
checks that these maps are group homomorphisms. Moreover they vanish
on $\mathcal{H}_0$, the connected component of the identity of 
$\mathcal{H}$. For the special case where $M=\mathbb{C}^n$,
$L_0=L_1=\mathbb{R}^n$, and $G$ acts on $\mathbb{C}^n$ by a 
linear injective representation $\rho$ one can show, see
\cite{frauenfelder2}, that the index map $I_{\mathcal{H}}$ is given by
$$I_{\mathcal{H}}(h)=\mathrm{deg}(\mathrm{det}^2_{\mathbb{C}}(\rho(h)))$$
for $h \in \mathcal{H}$.

We define 
$$\Gamma=\frac{\mathcal{H}}
{\ker I_\mathcal{H} \cap \ker E_\mathcal{H}}.$$
To the group $\Gamma$ we associate the Novikov ring
$\Lambda=\Lambda_\Gamma$ whose elements are formal sums
$$r=\sum_{\gamma \in \Gamma} r_\gamma \gamma$$
with coefficients $r_\gamma \in \mathbb{Z}_2$
which satisfy the finiteness condition
$$\#\{\gamma \in \Gamma:r_\gamma \neq 0, 
E_{\mathcal{H}}(\gamma) \geq \kappa\}<\infty$$
for every $\kappa>0$. The multiplication is given by
$$r*s=\sum_{\gamma \in \Gamma}\bigg
(\sum_{\substack{\gamma_1,\gamma_2 \in \Gamma\\
\gamma_1 \circ \gamma_2=\gamma}}r_{\gamma_1} s_{\gamma_2}\bigg)\gamma.$$
Since the coefficients $r_\gamma$ are taken in a field, the Novikov ring
is actually a field. 
The ring comes with a natural grading defined by
$$\deg(\gamma)=I_{\mathcal{H}}(\gamma)$$
and we shall denote by $\Lambda_k$ the elements of degree $k$. Note in
particular that $\Lambda_0$ is a subfield of $\Lambda$. Moreover,
the multiplication maps $\Lambda_j \times \Lambda_k \to \Lambda_{j+k}$.

\subsection[Definition of the homology]{Definition of the homology}

We assume that $H \in \mathrm{Ham}$ has the property that 
$\phi^1_{\bar{H}}(\bar{L}_0)$ and $\bar{L_1}$ intersect 
transversally and $J \in \mathcal{J}_{reg}=\mathcal{J}_{reg}(H)$, i.e.
the Fredholm operators obtained by linearizing the symplectic vortex
equations are surjective. Recall
\begin{equation}\label{critfloer}
\mathscr{C}=\mathscr{C}(\mathcal{A})=
\frac{\textrm{crit}(\mathcal{A})}{\ker I_\mathcal{H} \cap
\ker E_{\mathcal{H}}} 
\cong (\phi^1_{\bar{H}}(\bar{L}_0) \cap \bar{L}_1) \times \Gamma.
\end{equation}
We define the chain complex $CF_*(H,L_0,L_1,\mu)$ as a module over the
Novikov ring $\Lambda$. More precisely, $CF_k(H,L_0,L_1,\mu)$ are 
formal sums of the form
$$\xi=\sum_{\substack{c \in \mathscr{C}\\I(c)=k}}\xi_c c$$
with $\mathbb{Z}_2$-coefficients $\xi_c$ satisfying the 
finiteness condition
\begin{equation}\label{fino}
\#\{c:\xi_c \neq 0,E(c)\geq \kappa\}<\infty
\end{equation}
for every $\kappa>0$. The action of $\Gamma$ on $\mathscr{C}$ is
the induced action of $\mathcal{H}$ on 
$\mathrm{crit}(\mathcal{A})$. The Novikov ring acts on $CF_*$ by
$$r*\xi=\sum_{c \in \mathscr{C}}
\sum_{\substack{c' \in \mathscr{C}, \gamma' \in \Gamma\\
\gamma' c'=c}}\bigg(r_{\gamma'}\xi_{\gamma'}\bigg)c.$$
$CF_k$ is invariant under the action of $\Lambda_0$. In particular,
$CF_k$ which may be an infinite dimensional vector space over
the field $\mathbb{Z}_2$, is a finite dimensional vector space
over the field $\Lambda_0$. 

Recall that for $c_1,c_2 \in \mathscr{C}$ the moduli space 
is defined by
$$\tilde{\mathcal{M}}(c_1,c_2):=\{w \in (\{\ref{eq}\}\cap \mathcal{B})
/\mathcal{G}:
\mathrm{ev}_1(w) \in c_1, \,\mathrm{ev}_2(w) \in c_2\}.$$
Let $\mathcal{T}$ be the group
$$\mathcal{T}:=\frac{\ker I_\mathcal{H} \cap \ker E_\mathcal{H}}
{\mathcal{H}_0},$$
where $\mathcal{H}_0$ is the connected component of the identity
of $\mathcal{H}$. 
Then $\mathcal{T}\times \mathbb{R}$ act on $\tilde{\mathcal{M}}$ by
$$w(s) \mapsto g_*w(s+r), \quad (g,r) \in \mathcal{T} \times \mathbb{R}$$
and we define
$$\mathcal{M}(c_1,c_2):=
\frac{\tilde{\mathcal{M}}(c_1,c_2)}{\mathcal{T}\times \mathbb{R}}.$$
Assume that $c_1 \neq c_2$. Under this assumption 
$\mathcal{T}\times \mathbb{R}$ acts freely on $\tilde{\mathcal{M}}(c_1,c_2)$
and since $J \in \mathcal{J}_{reg}$ the moduli spaces are manifolds 
of dimension
$$\dim \mathcal{M}(c_1,c_2)=\dim \tilde{\mathcal{M}}(c_1,c_2)-1=
I(c_1)-I(c_2)-1.$$ 
Using the fact that the only obstruction to compactness for strips of
finite energy is the breaking off phenomenon, which cannot happen in the
case where the index equals zero, we conclude that for 
$c_1,c_2 \in \mathscr{C}$ with 
$I(c_1)-I(c_2)=1$ and $\kappa>0$ we have
\begin{equation}\label{finite}
\sum_{\substack{\gamma \in \Gamma\\E_\mathcal{H}(\gamma) \geq \kappa, 
\,\,
I_\mathcal{H}(\gamma)=0}} \#\mathcal{M}(c_1,\gamma c_2)<\infty.
\end{equation}
Set
$$n(c_1,c_2):=\#\mathcal{M}(c_1,c_2) \mod 2$$
and define the boundary operator
$\partial_k:CF_k \to CF_{k-1}$ as linear
extension of
$$\partial_k c=\sum_{I(c')=k-1}n(c,c')c'$$
for $c \in \mathscr{C}$ with $I(c)=k$. Note that 
(\ref{finite}) guarantees the
finiteness condition (\ref{fino}) for $\partial_k c$. 

As in the standard theory (see \cite{schwarz1,schwarz2, hofer-salamon} 
one shows that 
$$\partial^2=0.$$
This gives rise to homology groups
$$HF_k(H,J,L_0,L_1,\mu;\Lambda):=
\frac{\ker \partial_{k+1}}{\textrm{im} \partial_k}.$$
A standard argument (see \cite{schwarz1} or \cite{hofer-salamon})
shows that $HF_k(H,J,,L_0,L_1,\mu;\Lambda)$ is actually 
independent of the regular pair $(H,J)$. Hence we set for some regular 
pair $(H,J)$
$$HF_k(L_0,L_1,\mu;\Lambda):=HF_k(H,J,L_0,L_1,\mu;\Lambda).$$
We call the graded $\Lambda$ vector space $HF_*(L_0,L_1,\mu;\Lambda)$
the \textbf{moment Floer homology}.

\subsection[Computation of the homology]
{Computation of the homology}\label{agiv}

In this subsection we compute moment Floer homology for case where
the two Lagrangians coincide, i.e. $L_0=L_1=L$. We will see that
in this case, the moment Floer homology equals the singular homology
of the induced Lagrangian in the Marsden-Weinstein quotient $\bar{L}$
tensored with the Novikov ring introduced above. As a corollary
we get a proof of the Arnold-Givental conjecture for $\bar{L}$.

\begin{thm}\label{arnold} Assume that the two Lagrangians coincide, i.e.
$L_0=L_1=L$, then
$$HF_*(L,\mu;\Lambda)
:=HF_*(L,L,\mu;\Lambda)=H\bar{L}_*(\mathbb{Z}_2)\otimes_{\mathbb{Z}_2}
\Lambda.$$
\end{thm}   
\begin{cor}
The Arnold-Givental conjecture
holds for $\bar{L}$, i.e. under the transversality assumption  
$\bar{L} \pitchfork \phi^1_{\bar{H}} \bar{L}$, 
$$\#(\bar{L} \cap \phi^1_{\bar{H}}\bar{L}) \geq \sum_k 
b_k(\bar{L},\mathbb{Z}_2).$$  
\end{cor}
To prove Theorem~\ref{arnold} we consider the case where the Hamiltonian $H$
vanishes. In this case the Lagrangians in the quotient
$\bar{L}$ and $\phi^1_{\bar{H}}(\bar{L})=\bar{L}$ coincide. In
particular, they do not intersect transversally but still cleanly,
i.e. there intersection is still a manifold whose tangent space
is given by the intersection of the two tangent spaces. This is the 
infinite dimensional analogon of
a Morse-Bott situation. In our case the critical manifold
can be identified with $\bar{L} \times \Gamma$.
Following the approach explained
in the appendix, we still can define the homology in this case
by choosing a Morse function on the critical manifold. To define the
boundary operator one has to count flow lines with cascades. There
is a natural splitting of the boundary operator into two parts. The
first part takes account of the flow lines with zero cascades, i.e.
Morse flow lines on the critical manifold, and the second part takes
account of flow lines with at least one cascade. To prove the theorem
we have to show that the second part of the boundary operator vanishes. 
Using the antisymplectic involution we construct successive involutions
on the cascades which endowes the space of cascades with the structure
of a space of Arnold-Givental type which admits the structure of
a Kuranishi structure of stable Arnold-Givental type. Using this it
follows that the second part of the boundary vanishes.
\\ 

We now define moment Floer homology for the case where the 
Hamiltonian $H=0$.
We think of $\mathcal{H}(\mu^{-1}(0) \cap L)$ as the critical
manifold of the action functional $\mathcal{A}=\mathcal{A}_0$
of the unperturbed symplectic vortex equations. A Morse function on
the induced Lagrangian in the Marsden-Weinstein quotient
$\bar{L}=\mu^{-1}(0)/G_L$ will lift to a $\mathcal{H}$-invariant
Morse function on the critical manifold of $\mathcal{A}$.

We first describe the elements which are needed to define the chain complex.
Choose a Riemannian metric $g$ and a Morse-function $f$ on 
$\bar{L}$ which satisfy the Morse-Smale condition, i.e. stable
and unstable manifolds intersect transversally, 
and lift it to a $G_L$-equivariant metric $\tilde{g}$ and a
$G_L$-equivariant Morse-function 
$\tilde{f}$ on $\mu^{-1}(0) \cap L$. Recall that
for $x \in M$ and $\eta \in \mathfrak{g}$ the linear map
$L_x \colon \mathfrak{g} \to T_x M$ was defined by
$$L_x \eta=X_\eta(x)=\frac{d}{dr}\bigg|_{r=0}\exp(r\eta)(x).$$ 
Let 
$\tilde{\mathscr{C}}_0=\tilde{\mathscr{C}}_0(f)$ be 
the set of smooth maps
$(x,\eta):[0,1] \to M \times \mathfrak{g}$ satisfying
$$\dot{x}(t)+L_{x(t)} \eta(t)=0, \,\mu(x(t))=0,\,t \in [0,1],$$
$$x(j) \in \mathrm{crit}(\tilde{f}), \,
\eta(j) \in \mathfrak{g}_L^\perp, \,j \in \{0,1\}.$$
Note that $\eta$ is completely determined by $x$ through the formula
$$\eta(t)=-(L_{x(t)}^* L_{x(t)})^{-1}L_{x(t)}^* \dot{x}(t),$$
where $L_x^*$ is the adjoint of $L_x$ with respect to the fixed
invariant inner product on $\mathfrak{g}$ and the inner product
$\omega_x(\cdot,J_t(x) \cdot)$ on $T_{x(t)}M$. 
Moreover, it follows from
Proposition~\ref{zehnder} that
there exists an element $g_x$ of the gauge group $\mathcal{H}$ such that
$$x(t)=g_x(t)x(0).$$
Denote
\begin{equation}\label{critmorse}
\mathscr{C}_0:=\mathscr{C}_0(f):=\frac{\tilde{\mathscr{C}}_0}
{\ker I_\mathcal{H} \cap \ker E_\mathcal{H}}.
\end{equation}
Then the map
$$(x,\eta) \mapsto (x(0),g_x)$$
defines a natural bijection
$$\tilde{\mathscr{C}}_0 \cong\mathrm{crit}(\tilde{f}) \times \mathcal{H}$$
and induces a bijection in the quotient
$$\mathscr{C}_0 \cong \frac{\mathrm{crit}(\tilde{f})}{G_L}
\times \Gamma 
\cong \mathrm{crit}(f) \times\Gamma.$$
If $\mathrm{ind}_f$ is the Morse-index, then the  
index of a critical point $c=(q,h) \in \tilde{\mathscr{C}}_0$  
is defined to be
$$I(q,h):=\mathrm{ind}_{f}(\pi(q))+I_\mathcal{H}(h)$$
for the canonical projection onto
the Marsden-Weinstein quotient $\pi:\mu^{-1}(0) \to \bar{M}=\mu^{-1}(0)/G$.
We define the energy of a critical point by
$$E(q,h):=E_\mathcal{H}(h).$$
By abuse of notation we will also denote by $I$ and $E$ the induced
index and the induced energy on the quotient $\mathscr{C}_0$.

We next introduce flow lines with cascades which are needed to define
the boundary operator. 
\begin{fed}
For $c_1=(q_1,h_1),c_2=(q_2,h_2)  \in \tilde{\mathscr{C}}_0$ and
$m \in \mathbb{N}$ a \emph{\textbf{flow line 
from $c_1$ to $c_2$ with $m$ cascades}} 
$$v=
((w_k)_{1 \leq k \leq m},(T_k)_{1 \leq k \leq m-1})=$$
$$((u_k,\Psi_k,\Phi_k)_{1 \leq k \leq m},(T_k)_{1 \leq k \leq m-1})$$
consists of the triple of functions 
$(u_k,\Psi_k,\Phi_k) \in C^\infty_{loc}(\Theta,M\times \mathfrak{g})
\times C^\infty_\delta(\Theta,\mathfrak{g})$ 
\footnote{See (\ref{Cdelta}) for the definition of the space
$C^\infty_\delta$.}
and the nonnegative
real numbers
$T_k \in \mathbb{R}_\geq:=\{r \in \mathbb{R}:r \geq 0\}$
which have the following properties:
\begin{description}
 \item[(i)] $(u_k,\Psi_k, \Phi_k)$ are nonconstant, finite energy solutions of 
(\ref{eq}) with Hamiltonian equal to zero, i.e.  
\begin{equation}\label{eqd}
\begin{array}{c}
\partial_s u_k+X_{\Phi_k}(u_k)+J_t(u_k)(\partial_tu_k+
X_{\Psi_k}(u_k))=0\\
\partial_s \Psi_k-\partial_t\Phi_k+[\Phi_k,\Psi_k]+\mu(u_k)=0\\
u_k(s,j)\in L, \quad \Phi_k(s,j)\in \mathfrak{g}_L, 
\quad \Psi_k(s,j) \in \mathfrak{g}_L^\perp,
\end{array}
\end{equation}
where $j \in \{0,1\}$. 
\item[(ii)] There exist points $p_1 \in W^u_{\tilde{f}}(q_1)$ and
$p_2 \in W^s_{\tilde{f}}(q_2)$ such that
$\lim_{s \to -\infty}u_1(s,t)=h_1(t)p_1$ and
$\lim_{s \to \infty}u_m(s,t)=h_2(t)p_2$ uniformly in the 
$t$-variable.
\item[(iii)] For $1 \leq k \leq m-1$ there exist Morse-flow lines 
$y_k:(-\infty,\infty) \to \mu^{-1}(0) \cap L$, i.e. solutions 
of 
$$\dot{y}_k=-\nabla_{\tilde{g}} \tilde{f}(y_k),$$
and $g_k \in \mathcal{H}$ such that
$$\lim_{s \to \infty}u_k(s,t)=g_k(t)y_k(0)$$
and
$$\lim_{s \to -\infty}u_{k+1}(s,t)=g_k(t)y_k(T_k),$$
where the two limites are uniform in the $t$-variable. 
\end{description}
A \emph{\textbf{flow line with zero cascades from $c_1=(q_1,h)$ to
$c_2=(q_2,h)$}} is a tuple $(y,h)$ where y is just an ordinary 
Morse flow line from $q_1$ to $q_2$.
\end{fed}
Rccall from p. \pageref{banachgauge} that the gauge group $\mathcal{G}$
consists of smooth maps from the strip to $G$, which satisfy
appropriate boundary conditions and which decay exponentially.
For $m \in \mathbb{N}$ the group $\mathcal{G}_m$ 
consists of $m$-tuples
$$\textbf{g}=(g_k)_{1 \leq k \leq m}$$
where $g_k \in \mathcal{G}$, which have the
additional property that they form a chain, i.e.
$$\mathrm{ev}_2 (g_k)=\mathrm{ev}_1(g_{k+1}) ,\quad 1 \leq k \leq m-1.$$
For $m \geq 1$ the group $\mathcal{G}_m\times \mathbb{R}^m$ acts on the 
space of flow lines with $m$ cascades as follows
$$(u_k,\Psi_k,\Phi_k)(s) \mapsto (g_k)_*(u_k,\Psi_k,\Phi_k)(s+s_k)$$
for $1 \leq k \leq m$ and $(g_k,s_k) \in \mathcal{G} \times \mathbb{R}$.
For $m=0$ the group $\frac{\mathcal{H}}{\ker I_{\mathcal{H}} \cap
\ker E_{\mathcal{H}}} \times \mathbb{R}$ acts on 
the space of flow lines with zero cascades by
$$(y(s),h) \mapsto (y(s+s_0),g \circ h).$$
For $c_1,c_2 \in \mathscr{C}_0$
we denote the quotient of flow lines with $m$ cascades from $c_1$ to $c_2$
for $m \in \mathbb{N}_0$ by
$$\mathcal{M}_m(c_1,c_2).$$
We define the \textbf{space of flow lines with cascades from $c_1$ to $c_2$} by
$$\mathcal{M}(c_1,c_2):=
\bigcup_{m \in \mathbb{N}_0}\mathcal{M}_m(c_1,c_2).$$

Using the transversality result for the symplectic vortex equations
in \cite{frauenfelder2}, one proves in the same way as 
Theorem~\ref{manifold} that
the moduli spaces of flow lines with cascades are finite dimensional
manifolds for generic choice of the almost complex structure.  
\begin{thm}\label{mani0}
For each pair of a Morse function $f$ on $\bar{L}=(\mu^{-1}(0) \cap L)/G_L$ 
and a Riemannian metric $g$ on $\bar{L}$ which satisfy the
Morse-Smale condition, i.e. its stable and unstable manifolds intersect
transversally, there exists a subset of the space of admissible families
of almost complex structures
$$\mathcal{J}_{reg}=\mathcal{J}_{reg}(f,g) \subset \mathcal{J}$$
which is of the second category, 
i.e. $\mathcal{J}_{reg}$ is a countable intersection
of open and dense subsets of $\mathcal{J}$, and which is regular in
the following sense. 
For any two critical points
$c_1, c_2 \in \mathscr{C}_0$ the space 
$\mathcal{M}(c_1,c_2)=\mathcal{M}(c_1,c_2;J,f,g)$ 
is a smooth finite dimensional manifold. Its
dimension is given by
$$\mathrm{dim}(\mathcal{M}(c_1,c_2))=I(c_1)-I(c_2)-1.$$
If $I(c_1)-I(c_2)-1=0$, then $\mathcal{M}(c_1,c_2)$ is compact and hence
a finite set.
\end{thm}

We are now able to define moment Floer homology in the case where
the Hamiltonian vanishes.
Choose a triple $(f,g,J)$ where $f$ is a Morse function on 
$\bar{L}=(\mu^{-1}(0)\cap L)/G_L$, $g$ is a Riemannian metric on
$\bar{L}$, such that all the stable and unstable manifolds of $(f,g)$
intersect transversally, and $J \in \mathcal{J}_{reg}(f,g)$. 
As in the transversal case we define 
the chain complex $CF_k(f,g,J,L,\mu;\Lambda)$ as the $\mathbb{Z}_2$
vector space 
consisting of formal sums of the form
$$\xi=\sum_{\substack{c \in \mathscr{C}_0\\I(c)=k}}\xi_c c$$
with $\mathbb{Z}_2$-coefficients $\xi_c$ satisfying the 
finiteness condition
$$\#\{c:\xi_c \neq 0,E(c) \geq \kappa\}<\infty$$
for every $\kappa>0$. The Novikov ring $\Lambda=\Lambda_\Gamma$ 
acts naturally on $CF_*$.
Defining the boundary operator $\partial_k: CF_k
\to CF_{k-1}$ in the usual way, we obtain
homology groups
$$HF_k(f,g,J,L,\mu;\Lambda):=\frac{\mathrm{ker}\partial_{k+1}}
{\mathrm{im}\partial_k}.$$
As in theorem~\ref{continuation} one shows, that 
$HF_*(f,g,J,L,\mu;\Lambda)$ is canonically isomorphic to the moment
Floer homology groups $HF_*(L,\mu;\Lambda)$.

To compute moment Floer homology we show that contribution of the
cascades vanishes. To do that we endow the space of cascades with
the structure of a space of Arnold-Givental type.
Recall that $\dot{S}=dS(\mathrm{id})$ is the involution on the Lie algebra
induced by the antisymplectic involution $R$. 
First note that $R$ induces an involution $R_*$ on the
path space $\mathscr{P}$
by
$$R_*(x,\eta)(t):=(Rx,-\dot{S}(\eta))(1-t).$$
One easily checks that if $c \in \tilde{\mathscr{C}}_0$ then $c$ and
$R_* c$ represent the same element in $\mathscr{C}_0$.
Choose now an almost complex structure $J \in \mathcal{J}$ which is
independent of the $t$-variable and satisfies
$$R^*J=-J.$$
If $(u,\Psi,\Phi)$ is a cascade then
one verifies that
$$R_*(u,\Psi,\Phi)(s,t):=(Ru,-\dot{S}(\Psi),\dot{S}(\Phi))(s,1-t)$$
is also a cascade. Moreover, one verifies that
$$\mathrm{ev}_j(R_*(u,\Psi,\Phi))=R_*\mathrm{ev}_j(u,\Phi,\Psi), \quad
j \in \{0,1\}.$$
The following lemma shows us the relation between fixed gauge orbits
of $R_*$ and fixpoints of $R_*$.
\begin{lemma}\label{halb}
Assume that $(u,\Psi,\Phi)$ is a cascade and $g \in \mathcal{G}_{loc}$ is 
a gauge transformation such that
$$(u,\Psi,\Phi)=g_*(R_*(u,\Psi,\Phi)).$$
Then there exists $h \in \mathcal{G}_{loc}$ such that
\begin{equation}\label{bokyo}
h_*(u,\Psi,\Phi)=R_*(h_*(u,\Psi,\Phi)).
\end{equation}
\end{lemma}
\textbf{Proof:} Choose $h \in \mathcal{G}_{loc}$ such that 
$$h_* \Phi=0, \quad \lim_{s \to \infty}h_*(u,\Psi)(s,t)=(p,0)$$
where $p \in L \cap \mu^{-1}(0)$. Using the formula
$$h_*(u,\Psi,\Phi)=
(h g (Sh)^{-1})_* \circ R_* \circ h_*(u,\Psi,\Phi).$$
we get
$$0=h_*\Phi=(hg(Sh)^{-1})_*(0)=(hg(Sh)^{-1})^{-1}\partial_s (hg(Sh)^{-1})$$
and hence $hg(Sh)^{-1}$ is independent of $s$. Moreover, taking the
limit $s \to \infty$ and recalling $p=Rp$ we have
$$(hg(Sh)^{-1})(t)p=p.$$
Since $G$ acts freely on $\mu^{-1}(0)$ we obtain 
$$hg(Sh)^{-1} \equiv \mathrm{id}$$
and hence $h$ is the required gauge transformation, i.e. (\ref{bokyo})
holds. \hfill $\square$
\\ \\
\textbf{Proof of Theorem~\ref{arnold}:}
In view of the Lemma~\ref{halb} we can find in each fixed gauge orbit a
fixed  point. For fixed points we can define a sequence of involutions
whose domain is the fixed point set of the previous one in the same way
as for the pseudo-holomorphic disks in section~\ref{agi}. Using some
equivariant version of Theorem~\ref{kusta} it follows that the space of
cascades admits a Kuranishi structure of Arnold-Givental type. In
particular, the only contribution to the boundary comes from the flow
lines with zero cascades, i.e. the Morse-flow lines. This proves the theorem.
\hfill $\square$

\newpage

\appendix 

\section[Morse-Bott theory]{Morse-Bott theory}
\label{morsebott}

In this appendix we define a homology for Morse-Bott functions. Using
an idea of Piunikhin, Salamon and Schwarz, see
\cite{piunikhin-salamon-schwarz}, we define Morse-Bott
homology by counting flow lines with cascades. The homology is
independent of the choice of the Morse-Bott function and hence
isomorphic to the ordinary Morse homology.

\subsection[Morse-Bott functions]{Morse-Bott functions}

Let $(M,g)$ be a Riemannian manifold. A smooth function
$f \in C^\infty(M,\mathbb{R})$ is called \textbf{Morse-Bott} if 
$$\mathrm{crit}(f):=\{x \in M: df(x)=0\}$$
is a submanifold of $M$ and for each $x \in \mathrm{crit}(f)$ we have
$$T_x\mathrm{crit}(f)=\ker(\mathrm{Hess}(f)(x)).$$ 
\begin{ex}
Let $M=\mathbb{R}^{k_0} \times \mathbb{R}^{k_1} \times \mathbb{R}^{k_2}$. 
Write $x=(x_0,x_1,x_2)$ according to the splitting of $M$. Then
$$f(x_0,x_1,x_2)=x_1^2-x_2^2$$
is a Morse-Bott function on $M$. 
\end{ex}
\begin{ex}\label{nobot}
Let $M=\mathbb{R}$. Then 
$$f(x)=x^4$$
is no Morse-Bott function on $M$. 
\end{ex}
\begin{thm}\label{exp}
Let $(M,g)$ be a compact Riemannian manifold and $f$ a Morse-Bott function on
it. Let $y: \mathbb{R} \to M$ be a solution of 
$$\dot{y}(s)=-\nabla f(y(s)).$$
Then there exists $x \in \mathrm{crit}(f)$ and positive constants 
$\delta$ and $c$ such that
$$\lim_{s \to \infty}y(s)=x, \quad |\dot{y}(s)| \leq ce^{-\delta s}.$$
An analoguous result holds as $s$ goes to $-\infty$. 
\end{thm}
\begin{rem}
Without the Morse-Bott condition Theorem~\ref{exp} will in general not
hold. Let $M$ and $f$ be as in Example~\ref{nobot}. Then 
$$y(s):=\frac{1}{\sqrt{8s}}$$
is a solution of the gradient equation which converges to the critical
point $0$ as $s$ goes to $\infty$. But the convergence is not 
exponential. 
\end{rem}
\textbf{Proof of Theorem~\ref{exp}: }Since $M$ is compact and
$\mathrm{crit}(f)$ is normally hyperbolic there
exists $x \in \mathrm{crit}(f)$ such that
$$\lim_{s \to \infty} y(s)=x.$$
Set
$$A(s):=f(y(s))-f(x).$$
Then for some $\epsilon>0$ we have 
$$\dot{A}(s)=-|\dot{y}(s)|^2=-|\nabla f(y(s))|^2 \leq
-\epsilon A(s).$$
The last inequality follows from the Morse-Bott assumption. Hence there
exists a constant $c_0>0$ such that
$$A(s) \leq c_0 e^{-\epsilon s}.$$
This proves the theorem. \hfill $\square$

\subsection[Flow lines with cascades]{Flow lines with cascades}\label{flwc}

Let $(M,g)$ be a compact Riemannian manifold, $f$ a 
Morse-Bott function on $M$, $g_0$ a Riemannian metric on
$\mathrm{crit}(f)$, and
$h$ a Morse-function on $\mathrm{crit}(f)$. We assume that
$h$ satisfies the Morse-Smale condition, i.e. stable and unstable
manifolds intersect transversally. For a critical point
$c$ on $h$ let $\mathrm{ind}_f(c)$ be the number of negative eigenvalues
of $\mathrm{Hess}(f)(c)$ and $\mathrm{ind}_h(c)$ be the number of 
negative eigenvalues of $\mathrm{Hess}(h)(c)$. We define
$$\mathrm{Ind}(c):=\mathrm{Ind}_{f,h}(c):=
\mathrm{ind}_f(c)+\mathrm{ind}_h(c).$$ 
\begin{fed}\label{casc}
For $c_1,c_2 \in \mathrm{crit}(h),$ and $m \in \mathbb{N}$ a
\emph{\textbf{flow line from $c_1$ to $c_2$ with $m$ cascades}} 
$$(\textbf{x},\textbf{T})=((x_k)_{1 \leq k \leq m},(t_k)_{1 \leq k \leq m-1})$$
consists of $x_k \in C^\infty(\mathbb{R},M)$ and 
$t_k \in \mathbb{R}_{\geq}:=\{r \in \mathbb{R}: r \geq 0\}$ which satisfy
the following conditions:
\begin{description}
 \item[(i)] $x_k \in C^\infty(\mathbb{R},M)$ are nonconstant solutions of 
  $$\dot{x}_k=-\nabla f(x_k).$$
 \item[(ii)] There exists 
  $p \in W^u_h(c_1) \subset \mathrm{crit}(f)$ and $q \in W^s_h(c_2)$
  such that $\lim_{s \to -\infty}x_1(s)=p$ and 
  $\lim_{s \to \infty}x_m(s)=q$. 
 \item[(iii)] for $1 \leq k \leq m-1$ there are Morse-flow lines 
  $y_k \in C^\infty(\mathbb{R},\mathrm{crit}(f))$
  of $h$, i.e. solutions of 
  $$\dot{y}_k=-\nabla h(y_k),$$
  such that
  $$\lim_{s \to \infty}x_k(s)=y_k(0)$$
  and
  $$\lim_{s \to -\infty}x_{k+1}(s)=y_k(t_k).$$
\end{description}
A \emph{\textbf{flow line with zero cascades from $c_1$ to $c_2$}} 
is an ordinary Morse flow line of $h$ on $\mathrm{crit}(f)$ 
from $c_1$ to $c_2$. 
\end{fed}

\medskip
\begin{center}\input{f1.pstex_t}\end{center}
\medskip

\begin{rem}
In Definition~\ref{casc} we do not require that the Morse-flow lines
are nonconstant. It may happen that a cascade converges to a critical
point of $h$, but the flow line will only remain for a finite time on
the critical point. 
\end{rem}
We denote the space of flow lines with $m$ cascades from $c_1$ to 
$c_2 \in \mathrm{crit}(h)$ by
$$\tilde{\mathcal{M}}_m(c_1,c_2).$$
The group $\mathbb{R}$ acts by timeshift on 
the set of solutions connecting two critical points on the same
level $\tilde{\mathcal{M}}_0(c_1,c_2)$ 
and the group $\mathbb{R}^m$ acts on 
$\tilde{\mathcal{M}}_m(c_1,c_2)$ by timeshift
on each cascade, i.e.
$$x_k(s) \mapsto x_k(s+s_k).$$
We denote the quotient by
$$\mathcal{M}_m(c_1,c_2).$$
We define the \textbf{set of flow lines with cascades from $c_1$ to
$c_2$} by
$$\mathcal{M}(c_1,c_2):=\bigcup_{m \in \mathbb{N}_0}\mathcal{M}_m(c_1,c_2).$$
Immediately from the gradient equation the following lemma follows.
\begin{lemma}
If $f(c_1)<f(c_2)$ 
then $\mathcal{M}(c_1,c_2)$ is empty. If $f(c_1)=f(c_2)$ then 
$\mathcal{M}(c_1,c_2)$ contains only flow lines with zero cascades. If
$f(c_1)>f(c_2)$ then $\mathcal{M}(c_1,c_2)$ contains no flow line with
zero cascades. 
\end{lemma}

A sequence of flow lines with cascades may break up in the limit into
a connected chain of flow lines with cascades. To deal with this
phenomenon, we make the following definitions.
\begin{fed}
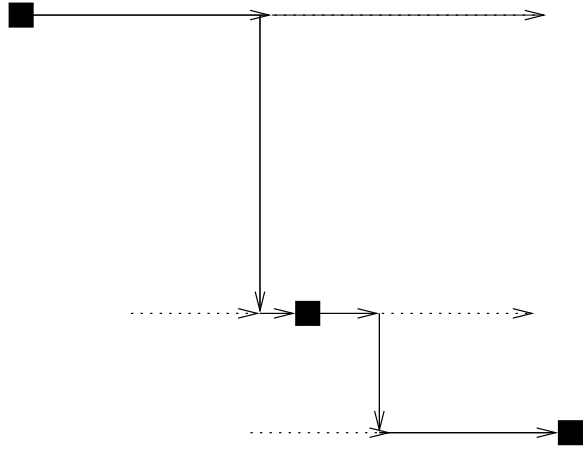
Let $c,d \in \mathrm{crit}(h)$. A \emph{\textbf{broken flow line with
cascades}} from $c$ to $d$
$$\textbf{v}=(v_j)_{1 \leq j \leq \ell}$$
for $\ell \in \mathbb{N}$ consists of flow lines with cascades
$v_j$ from $c_{j-1}$ to $c_j \in \mathrm{crit}(h)$
for $0 \leq j \leq \ell$ such that $c_0=c$ and $c_\ell=d$. 
\end{fed}

\medskip
\begin{center}\input{f2.pstex_t}\end{center}
\medskip

\begin{fed}\label{fgc}
Assume that $c,d \in \mathrm{crit}(h)$. Suppose that 
$v^\nu$ for $\nu \in \mathbb{N}$ is a sequence of flow lines
with cascades which satisfies the following
condition. There exists
$\nu_0 \in \mathbb{N}$ such that
for every $\nu\geq \nu_0$ it holds that
$v^\nu$ is a flow line with cascades from $c$ to $d$. There
are two cases. In the first case $c$ and $d$ lie on the same level and
hence
$v^\nu \in C^\infty(\mathbb{R},\mathrm{crit}(f))$
is a flow line with zero cascades for every $\nu \geq \nu_0$, in
the second case $c$ and $d$ lie on different levels and hence
$v^\nu=((x_k^\nu)_{1 \leq k \leq m^\nu},(t_k)_{1 \leq k \leq m^\nu-1})$
is a flow line with at least one cascade for every $\nu \geq \nu_0$.
We say that $v^\nu$ 
\emph{\textbf{Floer-Gromov converges}} to a broken flow line with
cascades $\textbf{v}=(v_j)_{1 \leq j \leq \ell}$ from $c$ to $d$ if 
the following holds. 
\begin{description}
 \item[(a)] In the first case, where 
  the $v^\nu$'s are flow lines with zero
  cascades for large enough $\nu$'s,
  all $v_j$'s are flow lines with
  zero cascades and there exists real numbers $s_j^\nu$
  for $\nu \geq \nu_0$ such that
  $(s^\nu_j)_*(v^\nu)(\cdot):=v^\nu(\cdot+s^\nu_j)$ converges
  in the $C^\infty_{loc}$-topology to $v_j$.
 \item[(b)] In the second case, where the $v^\nu$'s have at least
  one cascade for $\nu$ large enough, we require the following conditions.
  \begin{description}
   \item[(i)] If $v_j \in C^\infty(\mathbb{R},\mathrm{crit}(f))$
    is a flow line with zero cascades, then there exists 
    a sequence of solutions 
    $y^\nu_j  \in C^\infty(\mathbb{R},\mathrm{crit}(f))$ of 
    $\dot{y}^\nu_j=-\nabla h(y^\nu_j)$
    converging in $C^\infty_{loc}$ to $v_j$, a sequence of real numbers
    $s^\nu_j$, and a sequence of integers $k^\nu \in [1,m^\nu]$
    such that either  $\lim_{s \to -\infty}x^\nu_{k^\nu}(s)=y^\nu_j(s^\nu_j)$ or
    $\lim_{s \to \infty}x^\nu_{k^\nu}(s)=y^\nu_j(s^\nu_j)$.
  \item[(ii)] If $v_j$ is a flow line with at least one cascade,
   then we write 
   $v_j=((x_{i,j})_{1 \leq i \leq m_j},(t_{i,j})_{1 \leq i \leq m_j-1})
   \in \tilde{\mathcal{M}}_{m_j}(c_{j-1},c_j)$ for $m_j \geq 1$.
   We require that there exist surjective maps 
   $\gamma^\nu \colon [1, \sum_{p=1}^\ell m_p] \to [1,m^\nu]$,
   which are monotone increasing, 
   i.e. $\gamma^\nu(\lambda_1) \leq \gamma^\nu(\lambda_2)$ for
   $\lambda_1 \leq \lambda_2$, and real numbers $s^\nu_\lambda$
   for every $\lambda \in [1,\sum_{p=1}^\ell m_p]$, such that
   $$(s^\nu_\lambda)_*x^\nu_{\gamma^\nu(\lambda)}(\cdot)
   =x^\nu_{\gamma^\nu(\lambda)}(\cdot+s^\nu_\lambda)
   \stackrel{C^\infty_{loc}}{\longrightarrow} x_\lambda$$
   where $x_\lambda=x_{i,j}$ such that $\lambda=\sum_{p=1}^jm_p+i$.
   For $\lambda \in [1,\sum_{p=1}^\ell m_p-1]$ we set
   $$\tau_\lambda=\left\{\begin{array}{cc}
   t_{i,j} & \lambda=\sum_{p=1}^jm_p+i, \quad 0 <i<m_{j+1}\\
   \infty & \lambda=\sum_{p=1}^jm_p\end{array}\right.$$
   and
   $$\tau^\nu_\lambda=\left\{\begin{array}{cc}
   t^\nu_{\gamma^\nu(\lambda)} & \lambda=\mathrm{max}\{
   \lambda' \in [1,\sum_{p=1}^\ell m_p-1]: \gamma^\nu(\lambda')=
   \gamma^\nu(\lambda)\}\\
   0 & \mathrm{otherwise.} \end{array}\right.$$
   Now we require, that
   $$\lim_{\nu \to \infty} \tau^\nu_\lambda=\tau_\lambda.$$
   Here we use the convention that a sequence of real numbers
  $\tau^\nu$ converges to infinity, if for every $n \in \mathbb{N}$
  there exists a $\nu_0(n) \in \mathbb{N}$ such that $\tau^\nu \geq n$ for
  $\nu \geq \nu_0(n)$.
  \end{description}   
\end{description}
\end{fed}

\begin{thm}[Compactness]\label{fcc}
Let $v^\nu$ be a sequence of flow lines with cascades. Then there exists
a subsequence $v^{\nu_j}$ and a broken flow line with cascades
$\textbf{v}$ such that $v^{\nu_j}$ Floer-Gromov converges to 
$\textbf{v}$. 
\end{thm}
\textbf{Proof:} First assume that there exists a subsequence $\nu_j$
such that $v^{\nu_j}$ are flow lines with zero cascades. In this case 
Floer-Gromov convergence to a broken flow line without cascades
follows from the classical case, see 
\cite[Proposition 2.35]{schwarz1}. Otherwise pick
a subsequence $\nu_j$ of $\nu$ such that $v^{\nu_j}$ are flow
lines with at least one cascade. Since the number
of critical points of $h$ is finite, we can 
perhaps after passing over to a further
subsequence assume that all the $v^{\nu_j}$ are flow lines with 
cascades from a common critical point $c$ of $h$ to a common
critical point $d$ of $h$. Since the number of connected components of
$\mathrm{crit}(f)$ is finite, we can assume by passing to a further
subsequence that the number of cascades $m^\nu=m$ is independent on $\nu$.
\\ 
For simplicity of notation we denote the subsequence $\nu_j$ again by
$\nu$. We consider the sequence of points
$p^\nu:=\lim_{s \to \infty}x^\nu_1(s) \in \mathrm{crit}(f).$ 
Let $y^\nu \in C^\infty(\mathbb{R},\mathrm{crit}(f))$ be the unique solution
of the problem
$$\dot{y}^\nu(s)=-\nabla h(y^\nu(s)), \quad y^\nu(0)=p^\nu.$$ 
Note
that for every $\nu$ 
$$\lim_{s \to -\infty}y^\nu(s)=c.$$
Using again \cite[Proposition 2.35]{schwarz1} it follows
that perhaps passing over to a subsequence (also denoted by $\nu$) there
exists $p \in \mathrm{crit}(f)$, a nonnegative integer $m_0 \in \mathbb{N}$,
and a sequence of Morse-flow lines 
$v_j \in C^\infty(\mathbb{R},\mathrm{crit}(f))$ with respect
to the gradient flow of $h$ for $1 \leq j \leq m_0$
such that the following holds. 
\begin{description}
 \item[(i)] $\lim_{\nu \to \infty} p^\nu=p$,
 \item[(ii)] If $m_0=0$, then $p \in W^u_h(c)$,
 \item[(iii)] If $m_0 \geq 1$, then 
\begin{eqnarray*}
\lim_{s \to -\infty}v_1(s)&=&c, \\
\lim_{s \to \infty}v_j(s)&=&\lim_{s \to -\infty}v_{j+1}(s), \quad 1 \leq j \leq
m_0-1, \\
p &\in& W^u_h(\lim_{s \to \infty}v_{m_0}(s)).
\end{eqnarray*}
\end{description}
Using induction on $\mu \in [1,m]$, where $m$ is the number of cascades of each
$v^\nu$, the following claim follows as \cite[Proposition 2.35]{schwarz1}. 
\\ \\
\textbf{Claim: }\emph{There exist a subsequence of $\nu$ (still
denoted by $\nu$), a nonnegative integer
$\ell_\mu^1 \in \mathbb{N}_0$, a positive integer
$\ell_\mu^2 \in \mathbb{N}$,
a broken flow line with cascades 
$\textbf{v}_{\mu}=(v_j)_{1 \leq j \leq \ell^1_\mu}$
from $c$ to some critical point
$c_\mu$ of $h$, a sequence of cascades 
$x^\mu_k \in C^\infty(\mathbb{R},M)$ for $1 \leq k \leq \ell^2_\mu$, i.e.
of nonconstant solutions of the ordinary differential equation
$\dot{x}_k^\mu=-\nabla f(x_k^\mu)$, and a sequence of nonnegative
real number $t^\mu_k$ for $1 \leq k \leq \ell^2_\mu-1$
such that the following conditions are satisfied.
\begin{description}
  \item[(i)] If $\ell^1_\mu=0$, then $\lim_{s \to -\infty}x_1(s)
   \in W^u_h(c)$, otherwise $\lim_{s \to -\infty}x_1(s) \in
   W^u_h(c_\mu)$,
  \item[(ii)] For $1 \leq k \leq \ell^2_\mu-1$ there are Morse-flow
   lines $y^\mu_k \in C^\infty(\mathbb{R},\mathrm{crit}(f))$ of $h$, such that
   $$\lim_{s \to \infty} x^\mu_k(s)=y^\mu_k(0), \quad
   \lim_{s \to -\infty}x^\mu_{k+1}(s)=y^\mu_k(t^\mu_k),$$
  \item[(iii)] If $v_j \in C^\infty(\mathbb{R},\mathrm{crit}(f))$
    is a flow line with zero cascades, then there exists 
    a sequence of solutions 
    $y^\nu_j  \in C^\infty(\mathbb{R},\mathrm{crit}(f))$ of 
    $\dot{y}^\nu_j=-\nabla h(y^\nu_j)$
    converging in $C^\infty_{loc}$ to $v_j$, a sequence of real numbers
    $s^\nu_j$, and an integer $k \in [1,\mu]$
    such that $\lim_{s \to -\infty}x^\nu_k(s)=y^\nu_j(s^\nu_j)$.
  \item[(iv)] If $v_j$ is a flow line with at least one cascade,
   then we are able to write 
   $v_j=((x_{i,j})_{1 \leq i \leq m_j},(t_{i,j})_{1 \leq i \leq m_j-1})
   \in \tilde{\mathcal{M}}_{m_j}(c_{j-1},c_j)$ for $m_j \geq 1$.
   We require that there exists a surjective map
   $\gamma_\mu \colon [1, \sum_{p=1}^{\ell^1_\mu} m_p+
   \ell^2_\mu] \to [1,\mu]$,
   which is monotone increasing, 
   i.e. $\gamma_\mu(\lambda_1) \leq \gamma_\mu(\lambda_2)$ for
   $\lambda_1 \leq \lambda_2$, and real numbers $s^\nu_\lambda$
   for every 
   $\lambda \in [1,\sum_{p=1}^{\ell^1_\mu} m_p+\ell^2_\mu]$, such that
   $$(s^\nu_\lambda)_*x^\nu_{\gamma^\nu(\lambda)}(\cdot)
   =x^\nu_{\gamma^\nu(\lambda)}(\cdot+s^\nu_\lambda)
   \stackrel{C^\infty_{loc}}{\longrightarrow} x_\lambda$$
   where $x_\lambda=x_{i,j}$ if $\lambda=\sum_{p=1}^jm_p+i$
   for $0 \leq j \leq \ell^1_\mu$ and $0 \leq i \leq m_j-1$, or
   $x_\lambda=x^\mu_{\lambda-\sum_{p=1}^{\ell^1_\mu}m_p}$
   for $\lambda \in [\sum_{p=1}^{\ell^1_\mu}m_p+1,
   \sum_{p=1}^{\ell^1_\mu}m_p+\ell^2_\mu]$.
   For $\lambda \in [1,\sum_{p=1}^{\ell^1_\mu}m_p+\ell^2_\mu-1]$ we set
   $$\tau_\lambda=\left\{\begin{array}{cc}
   t_{i,j} & \lambda=\sum_{p=1}^jm_p+i, \quad 0 <i<m_{j+1}\\
   \infty & \lambda=\sum_{p=1}^jm_p\\
   t^\mu_{\lambda-\sum_{p=1}^{\ell^1_\mu}m_p} &
   \lambda \in [\sum_{p=1}^{\ell^1_\mu}m_p+1,\sum_{p=1}^{\ell^2_\mu}m_p
   +\ell^2_\mu]\end{array}\right.$$
   and
   $$\tau^\nu_\lambda=\left\{\begin{array}{cc}
   t^\nu_{\gamma^\nu(\lambda)} & \lambda=\mathrm{max}\{
   \lambda' \in [1,\sum_{p=1}^{\ell^1_\mu} m_p+\ell^2_\mu-1]: 
   \gamma^\nu(\lambda')=
   \gamma^\nu(\lambda)\}\\
   0 & \mathrm{otherwise.} \end{array}\right.$$
   We require, that
   $$\lim_{\nu \to \infty} \tau^\nu_\lambda=\tau_\lambda.$$
  \item[(v)] $\lim_{\nu \to \infty}\lim_{s \to \infty}x^\nu_\mu(s)$ exists
   and equals $\lim_{s \to \infty}x^\mu_{\ell^2_\mu}(s)$. 
\end{description}}
Given the claim for $\mu=m$, the Theorem follows now by applying
\cite[Proposition 2.35]{schwarz1} again. \hfill $\square$

\begin{thm}\label{manifold}
Let $c_1,c_2 \in \mathrm{crit}(h)$. For generic choice of the Riemannian metric
$g$ of $M$ the space $\mathcal{M}(c_1,c_2)$ is a smooth finite dimensional
manifold. Its dimension is given by
$$\dim \mathcal{M}(c_1,c_2)=\mathrm{ind}(c_1)-\mathrm{ind}(c_2)-1.$$
If $\mathrm{Ind}(c_1)-\mathrm{Ind}(c_2)=1$ then 
$\mathcal{M}(c_1,c_2)$ is compact and hence a finite set. 
\end{thm}

The rest of this subsection is devoted to the
proof of Theorem~\ref{manifold}.
Choose $0<\delta<\min(\{|\lambda|: \lambda \in \sigma(\mathrm{Hess}(f)(x))
\setminus \{0\}, 
x \in \mathrm{crit}(f)\}$. For a smooth cutoff function $\beta$ such that
$\beta(s)=-1$ if $s<0$ and $\beta(s)=1$ if $s>1$ define
$$\gamma_\delta:\mathbb{R} \to \mathbb{R}, \quad s \mapsto 
e^{\delta \beta (s)s}.$$
Let $\Omega$ be an open subset of $\mathbb{R}$. We define the
$|| \,||_{k,p,\delta}$-norm for a locally integrable function
$f: \Omega \to \mathbb{R}$ with weak derivatives up to order $k$ by
$$||f||_{k,p,\delta}:=\sum_{i=0}^k||\gamma_\delta \partial^i f||_p.$$
We denote
$$W^{k,p}_\delta(\Omega):=\{f \in W^{k,p}(\Omega):
||f||_{k,p,\delta}<\infty\}=\{f \in W^{k,p}(\Omega):\gamma_\delta f 
\in W^{k,p}(\Omega)\}.$$
We also set
$$L^p_\delta(\Omega):=W^{0,p}_\delta(\Omega).$$
Let 
$$T_h(t) \in \mathrm{Diff}(\mathrm{crit}(f))$$
be the smooth family of diffeomorphisms which assigns to 
$p \in \mathrm{crit}(f)$ the point
$x_p(t)$ where $x_p$ is the unique flow line of $h$ with 
$x_p(0)=p$. We define
$$\mathcal{B}:=\mathcal{B}_\delta^{1,p}(M,f,h)$$
as the Banach manifold consisting of all  tuples
$v=((x_j)_{1 \leq j \leq m},(t_j)_{1 \leq j \leq m-1})
\in (W^{1,p}_{loc}(\mathbb{R},M))^{m}\times \mathbb{R}_+^{m-1}$ where
$\mathbb{R}_+:=\{r \in \mathbb{R}: r>0\}$ and $m \in \mathbb{N}$
which satisfy the following conditions:
\begin{description}
 \item[(Asympotic behaviour)] For $1 \leq j \leq m$ there exists $p_j, q_j \in 
  \mathrm{crit}(f)$, $\xi_{1,j} \in W^{1,p}_\delta((-\infty,T], T_{p_j}M)$,
  $\xi_{2,j} \in W^{1,p}_\delta([T,\infty),T_{q_j}M)$ for $T \in \mathbb{R}$
  such that
  $$x_j(s)=\exp_{p_j}(\xi_{1,j}(s)), \quad s \leq -T, \qquad
  x_j(s)=\exp_{q_j}(\xi_{2,j}(s)), \quad s \geq T,$$
  where $\exp$ is taken with respect to the Riemannian metric $g$ of $M$.
 \item[(Connectedness)] $p_{j+1}=T_h(t_j) q_j$ for $1 \leq j \leq m-1$.
\end{description}
To define local charts on $\mathcal{B}$ choose
$v=((x_j)_{1 \leq j \leq m},(t_j)_{1 \leq j \leq m-1}) \in \mathcal{B}$
such that all the $x_j$ for $1 \leq j \leq m$ are smooth and define
a neighbourhood of $v$ in $\mathcal{B}$ via the exponential map of $g$
\footnote{Note that the differentiable structure of $\mathcal{B}$
is independent of the metric $g$ on $M$}.
There are two natural smooth evaluation maps 
$$\mathrm{ev}_1, \mathrm{ev}_2: \mathcal{B} \to \mathrm{crit}(f), \qquad
\mathrm{ev}_1(v)=p_1, \quad \mathrm{ev}_2(v)=q_m.$$
After choosing cutoff functions and smooth trivializations 
$\chi_{p_j}$ and $\chi_{q_j}$ of $TM$
near $p_j$ respectively $q_j$ the tangent space of $\mathcal{B}$ at $v$
can naturally be identified with tuples
$$\zeta=((\xi_{j,0},\xi_{j,1},\xi_{j,2})_{1 \leq j \leq m},(\tau_j)_{1 \leq j
\leq m-1}) \in$$
$$\bigoplus_{j=1}^m (W^{1,p}_\delta(\mathbb{R},x^*_j TM)
\times T_{p_j}\mathrm{crit}(f) \times T_{q_j}\mathrm{crit}(f)) 
\times \mathbb{R}^{m-1}$$
which satisfy
\begin{equation}
d T_h(t_j) \xi_{j,2}+\frac{d}{dt}(T_h(0)q_j)\tau_j=\xi_{j+1,1} 
\quad  1 \leq j \leq m-1.
\end{equation}
$T_v \mathcal{B}$ is a Banach space with norm
\begin{equation}\label{norm}
||\zeta||:=\sum_{j=1}^m(||\xi_{j,0}||_{1,p,\delta}+||\xi_{j,1}||
+||\xi_{j,2}||)+\sum_{j=1}^{m-1}|\tau_j|
\end{equation}
Let $\mathcal{E}$ be the Banachbundle over $\mathcal{B}$ whose fiber
at $v \in \mathcal{B}$ is given by
$$\mathcal{E}_v:=\bigoplus_{j=1}^m L^p_\delta(\mathbb{R},x^*_j TM).$$
Set
$$\tilde{\mathcal{M}}:=\mathcal{F}^{-1}(0)$$ where
$$\mathcal{F}:\mathcal{B} \to \mathcal{E}, \quad
v \mapsto (\dot{x}_k+\nabla f(x_k))_{1 \leq k \leq m}.$$
Where $\nabla=\nabla_g$ is the Levi-Civita connection of the Riemannian
metric $g$ on $M$.
Note that $\mathcal{F}=\mathcal{F}_g$ depends on $g$. Let
$$D_v:=D \mathcal{F}(v): T_v \mathcal{B} \to \mathcal{E}_v$$ 
be the vertical differential of $\mathcal{F}$ at $v \in \mathcal{F}^{-1}(0)$. 
If $p \in \mathrm{crit}(f)$ denote
by $\dim_p \mathrm{crit}(f)$ the local dimension of 
$\mathrm{crit}(f)$ at $p$. Note that it follows from the Morse-Bott condition
that $\dim_p \mathrm{crit}(f)$ equals the dimension of the kernel of
$\mathrm{Hess}(f)(p)$. Then we have
\begin{lemma}
$D_v$ is a Fredholm-operator of index
$$\mathrm{ind}(D_v)=
\mathrm{ind}_f(\mathrm{ev}_1(v))+\dim_{\mathrm{ev}_1(v)} \mathrm{crit}(f)
-\mathrm{ind}_f(\mathrm{ev}_2(v))+m-1,$$
where $m=m(v)$ is the number of cascades.
\end{lemma}
\textbf{Proof: }For $1 \leq j \leq m$ let
$$D_{v,j}: W^{1,p}_\delta(\mathbb{R},x^*_j TM) \to L^p_\delta(\mathbb{R},
x^*_j TM)$$
be the restriction of $D_v$ to $W^{1,p}_\delta(\mathbb{R},x^*_j TM)$. 
It suffices to show, that $D_{v,j}$ is a Fredholm-operator of index
$$\mathrm{ind}(D_{v,j})=\mathrm{ind}_f(p_j)-\mathrm{ind}_f(q_j)
-\dim_{q_j}\mathrm{crit}(f).$$
Write
$$D_{v,j}:=\partial_s +A(s).$$
Then 
$$A_1:=\lim_{s \to -\infty}A(s)=\mathrm{Hess}(f)(p_j), \quad
A_2:=\lim_{s \to \infty}A(s)=\mathrm{Hess}(f)(q_j).$$
Define the continuous isomorphisms
$$\phi_1: L^p_\delta \to L^p, \quad f \mapsto f \gamma_\delta, \qquad
\phi_2: W^{1,p}_\delta \to W^{1,p}, \quad f \mapsto f \gamma_\delta.$$
Define
$$\tilde{D}_v:W^{1,p}(\mathbb{R}, x^*_j TM) \to 
L^p(\mathbb{R},x^*_j TM)$$
by
$$\tilde{D}_v:=\phi_1 D_v \phi_2^{-1}.$$
$D_v$ is exactly a Fredholm-operator if $\tilde{D}_v$ is a Fredholm-operator.
In this case
$$\mathrm{ind}(D_v)=\mathrm{ind}(\tilde{D}_v).$$
For $\xi \in W^{1,p}(\mathbb{R},x_j^* TM)$ we calculate
\begin{eqnarray*}
\tilde{D}_v \xi&=&\phi_1 D_v \phi_2^{-1}\xi\\
&=&\phi_1 D_v(\xi \gamma_{-\delta})\\
&=&\phi_1(\partial_s(\xi \gamma_{-\delta})+A(s)\xi \gamma_{-\delta})\\
&=&\phi_1((\partial_s \xi)\gamma_{-\delta}+(\delta \beta(s)+
\delta \partial_s \beta(s) s)\xi\gamma_{-\delta}+A(s)\xi \gamma_{-\delta})\\
&=&\partial_s \xi+(A(s)+\delta(\beta(s)+\partial_s \beta(s) s)\mathrm{id})\xi.
\end{eqnarray*}
Hence $\tilde{D}_v$ is given by
$$\tilde{D}_v=\partial_s +B(s)$$
where
$$B(s)=A(s)+\delta(\beta(s)+\partial_s \beta(s)s)\mathrm{id}.$$
Let
$$B_j:=A_j+(-1)^j \delta, \quad j \in \{1,2\}.$$
Then 
$$\lim_{s \to -\infty} B(s)=B_1, \quad \lim_{s \to \infty} B(s)=B_2$$
and the $B_j$ are invertible by our choice of
$\delta$. If $p=2$ then it follows from
\cite{robbin-salamon1} that $\tilde{D}_v$ is a Fredholm-operator of the
required index. For general $p$ the lemma follows from 
\cite{salamon1}. See also \cite{schwarz2}, for an
alternative proof. \hfill $\square$
\\ \\
For $n \in \mathbb{N}$ we define the evaluation maps
$$\mathrm{EV}_n: \tilde{\mathcal{M}}^n \to 
\mathrm{crit}(f)^n \times \mathrm{crit}(f)^n \cong 
\mathrm{crit}(f)^{2n},$$
$$(v_1,\ldots,v_n) \mapsto (\mathrm{ev}_1(v_j)_{1 \leq j \leq n},
\mathrm{ev}_2(v_j)_{1 \leq j \leq n}),$$
where $\tilde{\mathcal{M}}=\mathcal{F}^{-1}(0)$ as introduced above. 
For critical points $c_1,c_2$ of $h$ we define
$A_n(c_1,c_2)$ to be the submanifold of 
$\mathrm{crit}(f)^n \times \mathrm{crit}(f)^n$ consisting of all tuples
$((p_j)_{1 \leq j \leq n},(q_j)_{1 \leq j \leq n})
\in \mathrm{crit}(f)^n \times \mathrm{crit}(f)^n$ such that
$p_1 \in W^u_h(c_1)$, $q_n \in W^s_h(c_2)$, and $p_{j+1}=q_j$ for
$1 \leq j \leq n-1$.
We shall prove the following theorem.
\begin{thm}\label{mani}
For generic Riemannian metric $g$ on $M$ the set $\tilde{\mathcal{M}}$ has
the structure of a finite dimensional manifold. Its local dimension
\begin{equation}\label{dim}
\dim_v \tilde{\mathcal{M}} = \mathrm{ind} D_v
=\mathrm{ind}_f(\mathrm{ev}_1(v))+\dim_{\mathrm{ev}_1(v)} \mathrm{crit}(f)
-\mathrm{ind}_f(\mathrm{ev}_2(v))+m-1,
\end{equation}
where $m=m(v)$ is the number of cascades. 
Moreover, for every
$n \in \mathbb{N}$ and for every $c_1,c_2 \in \mathrm{crit}(h)$ the
evaluation maps $\mathrm{EV}_n$ intersects $A_n(c_1,c_2)$ 
transversally.  
\end{thm}
\begin{fed}\label{fhg0reg}
Assume that $f$ is Morse-Bott function on a compact manifold $M$, 
$h$ is a Morse-function on $\mathrm{crit}(f)$, and 
$g_0$ is a Riemannian metric
on $\mathrm{crit}(f)$, such that $h$ and $g_0$ satisfy
the Morse-Smale condition, i.e. stable and unstable manifolds
of the gradient flow of $h$ with respect to $g_0$ on $\mathrm{crit}(f)$
intersect transversally. We say that a Riemannian metric $g$ on $M$
is \emph{\bf{$(f,h,g_0)$-regular}} if it satisfies the conditions
of Theorem~\ref{mani}.
\end{fed}
\textbf{Proof of Theorem~\ref{mani}: }For a positive integer $\ell$ let 
$\mathcal{R}^\ell$ be the Banach manifold of Riemannian metrics on $M$
of class $C^\ell$. Let
$$\mathcal{F}: \mathcal{R}^\ell \times \mathcal{B} \to \mathcal{E}$$
be defined by
$$\mathcal{F}(g,v)=(\dot{x}_k+\nabla_g f(x_k))_{1 \leq k \leq m}
=(\dot{x}_k+g^{-1}df(x_k))_{1 \leq k \leq m},$$
where $g^{-1} \colon T^*M \to TM$ is defined as the inverse of the map
defined by
$$g(\nabla_g f, \cdot)=df(\cdot).$$
We will prove that the universal moduli space
$$\mathcal{U}^\ell:=\{(v,g) \in \mathcal{B} \times \mathcal{R}^\ell:
\mathcal{F}(v,g)=0\}$$
is a separable manifold of class $C^\ell$. To show that, we have to verify
that
$$D_{v,g}:T_v \mathcal{B} \times T_g \mathcal{R}^\ell \to 
\mathcal{E}_v$$
given by
\begin{eqnarray*}
D_{v,g}(\zeta,A) &=&
(\partial_s \zeta_k+\nabla_{\zeta_k}\nabla f(x_k))_{1\leq k \leq m}
-(g^{-1}Ag^{-1}df(x_k))_{1 \leq k \leq m}\\
&=& D_v\zeta-(g^{-1}Ag^{-1}df(x_k))_{1 \leq k \leq m}
\end{eqnarray*}
is onto for every $(v,g) \in \mathcal{U}^\ell$. Here 
$\nabla_{\zeta_k}$ is the Levi-Civita connection of the metric $g$. 

The tangent space of $\mathcal{R}^\ell$ at $g \in \mathcal{R}^\ell$
consists of all symmetric $C^\ell$-sections from $M$ to
$T^*M \times T^* M$. Since $D_v$ is a Fredholm operator, it has a closed
range and a finite dimensional cokernel. Hence $D_{v,g}$ has a closed
range and a finite dimensional cokernel and it only remains to prove
that its range is dense. To see this, let 
$$\eta \in (\mathcal{E}_v)^*=\bigoplus_{j=1}^m L^q_{-\delta}
(\mathbb{R}, x^*_j TM), \quad \frac{1}{p}+\frac{1}{q}=1,$$
such that $\eta$ vanishes on the range of $D_{v,g}$, i.e.
\begin{equation}\label{reg1+}
\sum_{j=1}^m \int_{\mathbb{R}}\langle \eta_j, (D_v \zeta)_j\rangle ds =0
\end{equation}
for every $\zeta \in T_v \mathcal{B}$ and 
\begin{equation}\label{reg2+}
\sum_{j=1}^m \int_{\mathbb{R}}\langle \eta_j,g^{-1}A g^{-1}df(x_j)\rangle ds
=0
\end{equation}
for every $A \in T_g \mathcal{R}^\ell$. It follows from (\ref{reg1+}) that
$\eta_j$ is continuously differentiable
for $1 \leq j \leq m$. Now (\ref{reg2+}) implies
that $\eta$ vanishes identically. This proves that $D_{v,g}$ is onto for
every $(v,g) \in \mathcal{U}^\ell$. Now it follows from the implicit function
theorem that $\mathcal{U}^\ell$ is a Banach-manifold. 

The differential $d \pi^\ell$ of the projection
$$\pi^\ell: \mathcal{U}^\ell \to \mathcal{R}^\ell, \quad
(v,g) \mapsto g$$
at a point $(v,g) \in \mathcal{U}^\ell$ is just the projection
$$d \pi^\ell(v,g)\colon T_{(v,g)}\mathcal{U}^\ell \to T_g \mathcal{R}^\ell,
\quad (\zeta,A) \mapsto A.$$
The kernel of $d\pi^\ell(v,g)$ is isomorphic to the
kernel of $D_v$. Its image consists of all $A$ such that
$(g^{-1} A g^{-1} df(x_k))_{1 \leq k \leq m} \in \mathrm{im}D_v$.
Moreover, $\mathrm{im} d\pi^\ell(v,g)$ is a closed subspace of
$T_g \mathcal{R}^\ell$, and, since $D_{v,g}$ is onto, it has the
same (finite) codimension as the image of $D_v$. Hence
$d\pi^\ell(v,g)$ is a Fredholm operator of the same index as $D_{v,g}$.
In particular, the projection $\pi^\ell$ is a Fredholm map and
it follows from the Sard-Smale theorem that for
$\ell$ sufficiently large, the set $\mathcal{R}^\ell_{reg}$ of
regular values of $\pi^\ell$ is dense in $\mathcal{R}^\ell$. Note
that $g \in \mathcal{R}^\ell$ is a regular value of $\pi^\ell$
exactly if $D_v$ is surjective for every $v \in \mathcal{F}_g^{-1}(0)$.
Here $\mathcal{F}_g=\mathcal{F}(\cdot,g)$.

For $c>0$ let $\mathcal{U}^{c,\ell} \subset \mathcal{U}^\ell$ be the
set of pairs $(v,g) \in \mathcal{U}^\ell$ such that
$$||\partial_s x_j(s)||\leq ce^{-|s|/c}, \quad 1 \leq j \leq m,
\qquad \frac{1}{c} \leq t_k \leq c, \quad 1 \leq j \leq m-1.$$
The space 
$$\tilde{\mathcal{M}}^{\ell,c}(g):=\{v:(v,g) \in \mathcal{U}^{\ell,c}\}$$
is compact for every $g$. Indeed, the uniform exponential decay 
prevents the cascades from breaking up into several pieces. 
It follows that the set 
$\mathcal{R}^{\ell,c}_{reg}$ consisting of all $g \in \mathcal{R}^\ell$
such that $D_v$ is onto for all $(v,g) \in \tilde{\mathcal{M}}^{\ell,c}(g)$
is open and dense in $\mathcal{R}^\ell$. Hence the set
$$\mathcal{R}^{\infty,c}_{reg}:=\mathcal{R}^{\ell,c}_{reg} \cap \mathcal{R}$$
is dense in $\mathcal{R}^\ell$ with respect to the $C^\ell$-topology. 
Here $\mathcal{R}=\mathcal{R}^\infty$ denotes the Fr\'echet manifold of
smooth metrics on $M$. Since this holds for
every $\ell$ if follows that
$\mathcal{R}^{\infty,c}$ is dense in $\mathcal{R}$ with respect to the
$C^\infty$-topology. Using compactness again one obtains that
$\mathcal{R}^c_{reg}$ is also $C^\infty$-open. It follows that the set
$$\mathcal{R}_{reg}:=\bigcap_{c \in \mathbb{N}} \mathcal{R}^{\infty,c}_{reg}$$
is a countable intersection of open and dense subsets of
$\mathcal{R}$. 

To prove that the evaluation maps $\mathrm{EV}_n$ intersect
$A_n(c_1,c_2)$ transversally for generic $g$, we show that the evaluation
maps $\mathrm{ev}_j: \mathcal{U}^\ell \to \mathrm{crit}(f)$ for 
$j \in \{1,2\}$
are submersive. Let $(v,g) \in \mathcal{U}^\ell$ and
$\xi \in T_{\mathrm{ev}_j(v,g)}\mathrm{crit}(f)$. We have to
show that there exists $(\zeta,A) \in \mathrm{ker} D_{v,g}$ 
such that
\begin{equation}\label{submers}
d(\mathrm{ev}_j)(v,g)(\zeta,A)=\xi.
\end{equation}
Choose some arbitrary $(\zeta_0,A_0) \in T_v \mathcal{B} \times
T_g \mathcal{R}^\ell$ such that
$$d(\mathrm{ev}_j)(v,g)(\zeta_0,A_0)=\xi.$$ 
In the same way as one proved that $D_{v,g}$ is surjective
one can also show that already $D_{v,g}$ restricted
to $\{\zeta \in T_v \mathcal{B}: d (\mathrm{ev}_j) \zeta=0\}
\times T_g \mathcal{R}^\ell$ is surjective. Hence there
exist $(\zeta_1,A_1) \in T_v \mathcal{B} \times T_g \mathcal{R}^\ell$
such that
$$D_{v,g}(\zeta_0,A_0)=D_{v,g}(\zeta_1,A_1)$$
and
$$d(\mathrm{ev}_j)(v,g)(\zeta_1,A_1)=0.$$
Now set 
$$(\zeta,A):=(\zeta_0-\zeta_1,A_0-A_1).$$
Then $(\zeta,A)$ lies in the kernel of $D_{v,g}$ and satisfies 
(\ref{submers}). This proves the theorem. \hfill $\square$
\\ \\
For $c_1,c_2 \in \mathrm{crit}(h) \subset \mathrm{crit}(f)$ define
\begin{eqnarray*}
\tilde{\mathcal{M}}^-(c_1)&:=&\{v \in \tilde{\mathcal{M}}: 
\mathrm{ev}_1(v) \in W^u_h(c_1)\},\\
\tilde{\mathcal{M}}^+(c_2)&:=&\{v \in \tilde{\mathcal{M}}:
\mathrm{ev}_2(v) \in W^s_h(c_2)\},\\
\tilde{\mathcal{M}}(c_1,c_2)&:=&
\tilde{\mathcal{M}}^-(c_1) \cap \tilde{\mathcal{M}}^+(c_2).
\end{eqnarray*}
Immediately from Theorem~\ref{mani} the following Corollary follows.
\begin{cor}
For generic Riemannian metric $g$ on $M$ the spaces
$\tilde{\mathcal{M}}_u(c_1)$, $\tilde{\mathcal{M}}_s(c_2)$, and
$\tilde{\mathcal{M}}(c_1,c_2)$ are finite dimensional manifolds. Its
local dimensions are
\begin{eqnarray*}
\dim_v \tilde{\mathcal{M}}^-(c_1)&=&\mathrm{Ind}(c_1)
-\mathrm{ind}_f(\mathrm{ev}_2(v))+m-1\\
\dim_v \tilde{\mathcal{M}}^+(c_2)&=&\mathrm{ind}_f(\mathrm{ev}_1(v))
+\dim_{\mathrm{ev}_1(v)}\mathrm{crit}(f)-\mathrm{Ind}(c_2)+m-1\\
\dim_v \tilde{\mathcal{M}}(c_1,c_2)&=&\mathrm{Ind}(c_1)-\mathrm{Ind}(c_2)
+m-1.
\end{eqnarray*}
\end{cor}
It follows from the Corollary above that for generic Riemannian
metric $g$ the moduli space of flow lines
with cascades $\mathcal{M}(c_1,c_2)$ is 
a manifold of dimension $\mathrm{Ind}(c_1)-\mathrm{ind}(c_2)-1$ in
a neighbourhood of an element
$v=((x_j)_{1 \leq j \leq m},(t_j)_{1 \leq j \leq m-1})$ all of
whose $t_j>0$. It remains to consider the case where some of
the $t_j$'s are zero. In a neighbourhood of such an element the
number of cascades may vary and we need a gluing theorem to
parametrise such a neighbourhood. 

In view of the transversality, the subspace
$\mathcal{B}(c_1,c_2)$ of $\mathcal{B}$ for $c_1,c_2 \in \mathrm{crit}(h)$
defined by
$$\mathcal{B}(c_1,c_2):=\{v \in \mathcal{B}: 
\mathrm{ev}_1(v) \in W^u_h(c_1),\,\,\mathrm{ev}_2(v) \in W^s_h(c_2)\}$$
is actually a submanifold. 
Recall that the evaluation map
$\mathrm{EV}_n \colon \tilde{\mathcal{M}}^n \to \mathrm{crit}(f)^{2n}$
was defined by $\mathrm{EV}_n(v_1,\cdots,v_n) \mapsto
(\mathrm{ev}_1(v_j)_{1 \leq j \leq n}, \mathrm{ev}_2(v_j)_{1 \leq j \leq n})$
and $A_n(c_1,c_2) \subset \mathrm{crit}(f)^{2n}$ was defined as
the subspace of all tuples 
$((p_j)_{1 \leq j \leq n},(q_j)_{1 \leq j \leq n}) \in \mathrm{crit}(f)^{2n}$
satisfying $p_1 \in W^u_h(c_1), q_n \in W^s_h(c_2)$, and 
$p_{j+1}=p_j$ for $1 \leq j \leq n-1$. 
We next define for $T$ large enough the \textbf{pregluing map}
$$\#^0:\mathrm{EV}_n^{-1}(A_n(c_1,c_2)) \times (T,\infty)^{n-1}
\to \mathcal{B}(c_1,c_2).$$

\medskip
\begin{center}\input{f3.pstex_t}\end{center}
\medskip

Let 
$\textbf{v}=(v_k)_{1 \leq k \leq n}=
((x_{k,j})_{1 \leq j \leq m_k},(
t_{k,j})_{1 \leq j \leq m_k-1}))_{1 \leq k \leq n} 
\in \mathrm{EV}_n^{-1}(A_n(c_1,c_2))$ and
$\textbf{R}=(R_j)_{1 \leq j \leq n-1} \in (T,\infty)^{n-1}$. 
For $1 \leq p \leq n$ and
$$\sum_{k=1}^{p-1} m_k-p+2 \leq i < \sum_{k=1}^p m_k-p+1$$
set
$$s_i:=t_{p,i-\sum_{k=1}^{p-1}m_k+p-1}.$$
For 
$$\sum_{k=1}^{p-1} m_k-p+2 < i < \sum_{k=1}^p m_k-p+1$$
put
$$y_i:=x_{p,i-\sum_{k=1}^{p-1}m_k+p-1}.$$
Define in addition
$$y_1:=x_{1,1}, \quad y_{\sum_{k=1}^n m_k-n+1}:=x_{n,m_n}.$$
Recall that $x_{p,m_p}(s)$ converges as $s$ goes to $\infty$ to 
$\mathrm{ev}_2(v_p)$ and  $x_{p+1,1}(s)$ converges as $s$ goes
to $-\infty$ to $\mathrm{ev}_2(v_p)$ as well. In particular, 
it follows that there
exist $\xi_p(s), \eta_p(s) \in T_{\mathrm{ev}_2(v_p)} M$ such that
$x_{p,m_p}(s)=\exp_{\mathrm{ev}_2(v_p)}(\xi_p(s))$ for large negative $s$
and $x_{p+1,1}(s)=\exp_{\mathrm{ev}_2(v_p)}(\eta_p(s))$ for large positive
$s$ 
For $1 \leq p \leq n-1$ and 
$$i=\sum_{k=1}^p m_k-p+1$$
put
$$y_i:=x_{p,m_p} \#^0_{R_p} x_{p+1,1}$$
where
$$x_{p,m_p} \#^0_{R_p}x_{p+1,1}:=$$
$$\left\{ \begin{array}{cc}
x_{p,m_p}(s+R_p), & s \leq -R_p/2-1,\\
\exp_{\mathrm{ev}_2(v_p)}(\beta(-s-R_p/2)\eta_p(s+R_p)), 
& -R_p/2-1 \leq s \leq -R_p/2,\\
\mathrm{ev}_2(v_p), & -R_p/2 \leq s \leq R_p/2,\\
\exp_{\mathrm{ev}_2(v_p)}(\beta(s-R_p/2)\xi_p(s-R_p)), 
& R_p/2 \leq s \leq R_p/2+1,\\
x_{p+1,1}(s-R_p), & s \geq R_p/2+1.
\end{array}\right.$$
Here $\beta: \mathbb{R} \to [0,1]$ is a cutoff function equal to $1$ for
$s \geq 1$ and equal to $0$ for $s \leq 0$.
We abbreviate for the number of cascades of the image of the
pregluing map $\#^0$ 
$$\ell=\sum_{k=1}^n m_k-n+1.$$
Then we define
$$v^0_{\textbf{R}}:=\#^0(\textbf{v},\textbf{R}):=((y_i)_{1 \leq i \leq \ell},
(s_i)_{1 \leq i \leq \ell-1}).$$
$v^0_{\textbf{R}}$ will in general only lie ``close'' to  
$\tilde{\mathcal{M}}$. We will next construct $v_{\textbf{R}} \in
\tilde{\mathcal{M}}$ in a small neighbourhood of $v_{\textbf{R}}^0$.
It can be shown that
there exists a Riemannian metric $g$ on $M$ such that $\mathrm{crit}(f)$,
$W^u_h(c_1)$, and $W^s_h(c_2)$ are
totally geodesic with respect to $g$. For 
$\zeta \in T_v \mathcal{B}(c_1,c_2)$ let 
$$\rho(v,\zeta):\mathcal{E}_v
\to \mathcal{E}_{\exp_v(\zeta)}
\footnote{We denote the restriction of $\mathcal{E}$ to $\mathcal{B}(c_1,c_2)$
also by $\mathcal{E}$.}$$
be the parallel transport along the path $\tau \mapsto \exp_v(\tau \zeta)$
for $\tau \in [0,1]$. 
Define
$$F_v:T_v \mathcal{B}(c_1,c_2) \to \mathcal{E}_v, \quad 
\zeta \mapsto \rho(v,\zeta)^{-1}\mathcal{F}(\exp_v(\zeta)).$$
Note that
$$DF_v(0)=D\mathcal{F}(v).$$
For large $R \in \mathbb{R}$ 
and $i \in \mathbb{N}_0$ define $\gamma^i_{\delta,R}:\mathbb{R} \to \mathbb{R}$
by
$$\gamma^i_{\delta,R}(s):=\left\{\begin{array}{cc}
(1-\beta[-s])\gamma_\delta(s+R)+\beta(s)\gamma_\delta(s-R) & i \geq 1\\ 
(1-\beta[-s-R])\gamma_\delta(s+R)+\beta(s+R)\gamma_\delta(s-R) & i=0
\end{array}\right.$$
Here $\beta$ is a smooth cutoff function which equals $0$ for $s\leq -1$
and equals $1$ for $s \geq 0$.
For a locally integrable real valued function
$f\colon \mathbb{R} \to \mathbb{R}$ which has
weak derivatives up to order $k$ we define the norms
$$||f||_{k,p,\delta,R}:=
\left\{\begin{array}{cc}||\gamma^1_{\delta,R}f||_p & k=0\\
\sum_{i=0}^k||\gamma^i_{\delta,R}\partial^i f||_p+||f(0)|| & k>0.
\end{array}\right.$$
The $||\,||_{k,p,\delta, R}$-norm is equivalent to the
$||\,||_{k,p,\delta}$-norm, but their ratio diverges as $R$ goes to 
infinity. For $1 \leq i \leq \ell$ set
$$R^i:=\left\{ \begin{array}{cc}
R_q & i=\sum_{k=1}^q m_k-q+1 \\
0 & \mathrm{else}\end{array}\right.$$
We modify the Banachspace norm (\ref{norm}) on 
$T_{v^0_{\textbf{R}}}\mathcal{B}(c_1,c_2)$ in the following way. For
$\zeta \in T_{v^0_{\textbf{R}}}\mathcal{B}(c_1,c_2)$ we define
$$||\zeta||_{\textbf{R}}:=\sum_{j=1}^\ell(
||\xi_{j,0}||_{1,p,\delta,R^j}+||\xi_{j,1}||+||\xi_{j,2}||)+\sum_{j=1}^{\ell-1}
|\delta_j|.$$
Analoguously, we define for $\eta \in \mathcal{E}_{v^0_{\textbf{R}}}$ the
norm
$$||\eta||_{\textbf{R}}:=\sum_{j=1}^\ell||\eta_j||_{p,\delta,R^j}.$$ 
These norms were introduced in
\cite[Section 18]{fukaya-oh-ohta-ono} and are required to guarantee the
uniformity of the constant $c$ in (\ref{rightinv}) below. 
Abbreviate
$$D_{\textbf{R}}:=DF_{v^0_{\textbf{R}}}(0)=D_{v^0_{\textbf{R}}}.$$
Recall that $R_j>T$ for
$1 \leq j \leq n-1$. 
It is shown in \cite[Chapter 2.5]{schwarz1} that 
if $T$ is large enough then there exists a constant $c>0$ independent 
of $\textbf{R}$ and a right inverse $Q_{\textbf{R}}$ of 
$D_{\textbf{R}}$, i.e.
\begin{equation}\label{Q}
D_{\textbf{R}} \circ Q_{\textbf{R}}=\mathrm{id},
\end{equation}
such that
\begin{equation}\label{orth}
\mathrm{im} Q_{\textbf{R}}=\mathrm{ker}(D_{\textbf{R}})^\perp
\end{equation}
and 
\begin{equation}\label{rightinv}
||Q_{\textbf{R}} \eta||_{\textbf{R}} \leq c||\eta||_{\textbf{R}},
\end{equation}
for every $\eta \in \mathcal{E}_{v^0_{\textbf{R}}}$. Moreover, it follows
from the construction of the pregluing map that there exist constants
$c_0>0$ and $\kappa>0$ such that
\begin{equation}\label{est}
||F_{v^0_{\textbf{R}}}(0)||_{\textbf{R}} 
\leq c_0 e^{-\kappa ||\textbf{R}||}.
\end{equation}
Now (\ref{Q}), (\ref{orth}), (\ref{rightinv}), 
(\ref{est}) together with the Banach Fixpoint 
theorem imply that for $||\textbf{R}||$ large enough there exists
a constant $c>0$ and a
unique $\xi_{\textbf{R}}:=
\xi_{\textbf{v},\textbf{R}} \in \mathrm{ker}(D_{\textbf{R}})^\perp$
such that
\begin{equation}\label{glue}
F_{v^0_{\textbf{R}}}(\xi_{\textbf{R}})=0, \quad
||\xi_{\textbf{R}}||_{\textbf{R}} \leq c||F_{v^0_{\textbf{R}}}(0)||. 
\end{equation} 
For details see \cite[Chapter 2.5]{schwarz1}. Let
$U(\textbf{v})$ be a small neighbourhood of $\textbf{v}$ in 
$\mathrm{EV}^{-1}_n(A_n(c_1,c_2))$.
Define the gluing map
$$\#:U(\textbf{v}) \times (T,\infty)^{n-1}
\to \tilde{\mathcal{M}}(c_1,c_2)$$
by
$$\#(\textbf{w},\textbf{R}):=
\exp_{w^0_{\textbf{R}}}(\xi_{\textbf{w},\textbf{R}}).$$
Let 
$$N:=\sum_{k=1}^n m_k.$$
Then $\mathbb{R}^N$ acts on $U(\textbf{v})$ 
by timeshift on each cascade of the first factor. The gluing
map induces a map
$$\hat{\#}:(U(\textbf{v})/\mathbb{R}^N) \times (T,\infty)^{n-1} \to
\mathcal{M}(c_1,c_2)$$
which is an embedding for $T$ large enough.
\\ \\
\textbf{Proof of Theorem~\ref{manifold}:} For $m \in \mathbb{N}$ let
$\mathbb{N}_m:=\{1,\ldots,m\}$. For 
$I \subset \mathbb{N}_{m-1}$ let 
$$\mathcal{M}_{m,I}(c_1,c_2) \subset \mathcal{M}_m(c_1,c_2)$$
be the set of flow lines with cascades
$((x_k)_{1 \leq k \leq m},(t_k)_{1 \leq k \leq m-1})$
in $\mathcal{M}_m(c_1,c_2)$ such that
$$t_k>0 \quad \mathrm{if}\,\, k \in I, \qquad t_k=0 \quad \mathrm{if} 
\,\, k \notin I.$$
It follows from Theorem~\ref{mani} that for generic
Riemannian metric $g$ on $M$ the space $\mathcal{M}_{m,I}(c_1,c_2)$
is a smooth manifold. Note that
$$\mathcal{M}_m(c_1,c_2)=\bigcup_{I \subset \mathbb{N}_{m-1}}
\mathcal{M}_{m,I}(c_1,c_2)$$
and 
$$\mathrm{int}\mathcal{M}_m(c_1,c_2)=\mathcal{M}_{m,\mathbb{N}_{m-1}}
(c_1,c_2).$$
It follows from (\ref{dim}) in Theorem~\ref{mani} that
$$\dim\mathcal{M}_{m,\mathbb{N}_{m-1}}(c_1,c_2)=\mathrm{Ind}(c_1)
-\mathrm{Ind}(c_2)-1.$$
In particular, $\mathcal{M}_m(c_1,c_2)$ has for generic $g$ the structure
of a manifold with corners of dimension 
$\mathrm{Ind}(c_1)-\mathrm{Ind}(c_2)-1$. We put for $n \in \mathbb{N}$
$$\mathcal{M}_{\leq n}:=\bigcup_{1 \leq m \leq n} \mathcal{M}_m(c_1,c_2).$$
We show by induction on $n$ that for generic $g$ the set 
$\mathcal{M}_{\leq n}(c_1,c_2)$ can be endowed with the structure of a 
manifold of dimension $\mathrm{Ind}(c_1)-\mathrm{Ind}(c_2)-1$. This is clear
for $n=1$. It follows from gluing that $\mathcal{M}_{\leq n}(c_1,c_2)$ can
be compactified to a manifold with corners
$\bar{\mathcal{M}}_{\leq n}(c_1,c_2)$ such that
$$\partial \bar{\mathcal{M}}_{\leq n}(c_1,c_2)
=\bigcup_{I \subsetneq \mathbb{N}_n} \mathcal{M}_{n+1,I}(c_1,c_2)
=\partial \mathcal{M}_{n+1}(c_1,c_2).$$
Hence 
$$\mathcal{M}_{\leq n+1}(c_1,c_2)=\mathcal{M}_{\leq n}(c_1,c_2)
\cup \mathcal{M}_{n+1}(c_1,c_2)$$
can be endowed with the structure of a manifold such that
$$\dim \mathcal{M}_{\leq n+1}(c_1,c_2)=\dim\mathcal{M}_{\leq n}(c_1,c_2)
=\mathrm{Ind}(c_1)-\mathrm{Ind}(c_2)-1.$$
This proves the theorem. \hfill $\square$

\subsection[Morse-Bott homology]{Morse-Bott homology}

We assume in this subsection that $M$ is a compact manifold. 
\begin{fed}
We say that a quadruple $(f,h,g,g_0)$ consisting of
a Morse-Bott function $f$ on $M$ a Morse function
$h$ on $\mathrm{crit}(f)$, a Riemannian metric $g$ on
$M$ and a Riemannian metric $g_0$ on $\mathrm{crit}(f)$
is a 
\emph{\bf{regular Morse-Bott quadruple}} if the following
conditions hold.
\begin{description}
 \item[(i)] $h$ and $g_0$ satisfy the Morse-Smale condition, i.e.
  stable and unstable manifolds of the gradient of $h$ with
  respect to $g_0$ on $\mathrm{crit}(f)$ intersect transversally.
  \item[(ii)] $g$ is $(f,h,g_0)$-regular, in the sense of 
  Definition~\ref{fhg0reg}.
\end{description}
\end{fed}
Since the Morse-Smale condition is generic, see \cite[Chapter 2.3]{schwarz1},
it follows from Theorem~\ref{mani}, that regular Morse-Bott quadruples
exist in abundance. In particular, every pair $(f,g)$ consisting
of a Morse function $f$ on $M$ and a Riemannian metric $g$ on $M$ 
which satisfy the Morse-Smale condition gives a regular Morse-Bott
quadruple.

For a pair $(f,h)$ consisting of a Morse-Bott 
function $f$ on $M$ and a Morse-function $h$ on $\mathrm{crit}(f)$, 
we define the chain complex
$CM_*(M;f,h)$ as the $\mathbb{Z}_2$ vector space generated 
by the critical points of $h$ which is graded by the index. More precisely,
$CM_k(M;f,h)$ are formal sums of the form
$$\xi=\sum_{\substack{c \in \mathrm{crit}(h)\\
\mathrm{ind}(c)=k}}\xi_c c$$
with $\xi_c \in \mathbb{Z}_2$. 
For
generic pairs $(g,g_0)$ of a Riemannian metric $g$ on $M$ and a Riemannian
metric $g_0$ on $\mathrm{crit}(f)$, the moduli-spaces of of flow lines
with cascades $\mathcal{M}(c_1,c_2)$ is a smooth manifold of dimension
$$\dim{\mathcal{M}}(c_1,c_2)=\mathrm{ind}(c_1)-\mathrm{ind}(c_2)-1.$$
If $\dim{\mathcal{M}}(c_1,c_2)$ equals $0$, then
$\mathcal{M}(c_1,c_2)$ is compact by Theorem~\ref{fcc}. Set
$$n(c_1,c_2):=\#\mathcal{M}(c_1,c_2) \,\,\mathrm{mod}\,2.$$
We define a boundary operator
$$\partial_k: CM_k(M;f,h) \to CM_{k-1}(M;f,h)$$
by linear extension of 
$$\partial_k c=\sum_{\mathrm{ind}(c')=k-1}n(c,c')c'$$
for $c \in \mathrm{crit}(h)$ with $\mathrm{ind}(c)=k$.
The usual gluing and compactness arguments imply that
$$\partial^2=0.$$
This gives rise to homology groups
$$HM_k(M;f,h,g,g_0):=\frac{\mathrm{ker}\partial_{k+1}}
{\mathrm{im}\partial_k}.$$

\begin{thm}\label{continuation}
Let $(f^\alpha,h^\alpha,g^\alpha,g_0^\alpha)$ and 
$(f^\beta,h^\beta,g^\beta,g_0^\beta)$ be two regular 
quadrupels. Then the homologies
$HM_*(M;f^\alpha,h^\alpha,g^\alpha,g_0^\alpha)$ and
$HM_*(M;f^\beta,h^\beta,g^\beta,g_0^\beta)$ are naturally isomorphic. 
\end{thm}
\textbf{Proof:} Pick some $\ell \in \mathbb{N}$ and choose for
$1 \leq k \leq \ell$ smooth functions 
$f_k \in C^\infty(\mathbb{R} \times M, \mathbb{R}$ and
smooth families of Riemannian metrics $g_{k,s}$ on $TM$ with
$f_k(s, \cdot)$ and $g_{k,s}$ independent of $s$ for $|s| \geq T$ for some
large enough constant $T>0$  such that
$$f_1(-T)=f^\alpha, \quad f_\ell(T)=f^\beta, \quad
f_k(T)=f_{k+1}(-T),\,\, 1 \leq k \leq \ell-1$$
and
$$g_{1,-T}=g^\alpha, \quad g_{\ell,T}=g^\beta, \quad
g_{k,T}=g_{k+1,-T},\,\, 1 \leq k \leq \ell-1.$$
We assume further that $f_k(T)$ is Morse-Bott for $1 \leq k \leq \ell-1$. 
For $2 \leq k \leq \ell$ let $r_k \in \mathbb{R}_{\geq}$ be nonnegative
real numbers. Choose
smooth Morse functions $h_1 \in C^\infty((-\infty,0]
\times \mathrm{crit}(f^\alpha),\mathbb{R})$, 
$h_{\ell+1} \in C^\infty([0,\infty) \times \mathrm{crit}(f^\beta),\mathbb{R})$,
and $h_k \in C^\infty([0,r_k] \times \mathrm{crit}(f_{k+1}(T))$
and smooth families of Riemannian metrics $g_{0,1,s}$ on
$\mathrm{crit}(f^\alpha)$ for $s \in (-\infty,0]$ and 
$g_{0,\ell+1,s}$ on $\mathrm{crit}(f^\beta)$ for $s \in [0,\infty)$ 
and $g_{0,k,s}$ on $\mathrm{crit}(f_{k+1}(T))$ for $s \in [0,r_k].$
They are required to fulfill
$$h_1(s)=h^\alpha,\,\,g_{0,1,s}=g_0^\alpha,\,\,s \leq -T,$$
$$h_{\ell+1}(s)=h^\beta,\,\, g_{0,\ell+1,s}=g_0^\beta,\,\, s \geq T.$$
For $c_1 \in \mathrm{crit}(h^\alpha)$, $c_2 \in \mathrm{crit}(h^\beta)$,
$m_1,m_2 \in \mathbb{N}_0$ we consider the following 
flow lines from $c_1$ to $c_2$ with $m=m_1+m_2+\ell$ cascades 
$$(\textbf{x},\textbf{T})=((x_k)_{1 \leq k \leq m},(t_k)_{1 \leq k \leq m-1})$$
for $x_k \in C^\infty(\mathbb{R},M)$ and 
$t_k \in \mathbb{R}_{\geq}:=\{r \in \mathbb{R}: r \geq 0\}$ which satisfy
the following conditions:
\begin{description}
 \item[(i)] $x_k$ are nonconstant solutions of 
  $$\dot{x}_k(s)=-\nabla_{\tilde{g}_{k,s}} \tilde{f}_k(s,x_k),$$ 
  where 
  $$\tilde{f}_k=\left\{\begin{array}{cc}
  f^\alpha & 1 \leq k \leq m_1\\
  f_{k-m_1} & m_1+1 \leq k \leq m_1+\ell\\
  f^\beta & m_1+\ell+1 \leq k \leq m
  \end{array}\right.$$
  and
  $$\tilde{g}_k=\left\{\begin{array}{cc}
  g^\alpha & 1 \leq k \leq m_1\\
  g_{k-m_1} & m_1+1 \leq k \leq m_1+1+\ell\\
  g^\beta & m_1+\ell+1 \leq k \leq m.
  \end{array}\right.$$
 \item[(ii)] There exists $p_1 \in W^u_{h^\alpha}(c_1)$ and 
  $p_2 \in W^s_{h^\beta}(c_2)$
  such that $\lim_{s \to -\infty}x_1(s)=p_1$ and 
  $\lim_{s \to \infty}x_m(s)=p_2$. 
 \item[(iii)] denote
  $$\tilde{h}_k=\left\{\begin{array}{cc}
  h^\alpha & 1 \leq k \leq m_1-1\\
  h_{k-m_1+1} & m_1 \leq k \leq m_1+\ell-1\\
  f^\beta & m_1+\ell \leq k \leq m-1
  \end{array}\right.$$
  and 
  $$\tilde{g}_{0,k}=\left\{\begin{array}{cc}
  g^\alpha_0 & 1 \leq k \leq m_1-1\\
  g_{0,k-m_1+1} & m_1 \leq k \leq m_1+1+\ell-1\\
  g^\beta_0 & m_1+\ell \leq k \leq m-1.
  \end{array}\right.$$
  For $1 \leq k \leq m-1$ there are Morse-flow lines $y_k$
  of $h$, i.e. solutions of 
  $$\dot{y}_k(s)=-\nabla_{\tilde{h}_{0,k,s}}\tilde{h}_k(s,y_k),$$
  such that
  $$\lim_{s \to \infty}x_k(s)=y(0)$$
  and
  $$\lim_{s \to -\infty}x_{k+1}(s)=y(t_k).$$
 \item[(iv)] $t_{k+m_1-1}=r_k$ for $2 \leq k \leq \ell$.
\end{description}
For generic choice of the data the space of solutions of (i) to (iv)
is a smooth manifold whose dimension is given by the difference of
the indices of $c_1$ and $c_2$. If $I(c_1)=I(c_2)$ then this manifold
is compact and we denote by $n(c_1,c_2) \in \mathbb{Z}_2$ its cardinality modulo 2.
We define a map
$$\Phi^{\alpha \beta}: CM_*(M;f^\alpha,h^\alpha) \to
CM_*(M;f^\beta,h^\beta)$$
by linear extension of
$$\Phi^{\alpha \beta} c=\sum_{\substack{c' \in \mathrm{crit}(h^\beta)\\
\mathrm{ind}(c')=\mathrm{ind}(c)}} n(c,c') c'$$
where $c \in \mathrm{crit}(f^\alpha)$. Standard arguments, see
\cite[Chapter 4.1.3]{schwarz1}, show that
$\Phi^{\alpha\beta}$ induces isomorphisms on homologies
$$\phi^{\alpha \beta}: HM_*(M;f^\alpha,h^\alpha,g^\alpha,g_0^\alpha)
\to HM_*(M;f^\beta,h^\beta,g^\beta,g_0^\beta)$$
which satisfy
$$\phi^{\alpha \beta} \circ \phi^{\beta \gamma}=\phi^{\alpha \gamma}, 
\quad \phi^{\alpha \alpha}=\mathrm{id}.$$
This proves the theorem. \hfill $\square$
\\ \\
Theorem~\ref{continuation} allows us to define the
\textbf{Morse-Bott homology} of $M$ by
$$HM_*(M):=HM_*(M;f,h,g,g_0)$$
for some regular quadruple $(f,h,g,g_0)$. Either for the
special case of a Morse function or the case where
$f$ vanishes identically one obtains that Morse-Bott homology is isomorphic to
Morse homology. Hence we have proved the following Corollary.
\begin{cor}
Morse-Bott homology of a compact manifold $M$ is isomorphic to the Morse-homology
of $M$ and hence also to the singular homology of $M$.
\end{cor}

(U. Frauenfelder) Departement of Mathematics, Hokkaido University,
Sapporo 060-0810, Japan

E-mail address: urs@math.sci.hokudai.ac.jp

\end{document}

%% file: f1.pstex_t
\begin{picture}(0,0)%
\includegraphics{f1.pstex}%
\end{picture}%
\setlength{\unitlength}{3947sp}%
\begingroup\makeatletter\ifx\SetFigFont\undefined%
\gdef\SetFigFont#1#2#3#4#5{%
  \reset@font\fontsize{#1}{#2pt}%
  \fontfamily{#3}\fontseries{#4}\fontshape{#5}%
  \selectfont}%
\fi\endgroup%
\begin{picture}(3624,3312)(2839,-4636)
\put(4726,-3211){\makebox(0,0)[lb]{\smash{\SetFigFont{12}{14.4}{\familydefault}{\mddefault}{\updefault}{\color[rgb]{0,0,0}$\dot{y_1}=-\nabla h(y_1)$}%
}}}
\put(3001,-2311){\makebox(0,0)[lb]{\smash{\SetFigFont{12}{14.4}{\rmdefault}{\mddefault}{\updefault}{\color[rgb]{0,0,0}$\dot{x_1}=-\nabla f(x_1)$}%
}}}
\end{picture}

%% file: f2.pstex_t
\begin{picture}(0,0)%
\includegraphics{f2.pstex}%
\end{picture}%
\setlength{\unitlength}{3947sp}%
\begingroup\makeatletter\ifx\SetFigFont\undefined%
\gdef\SetFigFont#1#2#3#4#5{%
  \reset@font\fontsize{#1}{#2pt}%
  \fontfamily{#3}\fontseries{#4}\fontshape{#5}%
  \selectfont}%
\fi\endgroup%
\begin{picture}(3624,3237)(2839,-4561)
\end{picture}

%% file: f3.pstex_t
\begin{picture}(0,0)%
\includegraphics{f3.pstex}%
\end{picture}%
\setlength{\unitlength}{3947sp}%
\begingroup\makeatletter\ifx\SetFigFont\undefined%
\gdef\SetFigFont#1#2#3#4#5{%
  \reset@font\fontsize{#1}{#2pt}%
  \fontfamily{#3}\fontseries{#4}\fontshape{#5}%
  \selectfont}%
\fi\endgroup%
\begin{picture}(4374,2691)(1189,-2515)
\end{picture}